\documentclass{article}
\usepackage[square,sort,comma,numbers]{natbib}
\usepackage{fullpage}
\usepackage[utf8]{inputenc}
\usepackage[english]{babel}
\usepackage{amsmath}
\usepackage{amssymb}
\usepackage{amsthm}
\usepackage{hyperref}
\usepackage{graphicx}
\usepackage{subfig}
\usepackage{epsfig}

\newtheorem{theorem}{Theorem}[section]

\newtheorem{corollary}{Corollary}[section]
\newtheorem{lemma}[theorem]{Lemma}

\begin{document}
\title{ \bf{The polar-generalized normal distribution}}
\author{\vspace {.2cm} {\normalsize Masoud  Faridi,
Majid Jafari Khaledi$^a$}\\
{\scriptsize masoudfaridi@modares.ac.ir}\\
$^a$ {\scriptsize  Department of Statistics, Tarbiat Modares University, P.O. Box 14115-134, Tehran, Iran}\\
{\scriptsize jafari-m@modares.ac.ir}
}
\date{}
\maketitle
\begin{abstract}
This paper introduces  an extension to the normal distribution  through the polar method to capture bimodality and asymmetry, which are often observed characteristics  of empirical data. The later two features are entirely controlled by a separate scalar parameter. Explicit expressions for the cumulative distribution function, the density function and the moments were derived. The stochastic representation of the distribution facilitates implementing Bayesian estimation via the Markov chain Monte Carlo methods. Some real-life  data as well as simulated data are analyzed  to illustrate the flexibility of the distribution for modeling  asymmetric bimodality.
\end{abstract}

 {\it Keywords and phrases:} Skewness; Bimodality; Polar method; Beta distribution; Bayesian estimation.
\section{Introduction}

The normal distribution is often used  to model data with symmetric distributions.  However, the distribution of data were collected from a wide range of applications reveals asymmetry. Examples arise, for instance, in wind speed data (Zhu and Genton, 2012), temperature data (North et al., 2011), soil and air pollution data (Jafari Khaledi and Rivaz, 2009; Zareifard and Jafari Khaledi; 2013), precipitation data (Xu and Genton, 2017) and diabetes data (Kristensen et al., 2010).
The skew-normal (SN) distribution, introduced by Azzalini (1985), has been used successfully to model data sets that have asymmetric behavior. The normal distribution is a special case of this distribution. Despite attractive features of this distribution, it is not suitable when the distribution of the data are biomodal. Since many data sets from applied fields are bimodal (for example gene expression patterns in breast cancer (Wang et al., 2009; Ertel, 2010), material characterization in engineering (Dierickx et al.,  2000), water vapor and temperature in meteorology (Zhang et al. , 2003; Dahdouh and Jafari Khaledi, 2020),  there is extensive research on modeling bimodality. The mixture of normal distributions is often used in many different  areas to model data with biomodal distributions. This strategy is often criticized for its non-identifiability problem (McLachlan and  Peel, 2000; Marin, Mengersen and  Robert, 2005). Gomez et al. (2011) introduced a bimodal extension of the skew-normal distribution with application to a pollen dissemination data set, but it suffers complexity and estimation problems inherited from the skewed normal model.

The main objective of this paper is to provide a parsimonious probability distribution that could be used as an alternative to the mixture models and can  account for both asymmetry and bimodality. The basic idea relies on extending the standard polar method.
In this method, one generate a random pair $(X,Y)=\left(R \cos(2\pi U) , R \sin (2\pi U)\right)$ where
$R$ and $U$ are random polar coordinates. The random variables $R$ and $U$ are independent.  The random variable $U$ is uniformly distributed on the unit interval  $(0,1)$,
 and $R$ has a given distribution that is easy to sample from.  In the context of one-liners, we thus have

 \begin{eqnarray}
 \label{polarx}
    X &=& R \cos( 2\pi U).
 \end{eqnarray}
 In this case, if $R \sim \sqrt{\chi^2_{2}}$, $\chi^2_{2}$ has a chi-squared distribution with $2$ degrees of freedom,
   then the distribution of $X$ is normal. The well-known Box-Muller formula for normal random variables, $X = \sqrt{-2 \log U_2} \cos( 2\pi U_1)$
 is thus a simple one-linear, where $U_1$ and $U_2$ are independent random variables from the uniform distribution on the unit interval $(0, 1)$.

In our construction, instead of the uniform distribution, we propose to consider the beta distribution for the random variable $U$ in (\ref{polarx}). We call the resulting
construction the polar-generalized normal
(PGN) distribution. We refer to Alexander et al. (2012),  Mameli and  Musio (2013 and 2016), Alleva and  Giommi (2016) and Eugene, Lee and  Famoye (2002) for information on other generalizations of normal distribution that present bimodality and involve the beta distribution. The PGN distribution can model asymmetry and bimodality and include as special cases the normal distribution. In particular, it contains easily interpretable and estimable parameters. Our approach allows for testing normality and symmetry within the $PGN$ family via Beysian framework.

The remainder of this paper is organized as follows: Section 2 is devoted to  the extension of the polar method. Moreover, basic properties and other characteristics of the PGN distribution are presented.  In particular, some expressions are derived which simplify moments derivation.
In Section 3, the stochastic representation for a random variable following the PGN distribution is used  to implement Bayesian estimation. The results of a simulation study reported in Section 4 reveal satisfactory behaviour of the Bayesian estimates.  In Section 5, we provide some empirical applications. Finally, this paper concludes with a brief review of the main results.

\section{The PGN model  }
Below we consider an  extension of the polar method given in equation (\ref{polarx}).
Let $V$, $U$  be independent random variables.
Assume that $V$ has a chi-squared distribution with  $k$
 degrees of freedom, that is, the density is given by
 $ f_V(v) =\frac{1}{2^{\frac{k}{2}}\Gamma(\frac{k}{2})} v^{\frac{k}{2}-1} e^{-\frac{v}{2}}$ for $0<v< +\infty$ where $\Gamma(\cdot )$ is the gamma function and $k>0$. Suppose $U$ has a beta distribution with parameters $2\mu$ and $2(1-\mu)$, that is, the density is given by
 $f_U(u)=\frac{1}{B(2\mu,2(1-\mu))} u^{2\mu-1} (1-u)^{2(1-\mu)-1}$ for $0<u<1$ where $0<\mu<1$ and
 $B(\alpha ,\beta) =\frac{\Gamma(\alpha)\Gamma(\beta)}{\Gamma(\alpha +\beta)}$ is the beta function. We define the polar-generalized normal distribution as the distribution of the product
\begin{equation}
\label{definitionpolar}
X=\sqrt{V} \cos(\pi U).
\end{equation}
 We will use the short notation $PGN(\mu , k)$ distribution.
  The standard normal distribution is a special case of the $PGN(\mu , k)$  distribution when $\mu$ is $0.5$ and $k$  is $2$.
In the following subsections, we study distributional properties and
present some expansions for
moments of  $PGN(\mu , k)$.

\subsection{The distributional properties}

\theorem{
The cumulative distribution function (cdf) of $PGN(\mu , k)$ is given by
\begin{eqnarray}
\label{eqcdf}
F_{X}(x)&=&
\left\{
\begin{array}{cc}
\frac{I_{0.5}(2(1-\mu) , 2\mu )}{B(2\mu , 2(1-\mu))}
-\frac{1}{\Gamma(\frac{k}{2})}
\int_{0.5}^{1} \gamma(\frac{k}{2} , \frac{x^2}{2\cos^2(\pi u)} ) f_U(u) du
 & \text{if  } x<0, \\
\frac{I_{0.5}(2(1-\mu) , 2\mu )}{B(2\mu , 2(1-\mu))}+\frac{1}{\Gamma(\frac{k}{2})}\int_{0}^{0.5} \gamma(\frac{k}{2} , \frac{x^2}{2\cos^2(\pi u)} )f_U(u) du
&  \text{if  }  x\geq 0,
\end{array}
\right.
\end{eqnarray}
where $I_{c}(\alpha , \beta )$ is the regularized incomplete beta function,
$\gamma(s ,t)$ is the
regularized incomplete gamma function and $f_U(u)$ is the density of
beta distribution with parameters $2\mu$ and $2(1-\mu)$. Specifically,  $I_{c}(\alpha , \beta )=\int_0^c u^{\alpha -1}(1-u)^{\beta -1} du$ and
$\gamma(s ,t)=\int_0^{t} y^{s-1} e^{-y} dy$.
}
\proof{

\begin{eqnarray}
F_{X}(x)&=&P(X\leq x)=P(\sqrt{V} \cos(\pi U)\leq x)=
\int_{0}^{1} P(\sqrt{V} \cos(\pi u)\leq x)f_U(u) du
 \nonumber \\
 \label{eqnFXproof}
  &=&
\int_{0}^{0.5} P(\sqrt{V} \cos(\pi u)\leq x)f_U(u) du +
\int_{0.5}^{1} P(\sqrt{V} \cos(\pi u)\leq x)f_U(u) du.
\end{eqnarray}

The first term of the Equation (\ref{eqnFXproof}) is equal to

\begin{eqnarray}
\int_{0}^{0.5} P(\sqrt{V} \cos(\pi u)\leq x)f_U(u) du &=&
\left\{
\begin{array}{cc}
  0 & \text{if   } x<0 \\
  \int_{0}^{0.5} P(\sqrt{V} \leq \frac{x}{\cos(\pi u)} )f_U(u) du  & \text{if   } x\geq 0
\end{array}
\right.
\nonumber \\ &=&
\left\{
\begin{array}{cc}
  0 & \text{if   } x<0, \\
  \frac{1}{\Gamma(\frac{k}{2})}\int_{0}^{0.5} \gamma(\frac{k}{2} , \frac{x^2}{2\cos^2(\pi u)} )f_U(u) du & \text{if   } x\geq 0,
\end{array}
\right.
\end{eqnarray}

and the second term of the Equation (\ref{eqnFXproof}) is equal to

\begin{eqnarray}
\int_{0.5}^{1} P(\sqrt{V} \cos(\pi u)\leq x)f_U(u) du &=&
\left\{
\begin{array}{cc}
 \int_{0.5}^{1} P(\sqrt{V} \geq \frac{x}{\cos(\pi u)} )f_U(u) du & \text{if   } x<0 \\
  \int_{0.5}^{1}  f_U(u) du & \text{if   } x\geq 0
\end{array}
\right.
\nonumber \\ &=&
\left\{
\begin{array}{cc}
  \int_{0.5}^{1} f(u) du
-\int_{0.5}^{1} P(\sqrt{V} < \frac{x}{\cos(\pi u)} )f_U(u) du & \text{if   } x<0 \\
  \int_{0.5}^{1}  f_U(u) du & \text{if   } x\geq 0
\end{array}
\right.
\nonumber \\ &=&
\left\{
\begin{array}{cc}
\frac{I_{0.5}(2(1-\mu) , 2\mu )}{B(2\mu , 2(1-\mu))}
-\frac{1}{\Gamma(\frac{k}{2})}
\int_{0.5}^{1} \gamma(\frac{k}{2} , \frac{x^2}{2\cos^2(\pi u)} ) f_U(u) du & \text{if   } x<0, \\
  \frac{I_{0.5}(2(1-\mu) , 2\mu )}{B(2\mu , 2(1-\mu))}                     & \text{if   } x\geq 0.
\end{array}
\right.
\end{eqnarray}

Then, the cdf can be computed as follows:

\begin{eqnarray*}
\label{eqcdf}
F_{X}(x)&=&
\left\{
\begin{array}{cc}
\frac{I_{0.5}(2(1-\mu) , 2\mu )}{B(2\mu , 2(1-\mu))}
-\frac{1}{\Gamma(\frac{k}{2})}
\int_{0.5}^{1} \gamma(\frac{k}{2} , \frac{x^2}{2\cos^2(\pi u)} ) f_U(u) du
 & \text{if  } x<0, \\
\frac{I_{0.5}(2(1-\mu) , 2\mu )}{B(2\mu , 2(1-\mu))}+\frac{1}{\Gamma(\frac{k}{2})}\int_{0}^{0.5} \gamma(\frac{k}{2} , \frac{x^2}{2\cos^2(\pi u)} )f_U(u) du
&  \text{if  }  x\geq 0.
\end{array}
\right.  \quad \blacksquare
\end{eqnarray*}

}

\theorem{
The probability density function (pdf) of $PGN(\mu , k)$ has a two-peice form as
\begin{eqnarray}
\label{eqpdf}
f_X(x)=\frac{ (x^2)^{\frac{k-1}{2}}}{\Gamma(\frac{k}{2}) 2^{\frac{k}{2}-1}}\times
\left\{
\begin{array}{cc}
 \int_{0.5}^{1} \frac{1}{(\cos^2(\pi u))^{\frac{k}{2}}}
e^{-(\frac{x^2}{2\cos^2(\pi u)})} f_U(u) du
 & \text{if   } x<0, \\
 \int_{0}^{0.5} \frac{1}{(\cos^2(\pi u))^{\frac{k}{2}}}
e^{-(\frac{x^2}{2\cos^2(\pi u)})} f_U(u) du
& \text{if   } x\geq 0.
\end{array}
\right.
\end{eqnarray}

}
\proof{
The proof of this result follows directly by the derivation.
By the derivation of $F_X(x)$, $f_X(x)$ can be computed as follows:

\begin{eqnarray*}
f_{X}(x)&=&
\left\{
\begin{array}{cc}
-\frac{1}{\Gamma(\frac{k}{2})}
\int_{0.5}^{1} \frac{d}{dx} \gamma(\frac{k}{2} , \frac{x^2}{2\cos^2(\pi u)} ) f_U(u) du
 & x<0 \\
\frac{1}{\Gamma(\frac{k}{2})}\int_{0}^{0.5} \frac{d}{dx} \gamma(\frac{k}{2} , \frac{x^2}{2\cos^2(\pi u)} )f_U(u) du
&  x\geq 0
\end{array}
\right.  \\
&=&
\left\{
\begin{array}{cc}
 \frac{|x| (x^2)^{\frac{k}{2}-1}}{\Gamma(\frac{k}{2}) 2^{\frac{k}{2}-1}}
\int_{0.5}^{1} \frac{1}{(\cos^2(\pi u))^{\frac{k}{2}}}
e^{-(\frac{x^2}{2\cos^2(\pi u)})} f_U(u) du
 & x<0, \\
 \frac{|x| (x^2)^{\frac{k}{2}-1}}{\Gamma(\frac{k}{2}) 2^{\frac{k}{2}-1}}
\int_{0}^{0.5} \frac{1}{(\cos^2(\pi u))^{\frac{k}{2}}}
e^{-(\frac{x^2}{2\cos^2(\pi u)})} f_U(u) du
&  x\geq 0.
\end{array}
\right. \quad \blacksquare
\end{eqnarray*}

}

\corollary{
We could equivalently write the density of $PGN(\mu , k)$ as follows:

\begin{eqnarray}
\label{eqpdf2}
f_X(x) &=&\frac{ (x^2)^{\frac{k-1}{2}}}{\Gamma(\frac{k}{2}) 2^{\frac{k}{2}-1}}\times
\int_{0.5}^{1} \frac{1}{(\cos^2(\pi u))^{\frac{k}{2}}}
e^{-(\frac{x^2}{2\cos^2(\pi u)})} \left(f_U(u)1_{\{ x<0 \}} +f_U(1-u)1_{\{ x\geq 0 \}}\right) du
 \\ \label{eqpdf3} &=&
\frac{ (x^2)^{\frac{k-1}{2}}}{\Gamma(\frac{k}{2}) 2^{\frac{k}{2}-1}}\times
\int_{0}^{0.5} \frac{1}{(\cos^2(\pi u))^{\frac{k}{2}}}
e^{-(\frac{x^2}{2\cos^2(\pi u)})} \left(f_U(1-u)1_{\{ x<0 \}} +f_U(u)1_{\{ x\geq 0 \}}\right) du.
\end{eqnarray}

}

\proof{
By changing  variable $v$ to $1-u$ and knowing that $\cos(\pi (1-v))=-\cos(\pi v)$, we obtain

\begin{eqnarray*}
\int_{0.5}^{1} \frac{1}{(\cos^2(\pi u))^{\frac{k}{2}}}
e^{-(\frac{x^2}{2\cos^2(\pi u)})} f_U(u) du
&=&
\int_{0}^{0.5} \frac{1}{((\cos(\pi (1-v)))^2)^{\frac{k}{2}}}
e^{-(\frac{x^2}{2(\cos(\pi (1-v)))^2})} f_U(1-v) dv
\\ &=&
\int_{0}^{0.5} \frac{1}{(\cos^2(\pi v))^{\frac{k}{2}}}
e^{-(\frac{x^2}{2\cos^2(\pi v)})}  f_U(1-v) dv.
\end{eqnarray*}
In a similar way we have
\begin{eqnarray*}
\int_0^{0.5} \frac{1}{(\cos^2(\pi u))^{\frac{k}{2}}}
e^{-(\frac{x^2}{2\cos^2(\pi u)})} f_U(u) du
&=&
\int_{0.5}^{1} \frac{1}{(\cos^2(\pi v))^{\frac{k}{2}}}
e^{-(\frac{x^2}{2\cos^2(\pi v)})}  f_U(1-v) dv. \quad \blacksquare
\end{eqnarray*}

}

\begin{figure}[htp]
\captionsetup[subfloat]{captionskip=-27pt}
\begin{center}
  \label{figur:100}\caption{Shapes of $PGN(\mu , k)$ densities for different combinations of  parameters.}

  \subfloat[Shapes of $PGN(\mu , k)$ densities for different values of parameters. ]{\label{figur:101}\includegraphics[width=50mm, height=15em]{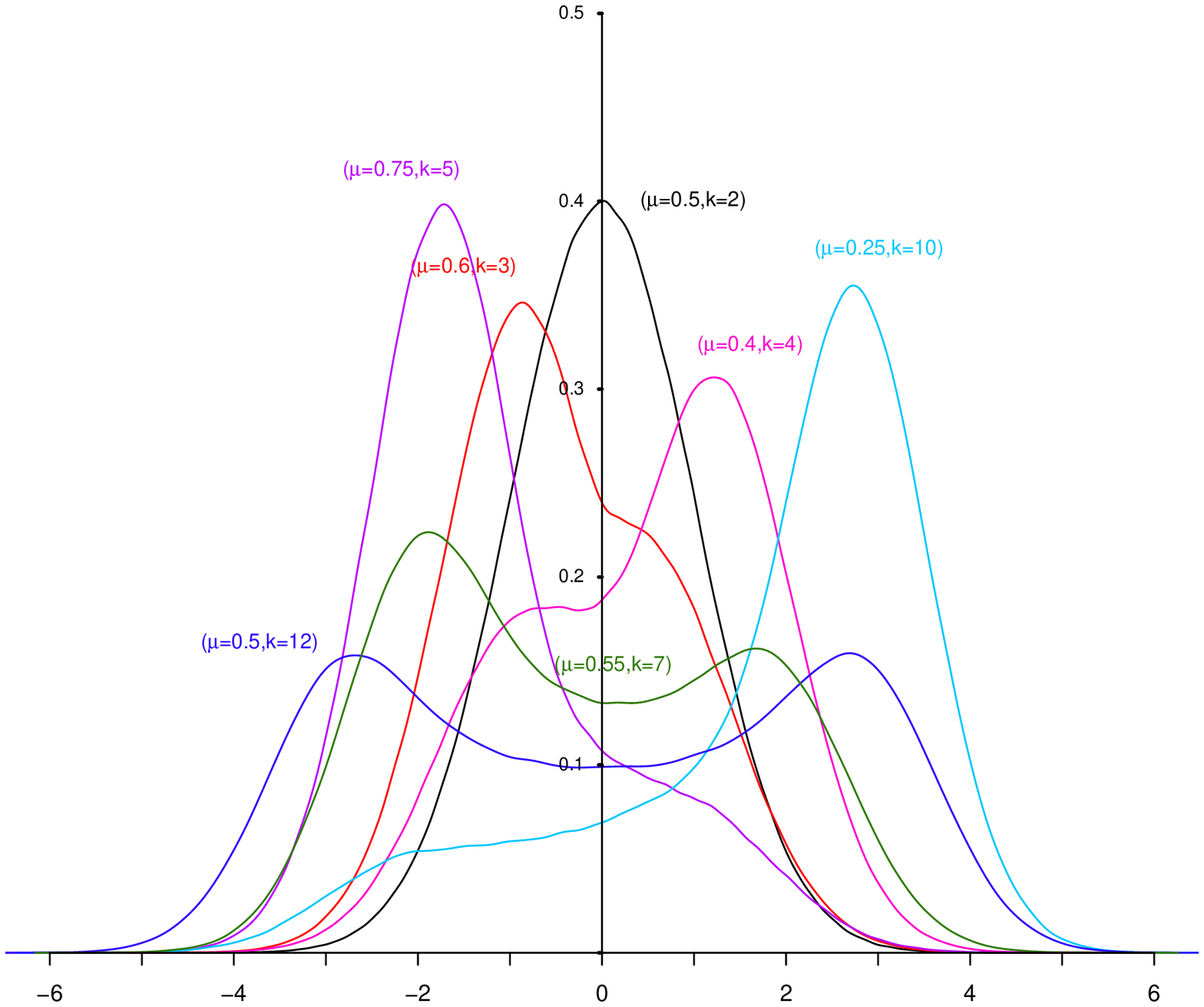}} \hspace{3cm} \vspace{-5pt}
    \subfloat[Shapes of $PGN(\mu=0.5 , k)$ densities for different values of $k$.]{\label{figur:102}\includegraphics[width=50mm, height=15em]{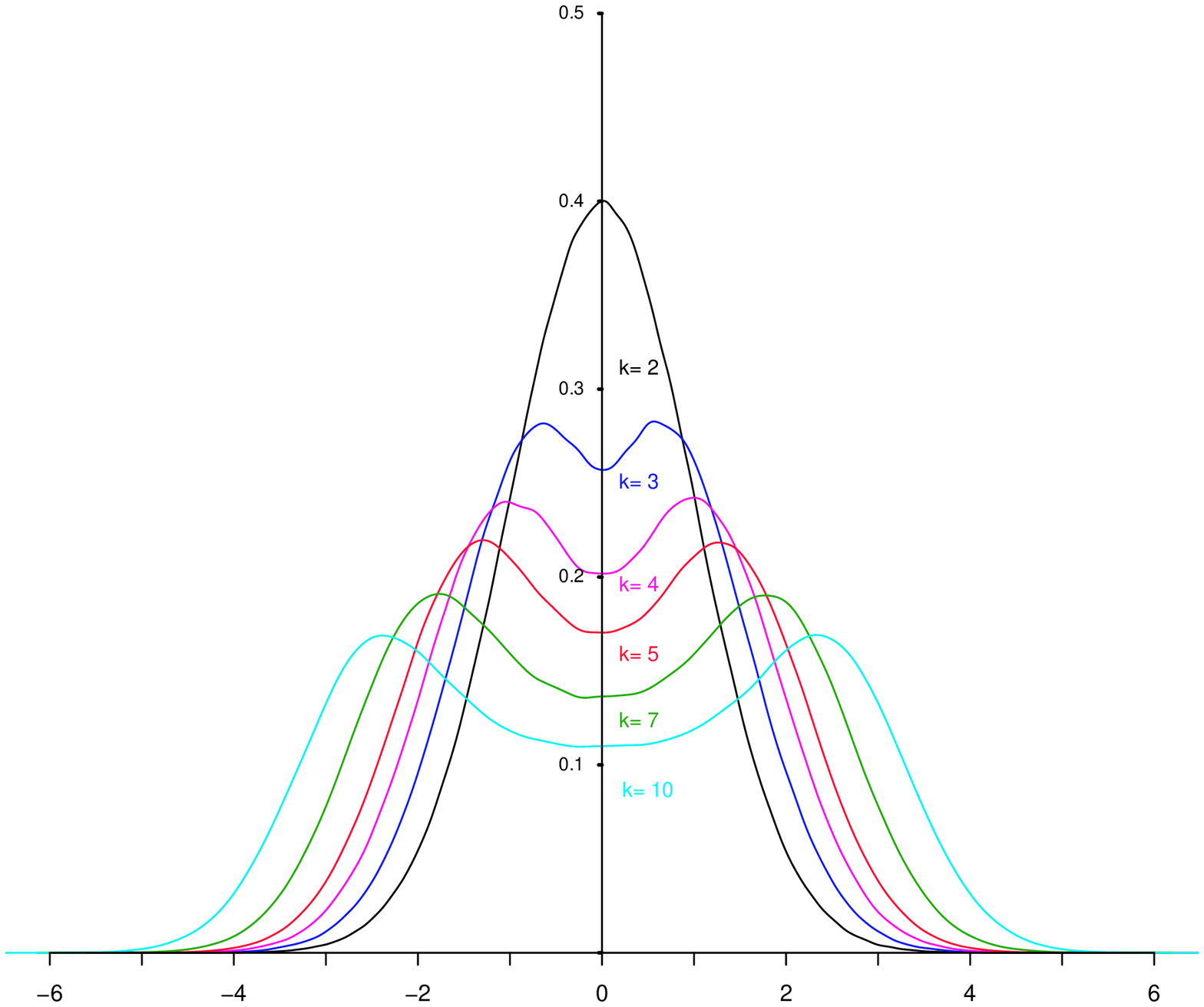}}
     \\
    \subfloat[Shapes of $PGN(\mu = 0.75 , k)$ densities for different values of $k$.]{\label{figur:103}\includegraphics[width=50mm, height=15em]{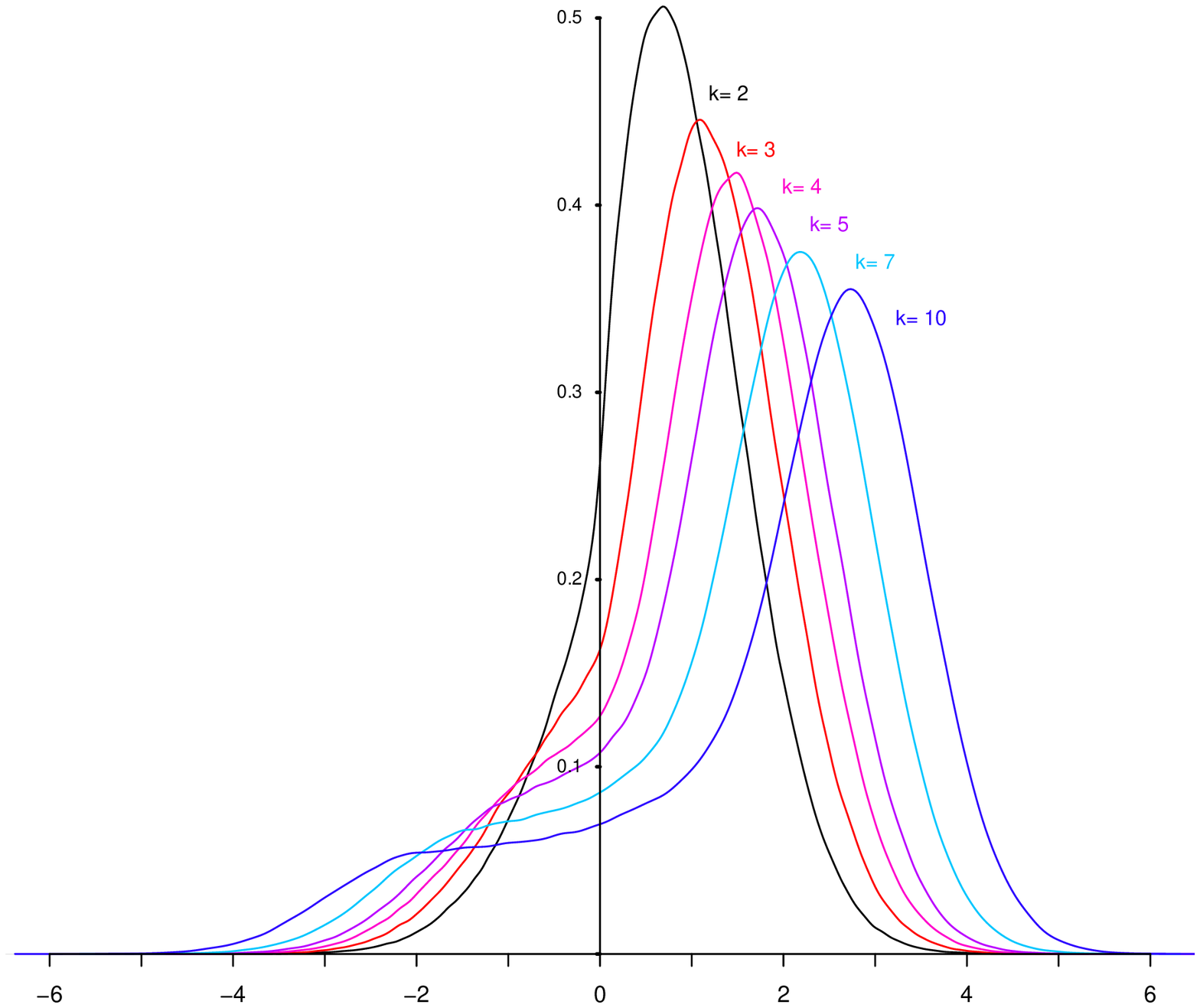}} \hspace{3cm} \vspace{-5pt}
  \subfloat[Shapes of $PGN(\mu = 0.25 , k)$ densities for different values of $k$.]{\label{figur:104}\includegraphics[width=50mm, height=15em]{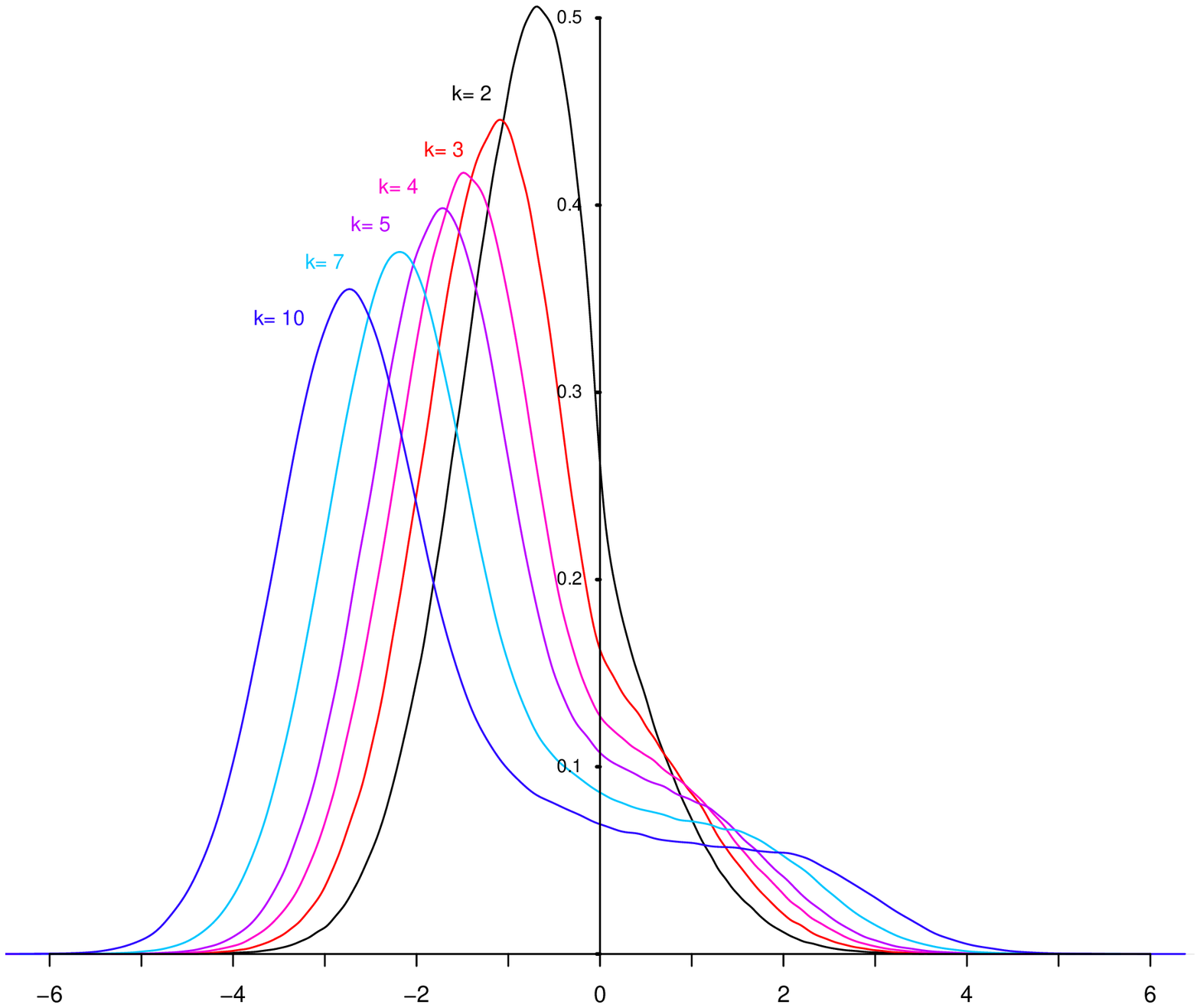}}
  \\
    \subfloat[Shapes of $PGN(\mu , k=5)$ densities for different values of $\mu$. ]{\label{figur:105}\includegraphics[width=50mm, height=15em]{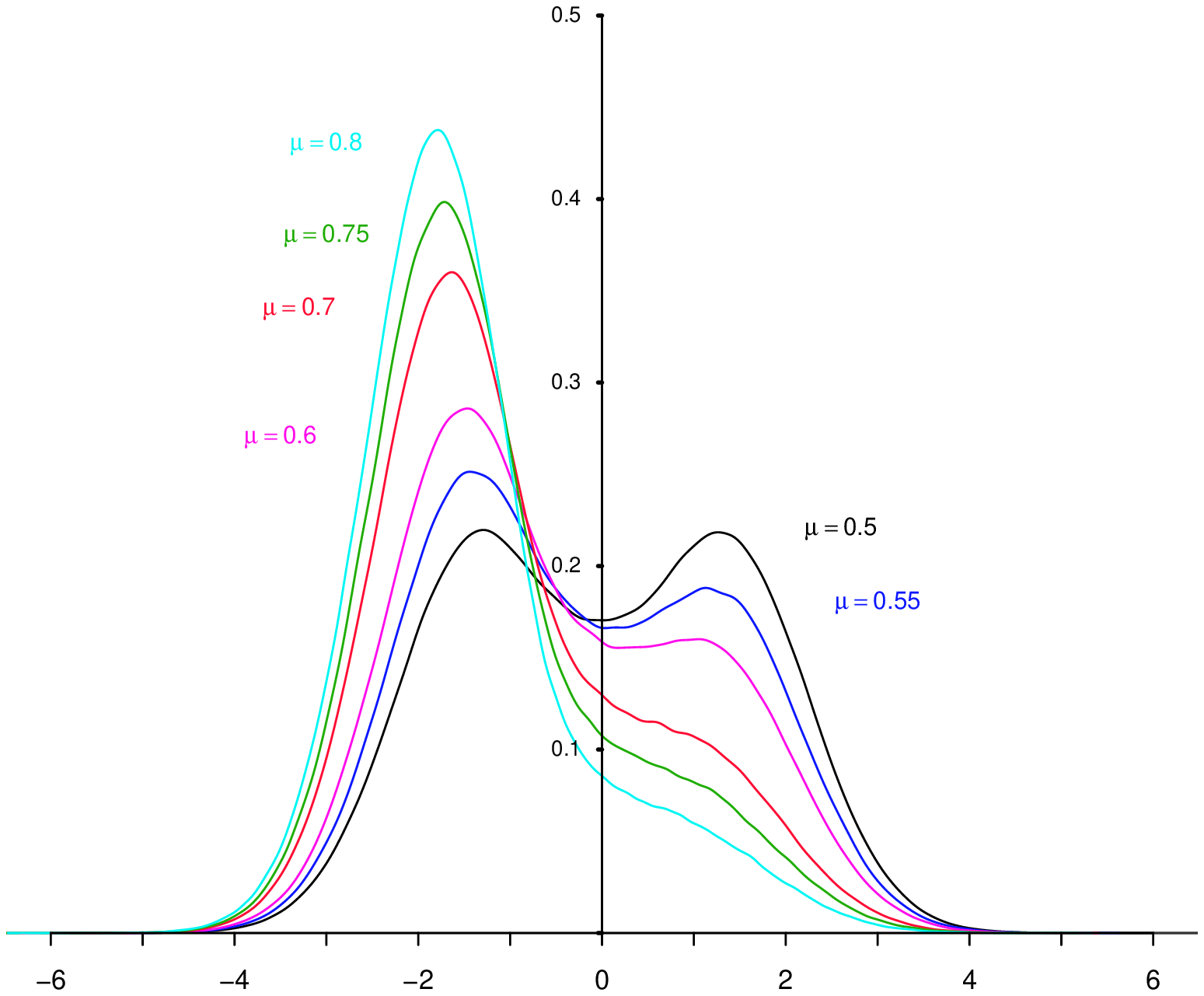}} \hspace{3cm} \vspace{-5pt}
    \subfloat[Shapes of $PGN(\mu  , k=5)$ densities for different values of $\mu$.]{\label{figur:106}\includegraphics[width=50mm, height=15em]{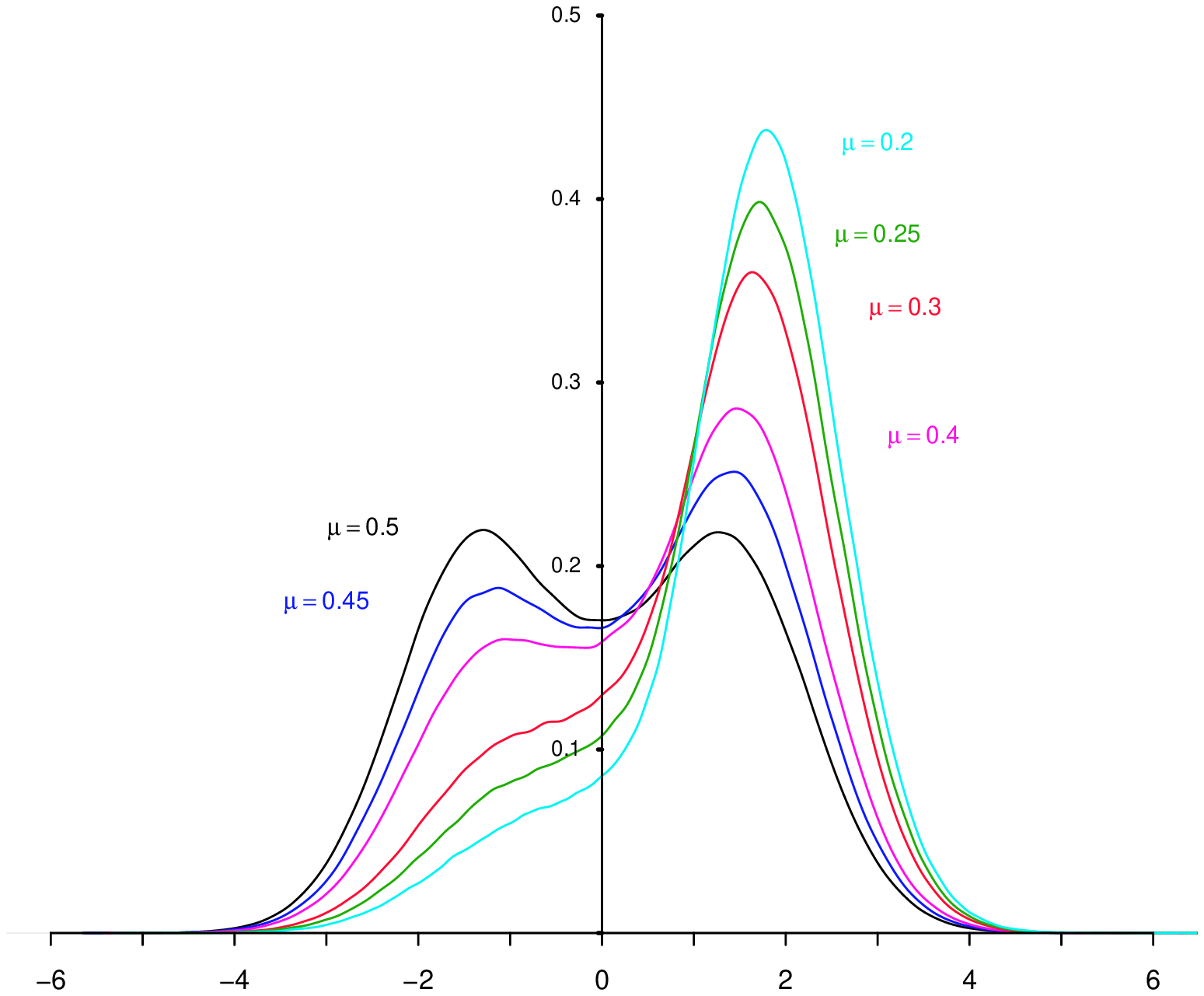}}
\end{center}
\end{figure}

\begin{figure}[htp]
\captionsetup[subfloat]{captionskip=-20pt}
\begin{center}
  \label{figur:200}\caption{Shapes of $PGN(\mu , k)$ densities for different combinations of  parameters.}

   \subfloat[Shapes of $PGN(\mu  , k=10)$ densities for different values of $\mu$.]{\label{figur:201}\includegraphics[width=50mm]{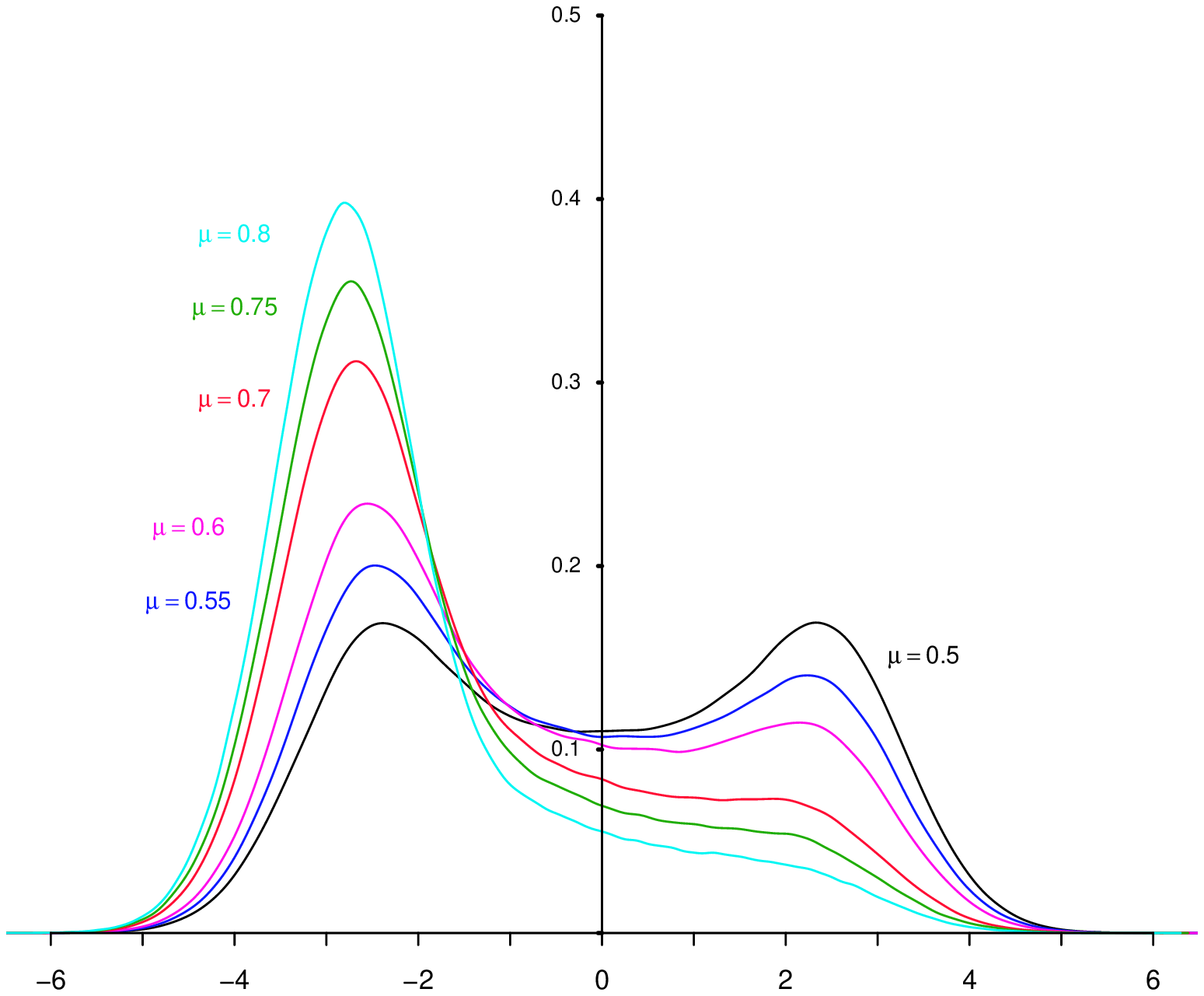}} \hspace{3cm} \vspace{-5pt}
  \subfloat[Shapes of $PGN(\mu  , k=10)$ densities for different values of $\mu$.]{\label{figur:202}\includegraphics[width=50mm]{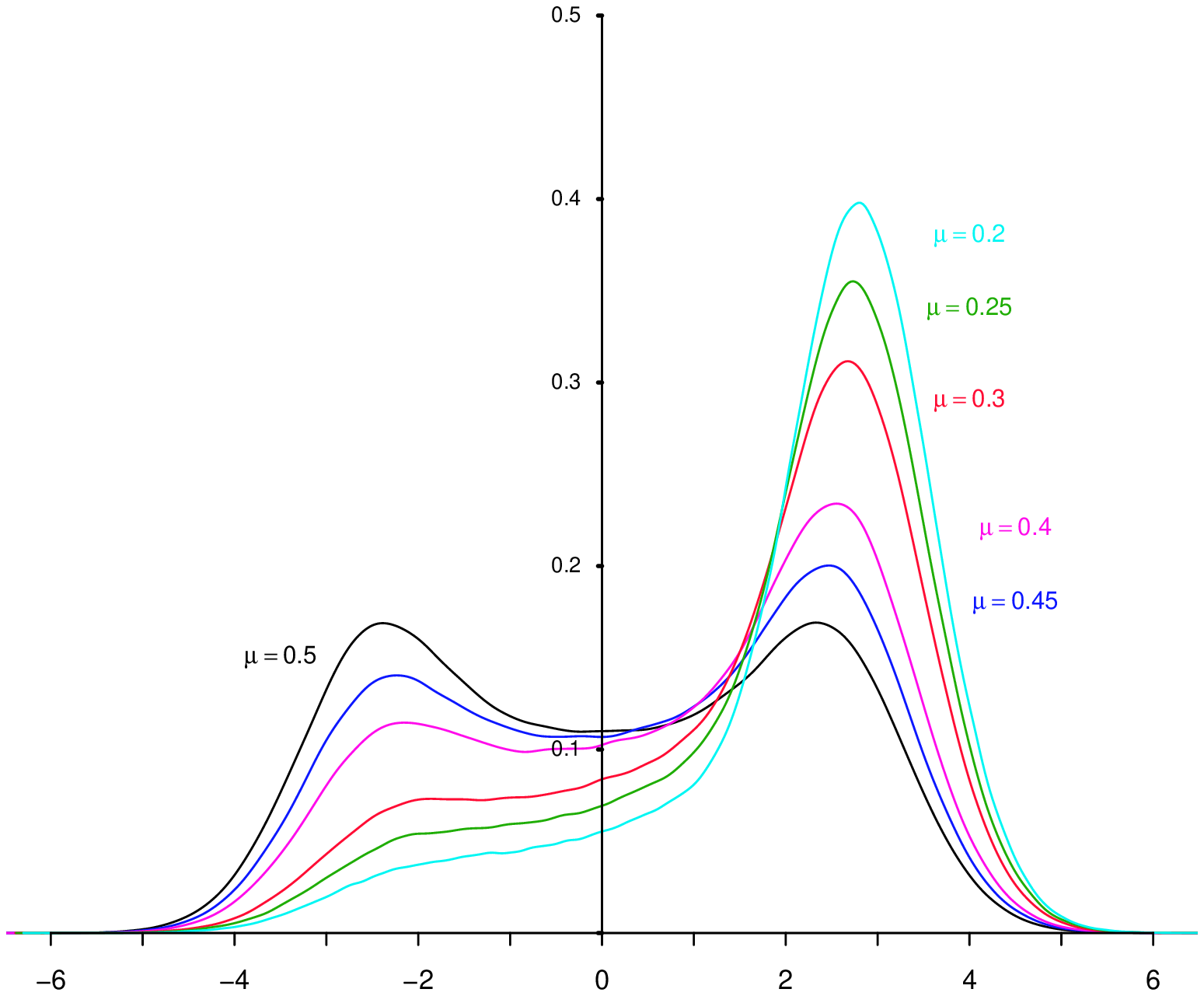}}
  \\
  \subfloat[Shapes of $PGN(\mu =0.5, k)$ densities for different values of $k <2$. ]{\label{figur:203}\includegraphics[width=50mm]{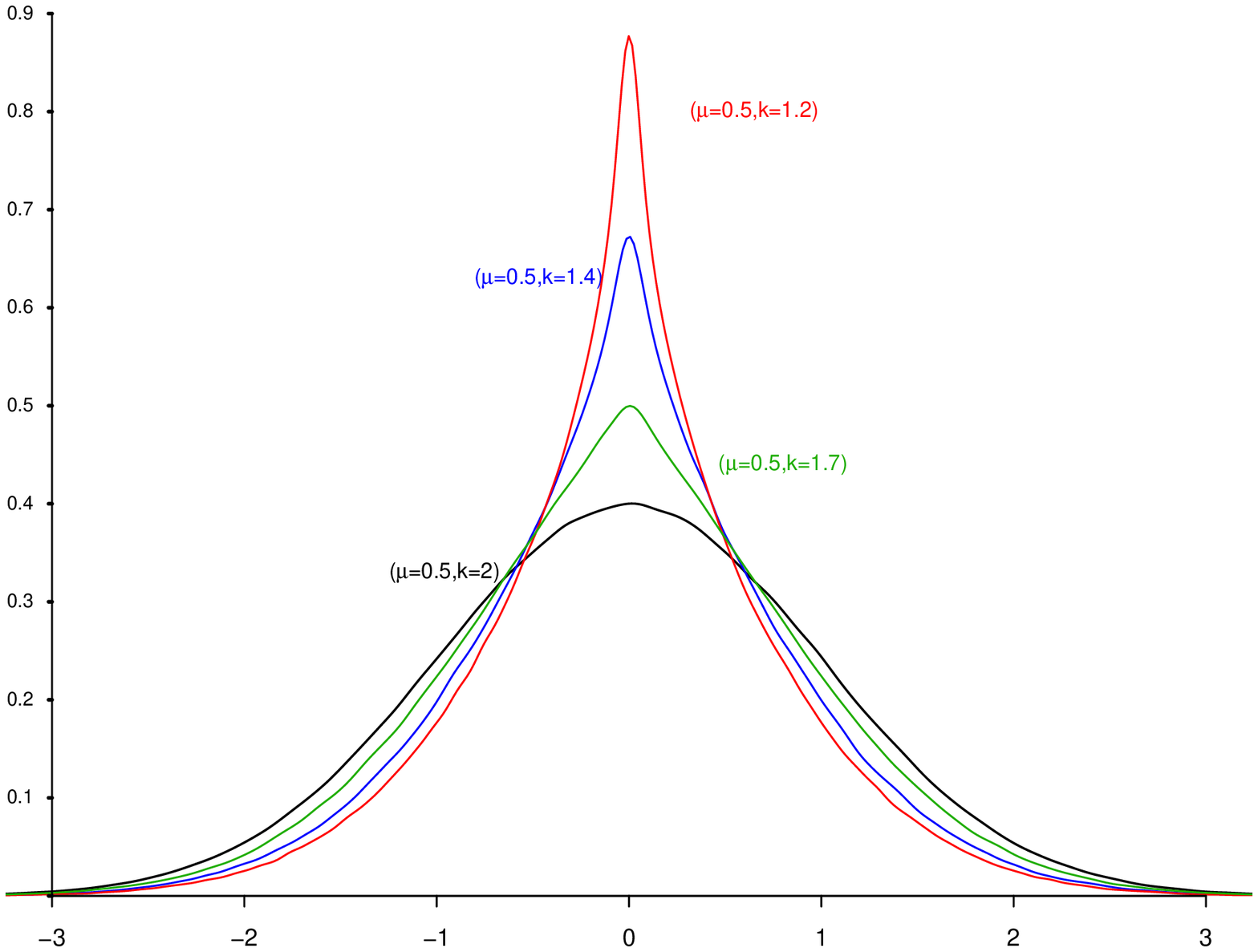}} \hspace{3cm} \vspace{-5pt}
  \subfloat[Shapes of $PGN(\mu =0.75 , k)$ densities for different values of $k <2$.]{\label{figur:204}\includegraphics[width=50mm]{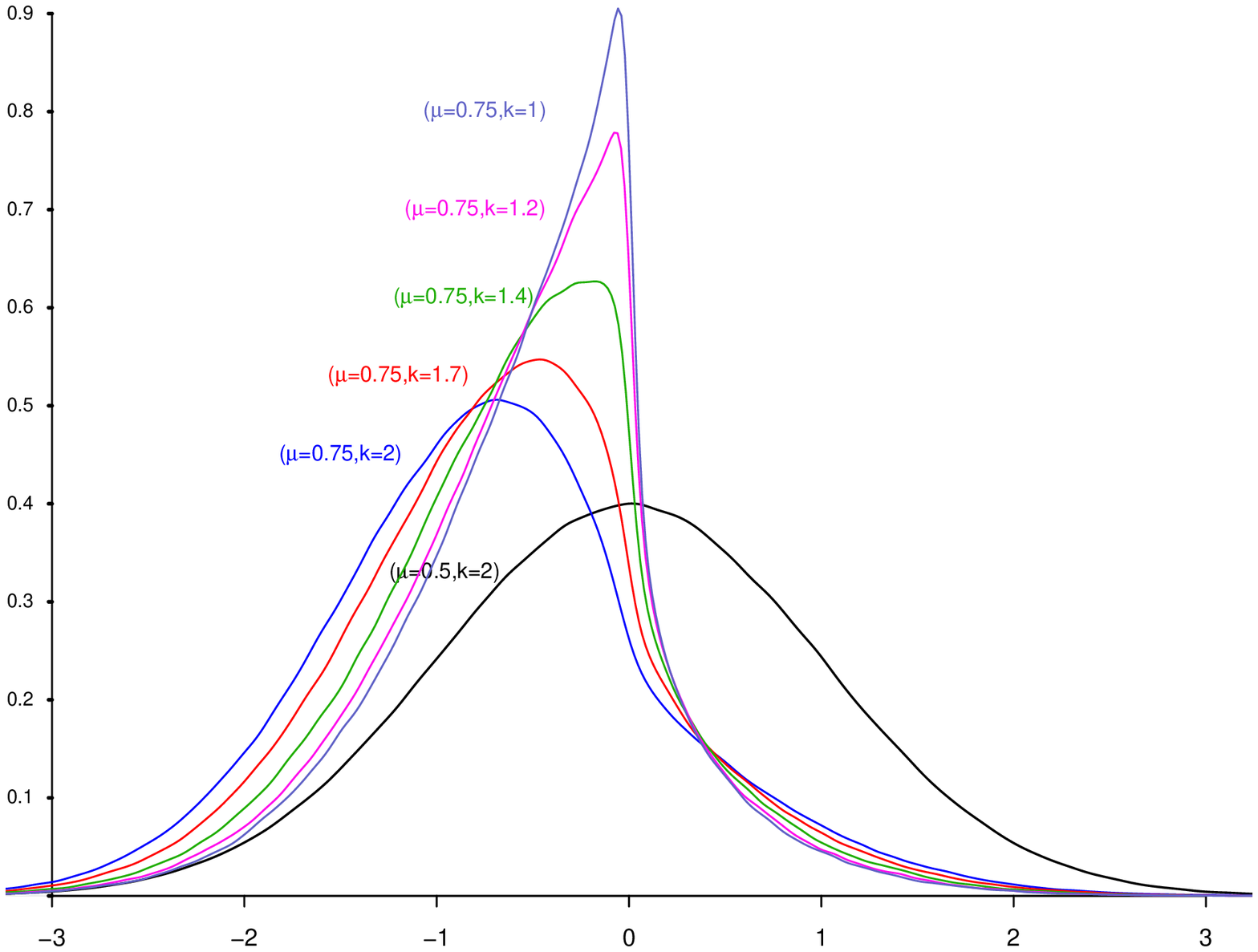}}
  \\
   \subfloat[Shapes of $PGN(\mu =0.25 , k)$ densities for different values of $k <2$.]{\label{figur:205}\includegraphics[width=50mm]{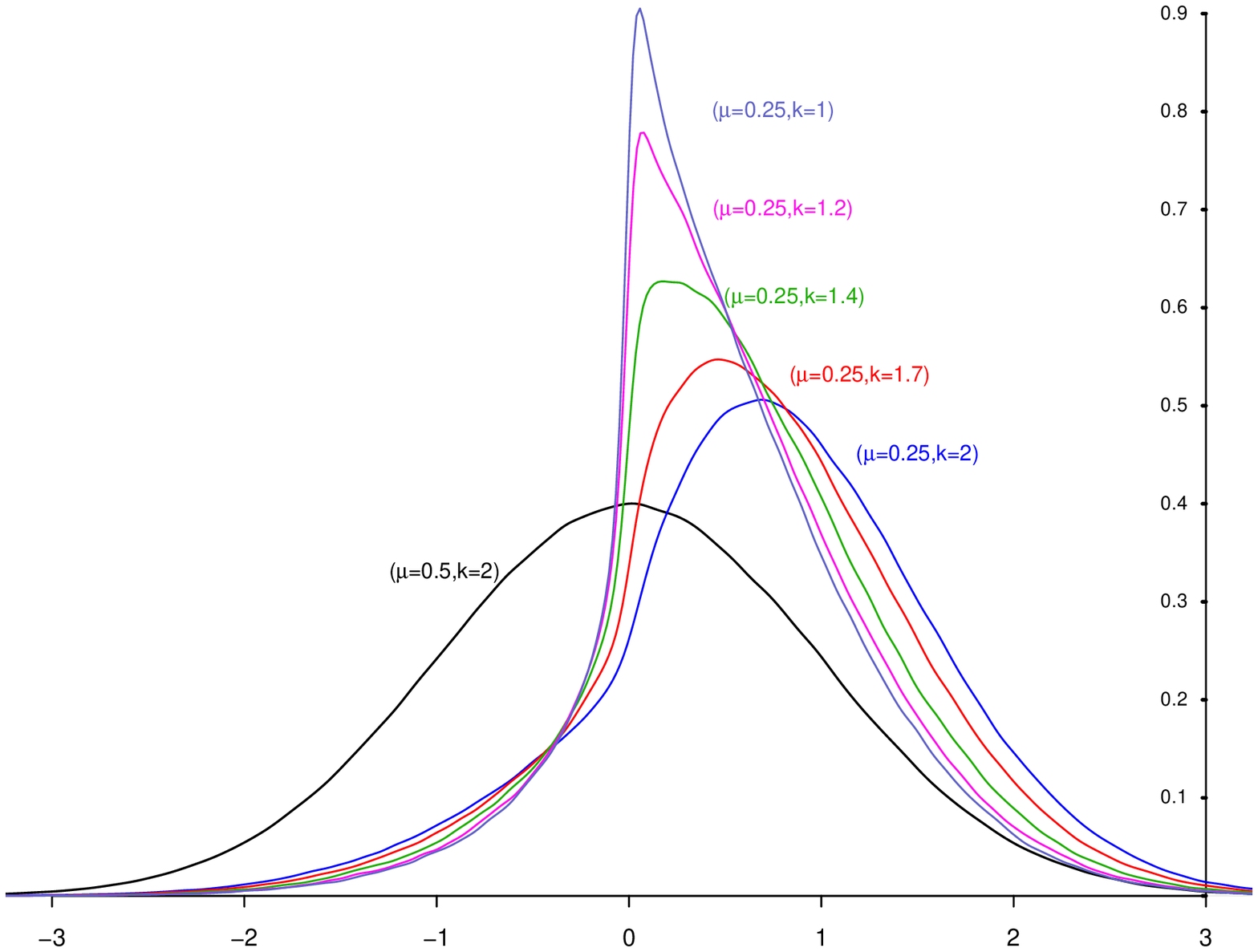}}
\end{center}
\end{figure}

The probability density function of $PGN(\mu , k)$ distribution for different choices of the parameters $\mu$ and $k$ are plotted in Figures 1 and 2.
%\ref{figur:100} and \ref{figur:200}
Figure \ref{figur:101} shows some specific densities of the distribution  $PGN(\mu , k)$ and illustrates  the flexibility of this family.
The distribution may be
symmetric, skewed to the left, or skewed to the right. At the same time, it can be unimodal or
bimodal.
It is noteworthy that the densities can display quite
different shapes depending on the values of the two parameters. In particular, it can
be symmetric when $\mu=0.5$   or asymmetric when $\mu \neq 0.5$.
For $\mu<0.5$, it is easy that to see $P(U<\frac{1}{2})>\frac{1}{2}$, or equivalently,
 $P(0<U<\frac{1}{2})>P(\frac{1}{2}<U<1)$. Accordingly,   $P(\cos (\pi U)>0)>P(\cos (\pi U)<0)$
and hence  $P(X>0)>P(X<0)$. So, $\mu$ controls the symmetry of the distribution. It confirmed
by the plots of $PGN(\mu , k)$ in Figures \ref{figur:103}, \ref{figur:104}, \ref{figur:105}, \ref{figur:106}, \ref{figur:201} and \ref{figur:202}. As we expected,
it is observed that for $\mu <0.5$ the distribution is left skewed and for $\mu >0.5$ is right skewed. The distribution also
approaches symmetric  when $\mu$ tends to $0.5$.  Asymmetry increases as $\mu$ tends to zero or one. So we call $\mu$  the asymmetric  parameter because it regulates the asymmetry of the density function.

We call $k$ the peak parameter because by increasing $k$ , the distribution  tends to have two peaks. These
peaks are distinguished for $\mu$ close to $0.5$, but as $\mu$ tends to $0$ or $1$ ,the small peak tends to smooth out.
From the plots of $PGN(\mu , k)$ in Figures \ref{figur:203}, \ref{figur:204} and \ref{figur:205}
it is observed that for $k <2$ and $\mu =0.5$  the distribution has a symmetric shape and the distribution tails  are shorter than the tails of a Normal distribution. It is also observed that by reducing $k$ the distribution  tends to has a sharp peak. For $k<2$ and $\mu \neq 0.5$ the distribution is skewed and by reducing $k$ the distribution  tends to has a sharp peak.

Figure \ref{figur:102} shows under $\mu=0.5$, the distribution is symmetric and for $k = 2$,  it is also a unimodal distribution (it has one peak). However,
by increasing $k$, the distribution
bifurcates, becomes bimodal and the distance between peaks increase.
This happens since as $k$  increases, the mode of the distribution of the random variable $\sqrt{V}$ increases, and consequently,  the   distance between
the mode and the zero point increases. Accordingly, by multiplying $\cos (\pi U)$, which is between $-1$ and $1$, to $\sqrt{V}$,  two distinct modes are created (at both sides of zero point) and the distance between two modes increases.

The $PGN$ distribution becomes
bimodal for certain values of the parameters $\mu$ and $k$, and the analytical solution of $\mu$ and $k$ ,
where the distribution becomes bimodal,cannot
be solved algebraically. However,  by observing density function of $PGN(\mu,k)$ for different values of the parameters $\mu$ and $k$,  we obtain Table (\ref{bimodalreg}). The table
 shows a grid of values where the distribution
is bimodal.
In the table, “1” indicates that PGN distribution is unimodal and “2” indicates that PGN distribution is bimodal.

\begin{table}[]
\centering 
\caption{Number of mods of $PGN(\mu, k)$ for different values of the parameters. “2” indicates where bimodality occurs and “1” indicates where unimodality occurs.} \label{bimodalreg}
\begin{tabular}{ccccccccccccccc}
  & k    & 1 & 2 & 3 & 4 & 5 & 6 & 7 & 8 & 9 & 10 & 11 & 12 & 13
 \\
 $\mu$ & &  &  &  &  &  &  &  &  &  &  &  &  & \\ \hline
0.50 & &1  & 1  & \textbf{2}  & 2  & 2  &  2 & 2  &  2 &  2 &  2  &  2  &  2  &  2  \\
0.51, 0.49 & &1  & 1  & \textbf{2}  & 2  & 2  &  2 & 2  &  2 &  2 &  2  &  2  &  2  &  2  \\
0.52, 0.48 & &1  & 1  & 2\textbf{}  & 2  & 2  &  2 & 2  &  2 &  2 &  2  &  2  &  2  &  2  \\
0.53, 0.47 & &1  & 1  & \textbf{2}  & 2  & 2  &  2 & 2  &  2 &  2 &  2  &  2  &  2  &  2  \\
0.54, 0.46 & &1  & 1  & \textbf{2}  & 2  & 2  &  2 & 2  &  2 &  2 &  2  &  2  &  2  &  2  \\
0.55, 0.45 & &1  & 1  & \textbf{2}  & 2  & 2  &  2 & 2  &  2 &  2 &  2  &  2  &  2  &  2  \\
0.56, 0.44& &1  & 1  & \textbf{2}  & 2  & 2  &  2 & 2  &  2 &  2 &  2  &  2  &  2  &  2  \\
0.57, 0.43 & &1  & 1  & 1   & \textbf{2}  & 2  & 2  & 2  &  2 & 2  &  2  &  2  &  2  &  2  \\
0.58, 0.42 & &1  & 1  & 1   & \textbf{2}  & 2  & 2  & 2  &  2 & 2  &  2  &  2  &  2  &  2  \\
0.59, 0.41 & &1  & 1  & 1   & \textbf{2}  & 2  & 2  & 2  &  2 & 2  &  2  &  2  &  2  &  2  \\
0.60, 0.40 & &1  & 1  & 1   & \textbf{2}  & 2  & 2  & 2  &  2 & 2  &  2  &  2  &  2  &  2  \\
0.61, 0.39 & &1  & 1  & 1   & 1  & \textbf{2}  & 2  & 2  &  2 & 2  &  2  &  2  &  2  &  2  \\
0.62, 0.38 & &1  & 1  & 1   & 1  & \textbf{2}  & 2  & 2  &  2 & 2  &  2  &  2  &  2  &  2  \\
0.63, 0.37 & &1  & 1  & 1   & 1  & 1  & \textbf{2}  & 2  &  2 & 2  &  2  &  2  &  2  &  2  \\
0.64, 0.36 & &1  & 1  & 1   & 1  & 1  & 1  & \textbf{2}  &  2 & 2  &  2  &  2  &  2  &  2  \\
0.65, 0.35 & &1  & 1  & 1   & 1  & 1  & 1  & 1  &  \textbf{2} & 2  &  2  &  2  &  2  &  2  \\
0.66, 0.34 & &1  & 1  & 1   & 1  & 1  & 1 & 1  & 1  & \textbf{2}  &  2  &  2  &  2  &  2  \\
0.67, 0.33 & &1  & 1  & 1  & 1  & 1  & 1  & 1  & 1  & 1  & 1   & \textbf{2}   &  2  & 2   \\
0.68, 0.32 & &1  & 1  & 1  & 1  & 1  & 1  & 1  & 1  & 1  & 1   &  1  & 1   &  \textbf{2}  \\
0.69-1.00,0.00-0.31 & & 1  & 1  & 1  & 1  & 1  & 1  & 1  & 1  & 1  &  1  &  1  & 1   &  1   \\
\end{tabular}
\end{table}

\theorem{
\label{theoremsymmetric}
If $\mu=0.5$, the cdf of the  $PGN(\mu,k)$ model  is  given by
\begin{eqnarray}
\label{pdfsymmetric}
f_X(x)=
 \frac{(x^2)^{\frac{k-1}{2}}e^{-\frac{x^2}{2}}}{\sqrt{\pi}\Gamma(\frac{k}{2}) 2^{\frac{k}{2}}}
  U(\frac{1}{2} , \frac{k+1}{2},\frac{x^2}{2}),
\end{eqnarray}
where $U(a,b,z)$ is Tricomi's (confluent hypergeometric) function, that is,

\begin{eqnarray*}
U(a,b,z)=
\frac{1}{\Gamma (a)}\int_0^{\infty} e^{-zu} u^{a-1} (1+u)^{b-a-1} du.
\end{eqnarray*}

}

\proof{

By Equation \ref{eqpdf3} and knowing that for $\mu=0.5$, $f_U(u)=1$, we have

\begin{eqnarray*}
f_X(x)&=&
\frac{|x| (x^2)^{\frac{k}{2}-1}}{\Gamma(\frac{k}{2}) 2^{\frac{k}{2}-1}}
\int_{0}^{0.5} \frac{1}{(\cos^2(\pi u))^{\frac{k}{2}}}
e^{-(\frac{x^2}{2\cos^2(\pi u)})}  du
\\
&=&
 \frac{ (x^2)^{\frac{k-1}{2}}e^{-\frac{x^2}{2}}}{\Gamma(\frac{k}{2}) 2^{\frac{k}{2}-1}}
\int_{0}^{0.5} (1+ \tan^2(\pi u))^{\frac{k}{2}}
e^{-\frac{x^2}{2}( \tan^2(\pi u))}  du
\\ &=&
 \frac{ (x^2)^{\frac{k-1}{2}}e^{-\frac{x^2}{2}}}{\Gamma(\frac{k}{2}) 2^{\frac{k}{2}-1}}
\int_{0}^{+\infty} (1+ t)^{\frac{k}{2}}
e^{-\frac{x^2}{2} t}  \frac{1}{\pi 2(1+t)\sqrt{t}}dt
\\
&=&
 \frac{ (x^2)^{\frac{k-1}{2}}e^{-\frac{x^2}{2}}}{\Gamma(\frac{k}{2}) 2^{\frac{k}{2}-1}}
\frac{1}{2\pi}
\int_{0}^{+\infty} e^{-\frac{x^2}{2} t} t^{\frac{1}{2}-1}(1+ t)^{\frac{k}{2}-1}
  dt
  \\
&=&
 \frac{ (x^2)^{\frac{k-1}{2}}e^{-\frac{x^2}{2}}}{\Gamma(\frac{k}{2}) 2^{\frac{k}{2}}}
\frac{1}{\pi}
\int_{0}^{+\infty} e^{-\frac{x^2}{2} t} t^{\frac{1}{2}-1}(1+ t)^{\frac{k+1}{2}-\frac{1}{2}-1}
  dt
  \\
&=&
 \frac{ (x^2)^{\frac{k-1}{2}}e^{-\frac{x^2}{2}}}{\Gamma(\frac{k}{2}) 2^{\frac{k}{2}}}
\frac{1}{\sqrt{\pi}}  U(\frac{1}{2} , \frac{k+1}{2},\frac{x^2}{2}). \quad \blacksquare
\end{eqnarray*}

}

In the rest of this section, we provide some
distributional properties of the $PGN(\mu,k)$
distribution.

\corollary{
\label{corollarysymmetric}
If $\mu=0.5$, the density of the $PGN(\mu,k)$ model is symmetric about zero.
}

\proof{ The proof is simple since the cdf \ref{pdfsymmetric} satisfy  $f_X(-x)=f_X(x)$. $\quad \blacksquare$

}

Table \ref{tricomi} shows some special cases of the  Tricomi confluent hypergeometric function (Wolfram Research, 2020).

\corollary{
 The standard normal distribution is a special case of the $PGN(\mu , k)$ distribution.
}

\proof{
 For the special case $\mu=0.5$ and  $k=2$,
$U(\frac{1}{2},\frac{3}{2},z)=\frac{1}{\sqrt{z}}$, we obtain the standard normal distribution since
$f_X(x)=\frac{(x^2)^{\frac{2-1}{2}}e^{-\frac{x^2}{2}}}{\sqrt{\pi}\Gamma(\frac{2}{2}) 2^{\frac{2}{2}}}\frac{1}{\sqrt{\frac{x^2}{2}}}=\frac{e^{-\frac{x^2}{2}}}{\sqrt{2\pi}}$. $\quad \blacksquare$

}

\corollary{
\label{col23}
For the special case $\mu=0.5$ and  $k=4$, the $PGN(\mu , k)$ distribution  is an equal mixture of two different families, normal distribution and Chi distribution with 3 degree of freedom.
}

\proof{

According to $U(\frac{1}{2},\frac{5}{2},z)=\frac{2z+1}{2z^{\frac{3}{2}}}$,
we have
\begin{eqnarray*}
  f_X(x) &=&  \frac{(x^2)^{\frac{4-1}{2}}e^{-\frac{x^2}{2}}}{\sqrt{\pi}\Gamma(\frac{4}{2}) 2^{\frac{4}{2}}}
  U(\frac{1}{2} , \frac{5}{2},\frac{x^2}{2}) \\
   &=& \frac{(x^2)^{\frac{3}{2}}e^{-\frac{x^2}{2}}}{\sqrt{\pi}\Gamma(2) 2^{2}}
   \frac{2\frac{x^2}{2}+1}{2(\frac{x^2}{2})^{\frac{3}{2}}}
   \\
   &=& \frac{(x^2+1)e^{-\frac{x^2}{2}}}{2\sqrt{2\pi} }
   =\frac{1}{2} \frac{x^2e^{-\frac{x^2}{2}}}{\sqrt{2\pi} }
   +\frac{1}{2} \frac{e^{-\frac{x^2}{2}}}{\sqrt{2\pi} },
\end{eqnarray*}
which is an equal mixture of two different families, normal distribution and Chi distribution with 3 degree of freedom. $\quad \blacksquare$
}

With regard to Corollary \ref{col23},  the distribution is bimodal, as you will see in Figure(\ref{figur:204}). From the table
we can see for $k\geq 3$ the pdf  can be written as a mixture of two pdf and this result could be seen in Figure(\ref{figur:204}).
Thus, the PGN distribution provides great flexibility in modeling symmetric,
able to accommodate both unimodal and bimodal cases.

\theorem{
For $\mu \rightarrow 0$, the $PGN(\mu,K)$ distribution becomes $\chi_{k}$, the Chi distribution with
$k$ degrees of freedom, that is, the density is given by
$f_X(x)=\frac{x^{k-1}e^{-\frac{x^2}{2}}}{2^{\frac{k}{2}-1}\Gamma(\frac{k}{2})} $ for $x>0$. For $\mu \rightarrow 1$ the $PGN(\mu,K)$ distribution becomes $-\chi_{k}$.
}

\proof{
For $\mu \rightarrow 0$, the beta distribution becomes a  Degenerate distribution at point $0$ and for $\mu \rightarrow 1$, the beta distribution becomes a  Degenerate distribution at point $1$. $\quad \blacksquare$
}

\begin{table}
\centering
\caption{Special cases of the  Tricomi confluent hypergeometric function. In this table, $K_a$ is modified Bessel function of the second kind, that is, $K_{\alpha}(x)=\frac{\pi}{2} \frac{I_{-\alpha}(x)-I_{\alpha}(x)}{\sin \alpha x}$, where $I_{\alpha}(x)=\sum_{m=0}^{\infty} \frac{1}{m! \Gamma (m+\alpha +1)} (\frac{x}{2})^{2m+\alpha}$ is the modified Bessel functions of the first kind.} \label{tricomi}
{\renewcommand{\arraystretch}{2.7}
\begin{tabular} {| l ||c  |}
 \hline
        & $U(\frac{1}{2},\frac{k+1}{2},z)$    \\ \hline \hline
$k=1$   & $U(\frac{1}{2},1,z)=\frac{e^{\frac{z}{2}}K_0(\frac{z}{2})}{\sqrt{\pi}}$            \\ \hline
$k=2$   & $U(\frac{1}{2},\frac{3}{2},z)=\frac{1}{\sqrt{z}}$   \\ \hline
$k=3$   & $U(\frac{1}{2},2,z)=\frac{e^{\frac{z}{2}}\left( K_0(\frac{z}{2}) +K_1(\frac{z}{2})\right)}{2\sqrt{\pi}}$             \\ \hline
$k=4$   & $U(\frac{1}{2},\frac{5}{2},z)=\frac{2z+1}{2z^{\frac{3}{2}}}$    \\ \hline
$k=5$   & $U(\frac{1}{2},3,z)
=\frac{e^{\frac{z}{2}}\left( z K_0(\frac{z}{2}) +(z+1)K_1(\frac{z}{2})\right)}{2\sqrt{\pi}z}$            \\ \hline
$k=6$   & $U(\frac{1}{2},\frac{7}{2},z)=\frac{4z(z+1)+3}{4z^{\frac{5}{2}}}$   \\ \hline
$k=7$   & $U(\frac{1}{2},4,z)=\frac{e^{\frac{z}{2}}\left( z(2z+1) K_0(\frac{z}{2}) +(z(2z+3)+4)K_1(\frac{z}{2})\right)}{4\sqrt{\pi}z^2}$             \\ \hline
$k=8$   & $U(\frac{1}{2},\frac{9}{2},z)=\frac{2z(4z^2+6z+9)+15}{8z^{\frac{7}{2}}}$   \\ \hline
$k=9$   & $U(\frac{1}{2},5,z)=\frac{e^{\frac{z}{2}}\left( z(2z(z+1)+3) K_0(\frac{z}{2}) +2(z(z(z+2)+4)+6)K_1(\frac{z}{2})\right)}{4\sqrt{\pi}z^3}$            \\ \hline
$k=10$  & $U(\frac{1}{2},\frac{11}{2},z)=
\frac{8z(z(2z(z+2)+9)+15)+105}{16z^{\frac{9}{2}}}$   \\ \hline
$k=11$  & $U(\frac{1}{2},6,z)=
\frac{e^{\frac{z}{2}}\left( z(z(4z^2+6z+15)+24) K_0(\frac{z}{2}) +(z(z(2z(2z+5)+27)+60)+96)K_1(\frac{z}{2})\right)}{8\sqrt{\pi}z^4}
$            \\ \hline
\end{tabular}
}
\end{table}

\theorem{
Let $X\sim PGN(\mu , k)$, the fraction $\frac{P(X>0)}{P(X<0)}$ only depends on $\mu$, we  denote
 $\varphi (\mu)= \frac{P(X>0)}{P(X<0)}$. We have

\begin{eqnarray}
\label{skewtheorem1}
  \varphi (\mu )>1  &\Longleftrightarrow & \mu < 0.5, \\
   \varphi (\mu )< 1  &\Longleftrightarrow & \mu > 0.5, \\
   \varphi (\mu )=1  &\Longleftrightarrow & \mu = 0.5.
\end{eqnarray}
}

\proof{
\begin{eqnarray*}
 \varphi (\mu )&=&\frac{P(X>0)}{P(X<0)}=\frac{P(\sqrt{V} \cos(\pi U)>0)}{P(\sqrt{V} \cos(\pi U)<0)} \\
   &=& \frac{P(U<0.5)}{P(U>0.5)}=\frac{\frac{1}{B(2\mu , 2(1-\mu))}\int_{0}^{0.5} u^{2\mu-1}(1-u)^{1-2\mu} du}{\frac{1}{B(2\mu , 2(1-\mu))}\int_{0.5}^{1} u^{2\mu-1}(1-u)^{1-2\mu} du} \\
   &=& \frac{\int_{0}^{0.5} (\frac{u}{1-u})^{2\mu-1} du}{\int_{0.5}^{1} (\frac{u}{1-u})^{2\mu-1} du} = \frac{\int_{0}^{0.5} (\frac{u}{1-u})^{2\mu-1} du}{\int_{0}^{0.5} (\frac{u}{1-u})^{1-2\mu} du} \\
   &=& \frac{\int_{0}^{0.5} h(u)^{2\mu-1} du}{\int_{0}^{0.5} h(u)^{1-2\mu} du},
\end{eqnarray*}
with $h(u)=\frac{u}{1-u}$ it is easy to check that $0<h(u)<1$ for every $u \in (0,0.5)$. Hence, if $2\mu -1<0$, i.e. if $\mu< 0.5$,
we have that $h(u)^{2\mu -1} > h(u)^{1-2\mu}$, $\forall u \in (0,0.5)$, so that $\varphi(\mu) >1$. On the other hand, the inequality is reversed if $\mu > 0.5$ and, finally, $\varphi(\mu) =1$
if $\mu=0.5$. $\quad \blacksquare$

}

\begin{figure}
  \centering
  \includegraphics[width=50mm]{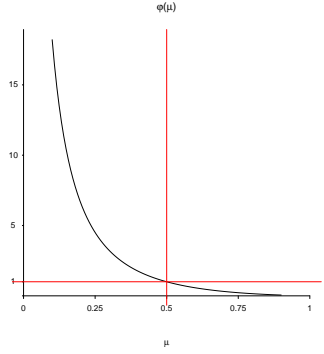}
  \caption{The shape of $\varphi (\mu)$.}\label{plotskewmu}
\end{figure}
Figure \ref{plotskewmu} shows $\varphi(\mu)$ for $\mu \in (0,1)$. Accordingly  we came to the conclusion that $\mu$ controls the
asymmetric of the distribution .
On the other hand, $\mu$ controls the allocation of mass to each side of the point $0$.
To be more specific, the asymmetry of the distribution  is only controlled by $\mu$.

\theorem{

 The  $PGN(\mu , k)$ distribution can not be symmetric around any $a\neq 0$.
 }

 \proof{
 To see this, let the  $PGN(\mu , k)$ distribution be symmetric around $a$. Let $a>0$, we have  $\forall x>0 \quad P(X<a-x) = P(X>a+x)$.
 For $x=a$ we obtain

 \begin{eqnarray}
\label{symmetry0}
     P(X<0) = P(X>2a).
 \end{eqnarray}

 The left hand side of \ref{symmetry0} only depend of $\mu$, so by changing $k$ it does not change, But as changing
 $k$, the  right hand side of \ref{symmetry0} will change. So the only way to the  $PGN(\mu , k)$ distribution
 be symmetric is for $a=0$. The proof of this for $a<0$ is similar to the proof for
 $a>0$. $\quad \blacksquare$
 }

\corollary{
From the above Proposition and Proposition (\ref{corollarysymmetric}), we deduce that the $PGN(\mu, k)$ is symmetric around 0, if and only if $\mu=0.5$
}

\lemma{
Let $X\sim PGN(\mu , k)$, then $-X\sim PGN((1-\mu) , k)$.
}
\proof{
Knowing that
$1-U \sim Beta(2(1-\mu),2\mu)$,
the proof of this result follows directly from definition (\ref{definitionpolar}) as follows:
\begin{eqnarray*}
-X=-\sqrt{V}\cos(\pi U)=\sqrt{V}\cos(\pi -\pi U)=\sqrt{V}\cos(\pi(1- U))\sim  PGN(1-\mu,k). \quad \blacksquare
\end{eqnarray*}
}

From the above proof we can deduce that if the distribution of $X$ is bimodal, the distribution of $-X$ is bimodal too. We can
see this in Table (\ref{bimodalreg}).
By the table the region of bimodality of the PGN distribution for $\mu$ and $1-\mu$ is same. Since by the above lemma,
if $X\sim PGN(\mu , k)$, $-X\sim PGN((1-\mu) , k)$ and $k$,that controls  bimodality of the distribution, does not change.

\lemma{
Let $V\sim \chi_{k}^{2}$ and $U \sim Beta (2\mu , 2(1-\mu))$ the following properties are deduced immediately from
the definition:

 \begin{eqnarray*}
    \sqrt{V}\sin(\pi ( U-\frac{1}{2})) &\sim&  PGN(\mu ,k). \\
    \sqrt{V}\sin(\pi (\frac{1}{2}- U)) &\sim&  PGN(1-\mu , k).\\
    \sqrt{V}	 \left( 2\cos^2(\frac{\pi U}{2}) -1   \right)		&\sim&  PGN(\mu ,k).\\
   \sqrt{V}	 \left(1- 2\sin^2(\frac{\pi U}{2})   \right)		&\sim&  PGN(\mu ,k).\\
    \sqrt{V}	 \left( \cos^2(\frac{\pi U}{2}) -\sin^2(\frac{\pi U}{2})   \right)		&\sim&   PGN(\mu ,k).\\
    \sqrt{V}	 \frac{1}{2}\left( e^{i\pi U}+e^{-i\pi U}  \right)		&\sim&  PGN(\mu ,k).
\end{eqnarray*}

}

\proof{
The proof is obvious from both Definition (\ref{definitionpolar})  and
 \begin{eqnarray*}
\cos(\pi U)&=&\frac{1}{2} \left( e^{i\pi U}+e^{-i\pi U}  \right)
\\ &=&
 2\cos^2(\frac{\pi U}{2}) -1
 \\ &=&
1- 2\sin^2(\frac{\pi U}{2})
 \\ &=&
 \cos^2(\frac{\pi U}{2}) -\sin^2(\frac{\pi U}{2}). \quad \blacksquare
\end{eqnarray*}
}

\subsection{Moments derivation}

In this section, some results on the
moment properties of $PGN(\mu ,k)$ distribution
are obtained.

\lemma{
\label{theoremIm}
Let $U$ has a beta distribution with parameters
$a$ and $b$. Define  $I_{(m)}(a,b)=E(\cos^m( \pi U))$, so

\small{
\begin{eqnarray}
 \label{eq39ii1}
I_{(1)}(a , b )&=&
 \frac{1}{2} \left( {}_1F_1(a;a+b;i\pi) +  {}_1F_1(a;a+b;-i\pi) \right), \\
\label{eq40ii2}
I_{(2)}(a , b )&=&
  \frac{1}{4} \left(2+
 {}_1F_1(a;a+b;i2\pi) +  {}_1F_1(a;a+b;-i2\pi)
 \right),
 \\
\label{eq41ii3}
I_{(3)}(a , b )&=&
 \frac{3}{8} \bigg(
  {}_1F_1(a;a+b;i\pi) +  {}_1F_1(a;a+b;-i\pi)  \bigg) +
 \frac{1}{8} \bigg(
 {}_1F_1(a;a+b;i3\pi) +  {}_1F_1(a;a+b;-i3\pi) \bigg),
  \\
\label{eq42ii4}
I_{(4)}(a , b )&=&
\frac{1}{16} \bigg(6+4
 {}_1F_1(a;a+b;i2\pi) + 4 {}_1F_1(a;a+b;-i2\pi)
 +
{}_1F_1(a;a+b;i4\pi) +  {}_1F_1(a;a+b;-i4\pi)
\bigg),
\end{eqnarray}
}
 where
  \begin{eqnarray}
 \label{eqhypergeoan}
(a)_{(n)}=\frac{\Gamma (a+n)}{\Gamma(a)}=\prod_{k=0}^{n-1} (a+k),
 \end{eqnarray}
 and
 ${}_1F_1(a;b;z)$ is the Kummer's (confluent hypergeometric) functions, that is,

 \begin{eqnarray}
 \label{eq1f1}
{}_1F_1(a;b;z)=\sum_{n=0}^{\infty} \frac{(a)_{(n)}}{(b)_{(n)}} \frac{z^n}{n!}
 = \int_0^{1} \frac{1}{B(a,b)} e^{zu} u^{a-1} (1-u)^{b-a-1} du.
 \end{eqnarray}

}

\proof{
Define   $J_{(m)}(a,b)=E(\cos(m \pi U))$ so
 \begin{eqnarray*}
J_{(m)}(a,b)&=&\frac{1}{B(a,b)}\int_{0}^{1} \cos (m\pi u) u^{a-1} (1-u)^{b-1}  du
\nonumber \\ &=&
\frac{1}{B(a,b)}\int_{0}^{1}  \frac{1}{2}(e^{i m \pi u}+e^{-i m \pi u}) u^{a-1} (1-u)^{b-1}  du
\nonumber \\ &=& \frac{1}{2} ( CF(i m \pi) +CF(i m \pi))
\nonumber \\ &=& \frac{1}{2} ( {}_1F_1(a;a+b;i m \pi)+{}_1F_1(a;a+b;-i m \pi)),
 \end{eqnarray*}
 where $CF(\cdot )$ is the characteristic function of the beta distribution. According to

\begin{eqnarray*}
\cos^2(x)&=&\frac{1}{2}(1+\cos(2x)),\\
\cos^3(x)&=&\frac{1}{4} (3cos(x)+cos(3x)), \\
\cos^4(x)&=&\frac{1}{8} (3+4cos(2x)+cos(4x)),
\end{eqnarray*}

we conclude that

\begin{eqnarray*}
I_{(1)}(a , b )&=&J_{(1)}(a , b )
= \frac{1}{2} \left( {}_1F_1(a;a+b;i\pi) +  {}_1F_1(a;a+b;-i\pi) \right), \\
I_{(2)}(a , b )&=&\frac{1}{2} \left(1+ J_{(2)}(a , b ) \right)
 \\
&=&\frac{1}{2} \left(1+
 \frac{1}{2}{}_1F_1(a;a+b;i2\pi) +  \frac{1}{2} {}_1F_1(a;a+b;-i2\pi)
 \right)
  \\
&=& \frac{1}{4} \left(2+
 {}_1F_1(a;a+b;i2\pi) +  {}_1F_1(a;a+b;-i2\pi)
 \right),
 \\
I_{(3)}(a , b )&=&\frac{3}{4} J_{(1)}(a , b )+\frac{1}{4}
 J_{(3)}(a , b )
 \\
&=& \frac{3}{8} \left(
  {}_1F_1(a;a+b;i\pi) +  {}_1F_1(a;a+b;-i\pi)  \right) +
 \frac{1}{8} \left(
 {}_1F_1(a;a+b;i3\pi) +  {}_1F_1(a;a+b;-i3\pi) \right),
  \\
I_{(4)}(a , b )&=&
\frac{1}{8} \left(3+4 J_{(2)}(a , b )+ J_{(4)}(a , b)  \right)  \nonumber \\
&=&
\frac{1}{8} \bigg(3+4
 \frac{1}{2} \left\{
{}_1F_1(a;a+b;i2\pi) +  {}_1F_1(a;a+b;-i2\pi) \right\}
+
 \frac{1}{2} \left\{
{}_1F_1(a;a+b;i4\pi) +  {}_1F_1(a;a+b;-i4\pi) \right\}
\bigg)
 \\
&=&
\frac{1}{16} \bigg(6+4
 {}_1F_1(a;a+b;i2\pi) + 4 {}_1F_1(a;a+b;-i2\pi)
 +
{}_1F_1(a;a+b;i4\pi) +  {}_1F_1(a;a+b;-i4\pi)
\bigg). \quad \blacksquare
%\end{aligned}
\end{eqnarray*}

}

We will now determine the moments of $PGN(\mu,k)$ using Lemma \ref{theoremIm}.

\theorem{
Let $X\sim PGN(\mu,k)$ , we can deduce the following properties:
\small{
\begin{eqnarray}
\label{eqex}
  E(X)&=&
\frac{\Gamma(\frac{k+1}{2})}{\sqrt{2}\Gamma(\frac{k}{2})}  \left( {}_1F_1(2\mu ;2;i\pi) +  {}_1F_1(2\mu ;2; -i\pi) \right),  \\
\label{eqex2}
   E(X^2)&=&  \frac{k}{4} \left(2+ {}_1F_1(2\mu ;2;i2\pi) +  {}_1F_1(2\mu ;2;-i2\pi) \right), \\
   \label{eqex3}
  E(X^3) &=& (k +1)   \frac{\sqrt{2}\Gamma(\frac{k+1}{2})}{8\Gamma(\frac{k}{2})}
\bigg\{3   {}_1F_1(2\mu;2;i\pi) + 3 {}_1F_1(2\mu;2;-i\pi)  \big)
 +   {}_1F_1(2\mu ;2;i3\pi) +  {}_1F_1(2\mu ;2;-i3\pi)  \bigg\}, \\
  E(X^4) &=& \frac{ k (k +2)}{16} \bigg(6+4
 {}_1F_1(2\mu ;2;i2\pi) + 4 {}_1F_1(2\mu ;2;-i2\pi)
+
{}_1F_1(2\mu ;2;i4\pi) +  {}_1F_1(2\mu ;2;-i4\pi)
\bigg).
\end{eqnarray}
}

}

\proof{
Since
 $V$ and $U$ are independent, we have
\begin{eqnarray*}
E(X^m)=E(\sqrt{V}^m) E(\cos^m(U))=\frac{\Gamma(\frac{k+m}{2})}{\Gamma(\frac{k}{2})}2^{\frac{m}{2}} E(\cos^m(\pi U))=\frac{\Gamma(\frac{k+m}{2})}{\Gamma(\frac{k}{2})}2^{\frac{m}{2}} I_{m}(2\mu,2(1-\mu)).
\end{eqnarray*}

\begin{eqnarray*}
E(X)&=&
\frac{\Gamma(\frac{k+1}{2})}{\Gamma(\frac{k}{2})} \sqrt{2} \frac{1}{2}  \left( {}_1F_1(2\mu ;2;i\pi) +  {}_1F_1(2\mu ;2; -i\pi) \right)
\\ &=&
\frac{\Gamma(\frac{k+1}{2})}{\sqrt{2}\Gamma(\frac{k}{2})}  \left( {}_1F_1(2\mu ;2;i\pi) +  {}_1F_1(2\mu ;2; -i\pi) \right).
\end{eqnarray*}

\begin{eqnarray*}
E(X^2)&=&E(\sqrt{V}^2) E(\cos^2(\pi U))= \frac{\Gamma(\frac{k+2}{2})}{\Gamma(\frac{k}{2})}2^{\frac{2}{2}} E(\cos^2(\pi U))
 \\ &=&
k I_{(2)}(2\mu,2(1-\mu))
 \\ &=&
k  \frac{1}{4} \left(2+ {}_1F_1(\mu ;1;i2\pi) +  {}_1F_1(\mu;2;-i2\pi) \right)
 \\ &=&
  \frac{k}{4} \left(2+ {}_1F_1(2\mu ;2;i2\pi) +  {}_1F_1(2\mu ;2;-i2\pi) \right).
\end{eqnarray*}

\small{
\begin{eqnarray*}
\begin{aligned}
E(X^3)&=E(\sqrt{V}^3) E(\cos^3(\pi U))=  \frac{\Gamma(\frac{k+3}{2})}{\Gamma(\frac{k}{2})}2^{\frac{3}{2}} E(\cos^k(\pi U))
=
(k +1)   \frac{\Gamma(\frac{k+1}{2})}{\Gamma(\frac{k}{2})} \sqrt{2}  I_{(3)}(2\mu,2(1-\mu))
\nonumber \\
&=
(k +1)   \frac{\Gamma(\frac{k+1}{2})}{\Gamma(\frac{k}{2})} \sqrt{2}
\bigg\{ \frac{3}{8} \big(
  {}_1F_1(2\mu ; 2;i\pi) +  {}_1F_1(2\mu ; 2;-i\pi)  \big)
 +  \frac{1}{8} \big(
 {}_1F_1(2\mu ; 2;i3\pi) +  {}_1F_1(2\mu ; 2;-i3\pi) \big) \bigg\}
\nonumber \\ &=
(k +1)   \frac{\Gamma(\frac{k+1}{2})}{\Gamma(\frac{k}{2})}
\frac{\sqrt{2}}{8}
\bigg\{3   {}_1F_1(2\mu ; 2;i\pi) + 3 {}_1F_1(2\mu ; 2;-i\pi)  \big)
 +   {}_1F_1(2\mu ; 2;i3\pi) +  {}_1F_1(2\mu ; 2;-i3\pi)  \bigg\}
\nonumber \\ &=
(k +1)   \frac{\Gamma(\frac{k+1}{2})}{\Gamma(\frac{k}{2})}
\frac{\sqrt{2}}{8}
\bigg\{3   {}_1F_1(2\mu ;2;i\pi) + 3 {}_1F_1(2\mu ;2;-i\pi)  \big)
 +   {}_1F_1(2\mu ;2;i3\pi) +  {}_1F_1(2\mu ;2;-i3\pi)  \bigg\}
 \\ &=
 (k +1)   \frac{\sqrt{2}\Gamma(\frac{k+1}{2})}{8\Gamma(\frac{k}{2})}
\bigg\{3   {}_1F_1(2\mu ;2;i\pi) + 3 {}_1F_1(2\mu ;2;-i\pi)  \big)
+    {}_1F_1(2\mu ;2;i3\pi) +  {}_1F_1(2\mu ;2;-i3\pi)  \bigg\}.
  \end{aligned}
  \end{eqnarray*}
  }

\small{
\begin{eqnarray*}
E(X^4)&=&E(\sqrt{V}^4) E(\cos^4(\pi U))=
\frac{\Gamma(\frac{k+4}{2})}{\Gamma(\frac{k}{2})}2^{\frac{4}{2}}
 E(\cos^4(\pi U))
=
k (k +2)  I_{(4)}(2\mu,2(1-\mu))
 \nonumber \\ &=&
 k (k +2)
\frac{1}{16} \bigg(6+4
 {}_1F_1(2 \mu ; 2;i2\pi) + 4 {}_1F_1(2 \mu ; 2;-i2\pi)
+
{}_1F_1(2 \mu ; 2;i4\pi) +  {}_1F_1(2 \mu ; 2;-i4\pi)
\bigg)
\nonumber \\ &=&
\frac{ k (k +2)}{16} \bigg(6+4
 {}_1F_1(2 \mu ; 2;i2\pi) + 4 {}_1F_1(2 \mu ; 2;-i2\pi)
+
{}_1F_1(2 \mu ; 2;i4\pi) +  {}_1F_1(2 \mu ; 2;-i4\pi)
\bigg). \quad \blacksquare
\end{eqnarray*}
}

}

\corollary{ The variance of the $PGN(\mu, k)$ is
\small{
\begin{eqnarray}
\label{eqvarx}
  Var(X)&=&\frac{k}{4} \left(2+ {}_1F_1(2\mu ;2;i2\pi) +  {}_1F_1(2\mu ;2;-i2\pi) \right)
-\left(
\frac{\Gamma(\frac{k+1}{2})}{\sqrt{2}\Gamma(\frac{k}{2})}  \left( {}_1F_1(2\mu ;2;i\pi) +  {}_1F_1(2\mu ;2; -i\pi) \right)
 \right)^2.
\end{eqnarray}
}
}

\theorem{If $X\sim PGN(\mu , k)$, the moment of order $m$ is given by
\small{
\begin{eqnarray}
\label{exm1}
E(X^m)&=& \frac{\Gamma(\frac{k+m}{2})}{\Gamma(\frac{k}{2})}2^{\frac{m}{2}} I_{m}(2\mu,2(1-\mu))
\\
\label{exm2} &=&
\frac{\Gamma(\frac{k+m}{2})}{\Gamma(\frac{k}{2})}2^{\frac{k-2m}{2}}
\times
\left\{
\begin{array}{cc}
	\left(  \sum_{r=0}^{2r<m}  { m \choose 2r }
 ( {}_1F_1(2 \mu ; 2;i  (m-2r)\pi)+{}_1F_1(2 \mu ; 2;-i  (m-2r)\pi))
  \right)		+		 {m \choose \frac{m}{2}}
&	m=2,4,\cdots , \\
	\left(  \sum_{r=0}^{2r<m}  { m \choose 2r }
 ( {}_1F_1(2 \mu ; 2;i  (m-2r)\pi)+{}_1F_1(2 \mu ; 2;-i  (m-2r)\pi))
      \right)	
& m=1,3,\cdots .  \\
\end{array}
\right. 
\end{eqnarray}
}

}

\proof{
With regard to

\begin{eqnarray*}
\cos^m(x)=
\left\{
\begin{array}{cc}
\frac{1}{2^{m-1}}	\left(  \sum_{r=0}^{2r<m}  { m \choose 2r }  \cos (m-2r)x        \right)		+		\frac{1}{2^{m}} {m \choose \frac{m}{2}}
&	m=2,4,\cdots .  \\
\frac{1}{2^{m-1}}	\left(  \sum_{r=0}^{2r<m}  { m \choose 2r }  \cos (m-2r)x        \right)	
& m=1,3,\cdots , \\
\end{array}
\right.
\end{eqnarray*}
we obtain
\begin{eqnarray*}
I_{(m)}(a,b)&=&\int_0^{1} \frac{1}{B(a,b)} \cos^m(\pi u) u^{a-1} (1-u)^{b-1} du
\\
&=&
\left\{
\begin{array}{cc}
\frac{1}{2^{m-1}}	\left(  \sum_{r=0}^{2r<m}  { m \choose 2r }
\int_0^{1} \frac{1}{B(a,b)}  \cos\big( (m-2r)\pi u\big)  u^{a-1} (1-u)^{b-1} du
  \right)		+		\frac{1}{2^{m}} {m \choose \frac{m}{2}}
&	m=2,4,\cdots  \\
\frac{1}{2^{m-1}}	\left(  \sum_{r=0}^{2r<m}  { m \choose 2r }
\int_0^{1} \frac{1}{B(a,b)}  \cos\big( (m-2r)\pi u\big)  u^{a-1} (1-u)^{b-1} du
      \right)	
& m=1,3,\cdots  \\
\end{array}
\right.
\\
&=&
\left\{
\begin{array}{cc}
\frac{1}{2^{m-1}}	\left(  \sum_{r=0}^{2r<m}  { m \choose 2r }
 \frac{1}{2} ( {}_1F_1(a;a+b;i  (m-2r)\pi)+{}_1F_1(a;a+b;-i  (m-2r)\pi))
  \right)		+		\frac{1}{2^{m}} {k \choose \frac{m}{2}}
&	m=2,4,\cdots  \\
\frac{1}{2^{m-1}}	\left(  \sum_{r=0}^{2r<m}  { m \choose 2r }
 \frac{1}{2} ( {}_1F_1(a;a+b;i  (m-2r)\pi)+{}_1F_1(a;a+b;-i  (m-2r)\pi))
      \right)	
& m=1,3,\cdots  \\
\end{array}
\right.
\\
&=&
\left\{
\begin{array}{cc}
\frac{1}{2^{m}}	\left(  \sum_{r=0}^{2r<m}  { m \choose 2r }
 ( {}_1F_1(a;a+b;i  (m-2r)\pi)+{}_1F_1(a;a+b;-i  (m-2r)\pi))
  \right)		+		\frac{1}{2^{m}} {m \choose \frac{m}{2}}
&	m=2,4,\cdots ,  \\
\frac{1}{2^{m}}	\left(  \sum_{r=0}^{2r<m}  { m \choose 2r }
 ( {}_1F_1(a;a+b;i  (m-2r)\pi)+{}_1F_1(a;a+b;-i  (m-2r)\pi))
      \right)	
& m=1,3,\cdots . \\
\end{array}
\right.
\end{eqnarray*}
Thus,
\small{
\begin{eqnarray*}
E(X^m)&=& E(\sqrt{V}^m)  E(\cos^m(\pi U)) = E(\sqrt{V}^m)  E(\cos^m(\pi U))
\\
&=&\frac{\Gamma(\frac{k+m}{2})}{\Gamma(\frac{k}{2})}2^{\frac{k}{2}} I_{(m)}(a,b)
\\
&=&
\frac{\Gamma(\frac{k+m}{2})}{\Gamma(\frac{k}{2})}2^{\frac{k}{2}}
\times
\left\{
\begin{array}{cc}
\frac{1}{2^{m}}	\left(  \sum_{r=0}^{2r<m}  { m \choose 2r }
 ( {}_1F_1(a;a+b;i  (m-2r)\pi)+{}_1F_1(a;a+b;-i  (m-2r)\pi))
  \right)		+		\frac{1}{2^{m}} {m \choose \frac{m}{2}}
&	m=2,4,\cdots  \\
\frac{1}{2^{m}}	\left(  \sum_{r=0}^{2r<m}  { m \choose 2r }
 ( {}_1F_1(a;a+b;i  (m-2r)\pi)+{}_1F_1(a;a+b;-i  (m-2r)\pi))
      \right)	
& m=1,3,\cdots  \\
\end{array}
\right.
\\
&=&
\frac{\Gamma(\frac{k+m}{2})}{\Gamma(\frac{k}{2})}2^{\frac{k-2m}{2}}
\times
\left\{
\begin{array}{cc}
	\left(  \sum_{r=0}^{2r<m}  { m \choose 2r }
 ( {}_1F_1(a;a+b;i  (m-2r)\pi)+{}_1F_1(a;a+b;-i  (m-2r)\pi))
  \right)		+		 {m \choose \frac{m}{2}}
&	m=2,4,\cdots , \\
	\left(  \sum_{r=0}^{2r<m}  { m \choose 2r }
 ( {}_1F_1(a;a+b;i  (m-2r)\pi)+{}_1F_1(a;a+b;-i  (m-2r)\pi))
      \right)	
& m=1,3,\cdots , \\
\end{array}
\right.
\end{eqnarray*}
}
the result is proved by replacing $a=2\mu$ and $b=2(1-\mu)$. $\quad \blacksquare$
}

From the Equation (\ref{exm1}), $E(X^m)$  can be easily calculated.
 Table \ref{Itable} in Appendix  gives the value of
 $I_m(2\mu,2(1-\mu))$ under different values of $\mu$ and $m$.
 The hypergeo package (Hankin, 2015) of R  routines  is used for numerical calculation of
hypergeometric functions.  For example, under $\mu=0.3$
and $k=5$,
\begin{eqnarray*}
  E(X) &=& \frac{\Gamma(\frac{k+1}{2})}{\Gamma(\frac{k}{2})} \sqrt{2} I_1(2\mu,2(1-\mu))= 2.1277 \times  0.46=0.9787, \\
  E(X^2)&=& \frac{\Gamma(\frac{k+2}{2})}{\Gamma(\frac{k}{2})}2^{\frac{2}{2}} EI_2(2\mu,2(1-\mu)) = 5\times 0.5995=2.9975,\\
  E(X^3)&=&  \frac{\Gamma(\frac{k+3}{2})}{\Gamma(\frac{k}{2})}2^{\frac{3}{2}} I_3(2\mu,2(1-\mu))= 12.7661\times  0.3942=5.0324,  \\
  E(X^4)&=&
\frac{\Gamma(\frac{k+4}{2})}{\Gamma(\frac{k}{2})}2^{\frac{4}{2}}
 I_4(2\mu,2(1-\mu))=35\times 0.4916=17.2060.
\end{eqnarray*}

\section{Bayesian analysis}
\label{sectionbayes}
Here we specify a Bayesian model that accounts for asymmetry and bimodality. In what follows,
the posterior inference from this model is discussed.
Let $X_1, X_2, \cdots, X_n$ be $n$ independent and identically distributed random variables from $PGN(\mu , k)$.
The Bayesian model specification requires a prior distributions for all the unknown
parameters, i.e., $\theta=(\mu , k)$.
 In the absence of prior information
and in order to guarantee posterior property, we adopt proper but diffuse priors. For
convenience but not always optimal, we suppose that elements of $\theta$ are independent
 so the joint prior distribution
$\pi (\theta)=\pi (\mu , k)=\pi (\mu ) \pi (k )$.
we adopt the following prior distributions
\begin{eqnarray}
\label{priormu}
\mu &\sim & Uniform(0,1), \\
\label{priork}
   k &\sim & \Gamma (k_0 , k_1).
\end{eqnarray}
Here, the hyper-parameters $k_0$ and $k_1$ are known (to be specified) positive scalers. We choose the non-informative uniform prior for parameter $\mu$ to allow the data to select appropriate value for $\mu$.
Then the posterior distribution is proportional
to
\begin{eqnarray}
\pi(\theta|\mathbf{x})&\propto & \pi(\mathbf{x} |\theta ) \pi(\theta)\\
&\propto &  (\prod_{i=1}^{n}p(x_i |\mu , k))\pi (\mu) \pi (k),
\end{eqnarray}
where $\mathbf{x}=(x_1,x_2, \cdots, x_n)$  and  $p$ denotes  the density function of $PGN(\mu , k)$ given in  (\ref{eqpdf}). The joint posterior distribution of all unknown quantities involved is given by
\begin{eqnarray}
\pi(\theta, \mathbf{u}|\mathbf{x})&\propto &  \prod_{i=1}^{n}p(x_i|u_i , k)p(u_i|\mu)\pi (\mu) \pi (k),
\end{eqnarray}
where $\mathbf{u}=(u_1,u_2, \cdots, u_n)$ and \begin{eqnarray}
\label{fullxi}
p(x_i|u_i , k)  &= &
\left\{
\begin{array}{ccc}
\frac{1}{|\cos (\pi u_i)|2^{\frac{k}{2}-1} \Gamma(\frac{k}{2})}( \frac{x_i}{\cos (\pi u_i)})^{\frac{k}{2}-1}  e^{- \frac{x_i^2}{2\cos^2 (\pi u_i)} }
& u_i>\frac{1}{2} & x_i<0\\
\frac{1}{|\cos (\pi u_i)|2^{\frac{k}{2}-1} \Gamma(\frac{k}{2})}( \frac{x_i}{\cos (\pi u_i)})^{\frac{k}{2}-1}  e^{- \frac{x_i^2}{2\cos^2 (\pi u_i)} } & u_i<\frac{1}{2} & x>0 \\
0& $O.W.$
\end{array}
\right.
 \nonumber \\
&= &
 \frac{1}{|\cos (\pi u_i)|2^{\frac{k}{2}-1} \Gamma(\frac{k}{2})}(| \frac{x_i}{\cos (\pi u_i)}|)^{\frac{k}{2}-1}  e^{- \frac{x_i^2}{2\cos^2 (\pi u_i)} } 1_{\{x_i\times(\frac{1}{2}-u_i)\}}.
\end{eqnarray}
The posterior distribution is in a complicated form and we will use a Markov chain Monte Carlo (MCMC) method to generate samples from posterior distributions. To generate samples from the posterior distribution, we exploit the full conditional distributions in a Gibbs sampling framework.
The full conditional $p(u_i|x_i, \mu, k) $ does not defines a standard probability distribution, so sampling from this distribution is not simply practicable. Some methods like rejection sampling, importance sampling or Metropolis-Hastings sampler can be used for sampling from this full conditional. We  use Metropolis-Hastings sampler with generate proposals of beta distribution that the mean of it, centered on the current value. That is, in the Gibbs  sampling algorithm, at step  $t$, $u_i^{(t+1)}$ is generated from $Beta(2u_i^{(t)},2(1-u_i^{(t)}))$.

The full conditional posterior distribution  of the parameter $k$, i. e. $\pi (k |  \mathbf{x}, \mathbf{u}, \mu)$,  can be computed as follows:
\begin{eqnarray}
\label{fullk}
\pi (k | \mathbf{x}, \mathbf{u} , \mu) &\propto &
p(x_1|u_1, k)
p(x_2|u_2, k) \cdots p(x_n|u_n, k)  \pi (k) \nonumber \\
&\propto & \pi (k) \prod_{i=1}^{n}
\frac{1}{|\cos (\pi u_i)|^{\frac{k}{2}}2^{\frac{k}{2}-1} \Gamma(\frac{k}{2})}  e^{- \frac{x_i^2}{2\cos^2 (\pi u_i)} }
 (\frac{u_i}{1-u_i})^{2\mu-1}
  1_{\{x_i\times(\frac{1}{2}-u_i)\}}
 \nonumber \\ &\propto &
   \pi (k  ) \prod_{i=1}^{n}
\frac{1}{|2\cos (\pi u_i)|^{\frac{k}{2}} \Gamma(\frac{k}{2})}
  1_{\{x_i\times(\frac{1}{2}-u_i)\}}
   \nonumber \\ &\propto &
   \frac{ \pi (k  )}{2^{\frac{n k}{2}}|\prod_{i=1}^{n} \cos (\pi u_i)|^{\frac{k}{2}} \Gamma^n(\frac{k}{2})}
  \prod_{i=1}^{n}1_{\{x_i\times(\frac{1}{2}-u_i)\}}.
\end{eqnarray}
The full conditional of $k$ does not define a standard probability distribution too.
The Metropolis–Hastings sampler can be used for sampling from this full conditional with generate proposals
of an exponential distribution, centered on the current value. That is, in the Gibbs  sampling algorithm, at step $t$, $k^{(t+1)}$ is generated from $Exp(k^{(t)})$.

The full conditional posterior distribution  of the parameter $k$, $\pi (\mu | \mathbf{x}, \mathbf{u}, k)$,  is given by
\begin{eqnarray}
\label{fullmu}
\pi (\mu | \mathbf{x}, \mathbf{u} , k)
 &\propto &
 p(u_1| \mu )
p(u_2| \mu ) \cdots p(u_n| \mu )  \pi (\mu  )
\nonumber \\ &\propto &
  \frac{\pi (\mu  )}{\Gamma^{n}(2\mu) \Gamma^{n}(2(1-\mu))}(\prod_{i=1}^{n}\frac{u_i}{1-u_i})^{2\mu-1}.
\end{eqnarray}

The full conditional
does not defines a standard probability distribution too. The Metropolis–Hastings sampler can be used for sampling from this full conditional with generate proposals
of beta distribution, centered on the current value.That is, in the Gibbs  sampling algorithm, at step $t$, $\mu^{(t+1)}$ is generated from $Beta(2\mu^{(t)},2(1-\mu^{(t)}))$.

Based on these full conditional distributions, we use the following algorithm for sampling from the
joint posterior distribution:
The main steps of the Gibbs  sampling algorithm
at step $(t + 1)$ as

\begin{enumerate}
  \item  Draw $u^{(t+1)}_i$ from $p(u_i|x_i,\mu^{(t)},k^{(t)})$ for $i=1,2,\cdots, n$.
  \item  Draw $k^{(t+1)}$ from $\pi (k|\mathbf{x},\mathbf{u}^{(t+1)},\mu^{(t)})$.
  \item  Draw $\mu^{(t+1)}$ from $\pi (\mu|\mathbf{x},\mathbf{u}^{(t+1)},k^{(t+1)})$.
  \item  Iterate above steps until we get the appropriate number of MCMC samples.
\end{enumerate}

\section{Simulation study}

\begin{table}[ht]
\centering
\small
\caption{Biases and RMSE  of the parameter estimators of the PGN model  for various values of $\mu$ and $k=1,2,5$.  }
\label{tablesimu1}
\begin{tabular} {|c |c |c |c |c |c |c |c |}
 \hline
  $\mu$ & $\hat{\mu}$ & Bias & RMSE & $k$ & $\hat{k}$ & Bias & RMSE \\ \hline \hline
0.10 & 0.2640 & -0.1640 & 0.0113 & 1 & 1.2980 & -0.2979 & 0.1971 \\ \hline
0.20 & 0.2955 & -0.0955 & 0.0194 & 1 & 1.1630 & -0.1629 & 0.1790 \\ \hline
0.25 & 0.3188 & -0.0688 & 0.0222 & 1 & 1.1100 & -0.1103 & 0.1672 \\ \hline
0.30 & 0.3460 & -0.0460 & 0.0269 & 1 & 1.0760 & -0.0755 & 0.1667 \\ \hline
0.40 & 0.4173 & -0.0173 & 0.0345 & 1 & 1.0350 & -0.0352 & 0.1635 \\ \hline
0.50 & 0.4997 & 0.0003 & 0.0415 & 1 & 1.0070 & -0.0075 & 0.1644 \\ \hline
0.60 & 0.5827 & 0.0173 & 0.0343 & 1 & 1.0370 & -0.0367 & 0.1651 \\ \hline
0.70 & 0.6538 & 0.0462 & 0.0265 & 1 & 1.0740 & -0.0744 & 0.1647 \\ \hline
0.75 & 0.6812 & 0.0688 & 0.0218 & 1 & 1.1100 & -0.1103 & 0.1680 \\ \hline
0.80 & 0.7046 & 0.0954 & 0.0195 & 1 & 1.1610 & -0.1610 & 0.1834 \\ \hline
0.90 & 0.7352 & 0.1648 & 0.0111 & 1 & 1.2980 & -0.2979 & 0.1935 \\ \hline \hline
0.10 & 0.2592 & -0.1592 & 0.0117 & 2 & 2.6750 & -0.6754 & 0.3754 \\ \hline
0.20 & 0.2863 & -0.0863 & 0.0168 & 2 & 2.3430 & -0.3427 & 0.4539 \\ \hline
0.25 & 0.3177 & -0.0678 & 0.0229 & 2 & 2.2910 & -0.2914 & 0.3564 \\ \hline
0.30 & 0.3442 & -0.0442 & 0.0259 & 2 & 2.1930 & -0.1933 & 0.3870 \\ \hline
0.40 & 0.4217 & -0.0217 & 0.0348 & 2 & 2.0600 & -0.0603 & 0.3241 \\ \hline
0.50 & 0.4839 & 0.0161 & 0.0349 & 2 & 1.8600 & 0.1395 & 0.2312 \\ \hline
0.60 & 0.5776 & 0.0224 & 0.0375 & 2 & 2.0210 & -0.0213 & 0.3166 \\ \hline
0.70 & 0.6556 & 0.0444 & 0.0259 & 2 & 2.1980 & -0.1975 & 0.3843 \\ \hline
0.75 & 0.6822 & 0.0678 & 0.0227 & 2 & 2.2930 & -0.2931 & 0.3558 \\ \hline
0.80 & 0.7042 & 0.0958 & 0.0206 & 2 & 2.3590 & -0.3586 & 0.3204 \\ \hline
0.90 & 0.7411 & 0.1589 & 0.0123 & 2 & 2.6160 & -0.6164 & 0.3735 \\ \hline \hline
0.10 & 0.2471 & -0.1471 & 0.0116 & 5 & 6.3180 & -1.3180 & 0.6207 \\ \hline
0.20 & 0.2852 & -0.0852 & 0.0183 & 5 & 5.8600 & -0.8605 & 0.6903 \\ \hline
0.25 & 0.3132 & -0.0632 & 0.0228 & 5 & 5.6490 & -0.6493 & 0.6628 \\ \hline
0.30 & 0.3398 & -0.0398 & 0.0267 & 5 & 5.4490 & -0.4491 & 0.8543 \\ \hline
0.40 & 0.4134 & -0.0134 & 0.0330 & 5 & 5.1830 & -0.1826 & 0.8613 \\ \hline
0.50 & 0.5006 & -0.0006 & 0.0389 & 5 & 5.1970 & -0.1971 & 0.8556 \\ \hline
0.60 & 0.5867 & 0.0133 & 0.0332 & 5 & 5.1880 & -0.1876 & 0.8483 \\ \hline
0.70 & 0.6601 & 0.0399 & 0.0264 & 5 & 5.4460 & -0.4464 & 0.8517 \\ \hline
0.75 & 0.6869 & 0.0631 & 0.0226 & 5 & 5.6480 & -0.6476 & 0.6690 \\ \hline
0.80 & 0.7148 & 0.0852 & 0.0189 & 5 & 5.8550 & -0.8554 & 0.6862 \\ \hline
0.90 & 0.7530 & 0.1470 & 0.0117 & 5 & 6.3280 & -1.3280 & 0.6187 \\ \hline
\end{tabular}
\end{table}

Here we assess the finite sample behaviour of the Bayesian estimators of the parameters in a PGN model using a simulation study. We perform 100 Monte Carlo replications.  In each replication a random sample of size $n=50$ is drawn from $PGN(\mu , k)$. In order to evaluate the accuracy of posterior estimates under different scenarios, we consider various different values for the parameters $\mu$ and $k$.  The true parameter values used in the data generating process are $\mu=0.1,0.2,0.25,0.3,0.4,0.5,0.6,0.7,0.75,0.8,0.9$ and $k=1,2,5,10,15$. The parameters are estimated from the sample by the Bayesian approach.  For a comparison, we
computed two criteria: the average of bias (Bias) and the root of mean-square error
(RMSE). The OpenBUGS is used  for fitting the model. We generate 2,000 MCMC samples after a burn-in 1,000 iterations in one chain with different starting values. We resampled the components every 10 iterations to improve mixing. Visual inspection of the trace plots does not
reveal any obvious mixing or convergence problems, and autocorrelations seem within reasonable levels.

The posterior results are summarized in Tables \ref{tablesimu1} and  \ref{tablesimu2}.
As seen,   the results indicate overwhelming support for the proposed model.
Overall, it is observed that the Bias and RMSE   of the parameter estimates are relatively small.
This illustrates the effectiveness of the estimation procedure outlined in Section 3.
As expected, the smallest Bias's correspond to   $\mu=0.5$. This is a
direct consequence of the fact that we have maximum information when the density is symmetric. There is a small bias in the estimation of the parameters when $\mu$ is close to $0$ or $1$.   It must be noted that  the results have been reported based
on small sample sizes.

\begin{table}[ht]
\centering
\small
\caption{Biases and RMSE  of the parameter estimators of the PGN model  for various values of $\mu$ and $k=10,15$. }
\label{tablesimu2}
\begin{tabular} {|c |c |c |c |c |c |c |c |}
 \hline
  $\mu$ & $\hat{\mu}$ & Bias & RMSE & $k$ & $\hat{k}$ & Bias & RMSE \\ \hline \hline
0.10 & 0.2397 & -0.1397 & 0.0115 & 10 & 11.8100 & -1.8060 & 0.9044 \\ \hline
0.20 & 0.2805 & -0.0805 & 0.0185 & 10 & 11.1200 & -1.1220 & 1.0240 \\ \hline
0.25 & 0.3088 & -0.0588 & 0.0223 & 10 & 10.9100 & -0.9116 & 1.1080 \\ \hline
0.30 & 0.3379 & -0.0379 & 0.0253 & 10 & 10.6100 & -0.6143 & 1.3300 \\ \hline
0.40 & 0.4125 & -0.0125 & 0.0332 & 10 & 10.1200 & -0.1246 & 1.4200 \\ \hline
0.50 & 0.5005 & -0.0005 & 0.0397 & 10 & 10.2000 & -0.1959 & 1.4120 \\ \hline
0.60 & 0.5867 & 0.0133 & 0.0330 & 10 & 10.1400 & -0.1384 & 1.4130 \\ \hline
0.70 & 0.6618 & 0.0382 & 0.0252 & 10 & 10.6000 & -0.6005 & 1.3070 \\ \hline
0.75 & 0.6912 & 0.0588 & 0.0222 & 10 & 10.9100 & -0.9125 & 1.1240 \\ \hline
0.80 & 0.7188 & 0.0812 & 0.0191 & 10 & 11.1700 & -1.1670 & 1.0520 \\ \hline
0.90 & 0.7591 & 0.1409 & 0.0110 & 10 & 11.9000 & -1.9010 & 0.8830 \\ \hline   \hline
0.10 & 0.2351 & -0.1351 & 0.0122 & 15 & 17.3800 & -2.3850 & 1.0470 \\ \hline
0.20 & 0.2773 & -0.0773 & 0.0190 & 15 & 16.3800 & -1.3840 & 1.4030 \\ \hline
0.25 & 0.3070 & -0.0570 & 0.0222 & 15 & 16.2900 & -1.2910 & 1.4550 \\ \hline
0.30 & 0.3374 & -0.0374 & 0.0261 & 15 & 15.8300 & -0.8268 & 1.6610 \\ \hline
0.40 & 0.4132 & -0.0132 & 0.0326 & 15 & 15.1400 & -0.1433 & 1.8900 \\ \hline
0.50 & 0.5009 & -0.0009 & 0.0394 & 15 & 15.2100 & -0.2117 & 1.7480 \\ \hline
0.60 & 0.5870 & 0.0130 & 0.0322 & 15 & 15.1500 & -0.1490 & 1.8940 \\ \hline
0.70 & 0.6627 & 0.0373 & 0.0261 & 15 & 15.8300 & -0.8327 & 1.6950 \\ \hline
0.75 & 0.6929 & 0.0571 & 0.0219 & 15 & 16.2900 & -1.2870 & 1.4540 \\ \hline
0.80 & 0.7223 & 0.0777 & 0.0191 & 15 & 16.4100 & -1.4050 & 1.4350 \\ \hline
0.90 & 0.7648 & 0.1352 & 0.0126 & 15 & 17.3900 & -2.3930 & 1.0530 \\ \hline
\end{tabular}
\end{table}

\section{Illustrative example}
In this section, we apply the location and scale extension of the standard PGN model to five illustrative data sets. This extension, denoted by $PGN(\beta , \sigma^2 , \mu ,k)$, is given by the distribution of
$Y = \beta  + \sigma X$, where $X \sim PGN(\mu ,k)$ and  $\beta \in R$ , $\sigma  > 0$ are the location and scale parameters, respectively.

The data sets  contain clinical annotations and the log expression of $150$ genes for a set of $444$ lung
cancer patients. The  data (lung.dataset) are  available in $oompaData$ package from \href{https://cran.r-project.org/package=oompaData}{https://cran.r-project.org/}.
The original data are also available in FTP from \href{ftp://caftpd.nci.
nih.gov/pub/caARRAY/experiments/caArray_beer-00153/}{here} and in the Gene Expression Omnibus at \href{http://www.ncbi.nlm.nih.gov/geo/query/acc.cgi?acc=GSE68571}{here}. The data were
log transformed by mapping the expression value $x$ to $\log 2(1+x)$. The PGN model is fitted to  the gene $\sharp 37145$ (row 2), the gene $\sharp 202459$ (row 12), the gene $\sharp 208288$ (row 63),  the gene $\sharp 215456$ (row 106) and the gene $\sharp 216437$ (row 113). These data sets display bimodality and asymmetry (see Figure 6).
The parameters are estimated via the Bayesian approach as described in Section 3.
The priors of the parameters  $\beta$ and $\sigma^2$
are specified as:
\begin{eqnarray}
  \label{priorbetac}
  \beta | \sigma^2  &\sim& \text{Normal}(\beta_0 , c\sigma^2) ,\\
  \label{priortauc}
\tau =\frac{1}{\sigma^2} &\sim& \Gamma(\tau_0 ,\tau_1 ).
\end{eqnarray}
Here, the hyper-parameters $c$, $\tau_0$ and $\tau_1$   are known (to be specified) positive scalers and $\beta_0$  is a real-valued scaler. The hyper-parameter $\beta_0$ are set to be close to mean of the data. The hyper-parameters $\tau_0$ and $\tau_1$ are set to be some values such that
to ensure  large prior variance for the parameter $\tau$.
If the histogram of the data be unimodal the hyper-parameters $k_0$ and $k_1$ are chosen in the sense  that $E(k) \approx 2$. For the bimodality cases, these hyperparameters are considered such that $E(K)>2$.

We generate 2,000 MCMC samples after
a burn-in 1,000 iterations in one chain with different starting values. We resampled the components every 10 iterations.
The posterior means and the $95\%$
 highest probability density (HPD) credible intervals of the model parameters are given in Table \ref{tableestimate}.
  The HPD interval is one of the most useful tools to measure the posterior uncertainty.
The short HPD intervals for the parameters, in all data set,
indicate the accuracy of the estimates. Our approach allows for testing symmetry and normality  within the $PGN$  family via the HPD intervals. The central position $\mu=0.5$ allows
 for testing symmetry because under
$\mu=0.5$ the $PGN$ family is symmetric.
 In a similar way,  a normality test
 is $H_0: \mu=0.5 \& \, \, \, k=2$. The HPD intervals are used  to make inferences about hypothesis testing of parameters. Specifically, in the first four data sets, the HPD reject the hypothesis
$\mu=0.5$, so these data sets are not symmetric. However, for the data set 5,
although we reject the normality assumption but  we  can not reject that the distribution is symmetric.

The observed
and estimated densities plotted in Figure 6. As seen, the PGN model is able to capture different distribution shapes which are observed
in the  data sets. This result confirms the great flexibility of the model in accommodating bimodality and asymmetry.

\begin{table}[ht]
\centering
\small
\caption{The  HPD intervals and the estimated values of the parameters under the location and scale PGN model. }
\label{tableestimate}
\begin{tabular} {|c| cc| cc |cc |cc |}
 \hline
 \text{Data set} & $\hat{\beta}$ &  \text{HPD} $ \beta$ & $\hat{\sigma}$ &  \text{HPD} $ \sigma$ &
   $\hat{\mu}$ &  \text{HPD} $ \mu$ & $\hat{k}$ &  \text{HPD} $ k$ \\ \hline \hline
 \text{Data set 1} & $5.7520$ &  $(5.7400 , 5.7720)$ & $1.5655$ &  $(1.3700 , 1.7380)$  &   $0.6380$ &  $(0.6117 , 0.6647
)$  & $2.9210$  & $(2.3000 , 3.6860)$  \\ \hline
 \text{Data set 2} & $8.6950$ &  $(8.6920 , 8.6970
)$  & $1.0490$ &  $(0.9273 , 1.1570)$  &  $0.4367$ &  $(0.4106 , 0.4651)$  & $1.3970$ &  $(1.1870 , 1.6420
)$   \\ \hline
 \text{Data set 3} & $3.5410$ &  $(3.5170 , 3.5470)$  & $0.6692$ &  $( 0.5961 , 0.7790)$  &  $0.3959$ &  $( 0.3686 , 0.4246)$  & $7.6185$ &  $(5.2380 , 9.4610)$   \\ \hline
 \text{Data set 4} & $0.2123$ &  $(0.2086 , 0.2156
)$  & $0.6879$ &  $(0.5605 , 0.7654)$  &  $0.5605$ &  $(0.5323 , 0.5866)$  & $4.3435$  &  $(3.2770 , 6.3830

)$  \\ \hline
 \text{Data set 5} & $0.0553$ &  $(0.0521 , 0.0585)$  & $0.8850$ &  $(0.8114 , 1.0000)$  &  $0.5111$ &  $(0.4817 , 0.5388
)$  & $2.5330$ &  $( 2.0400 , 2.9790)$  \\ \hline
\end{tabular}
\end{table}

\begin{figure}[htp]
  \centering
  \label{figurchain:0}\caption{Markov chain history for   asymmetric and the peak parameters.}

  \subfloat[Asymmetric  parameter (Data set 1). ]{\label{figurchain:1}\includegraphics[width=50mm,height=30mm]{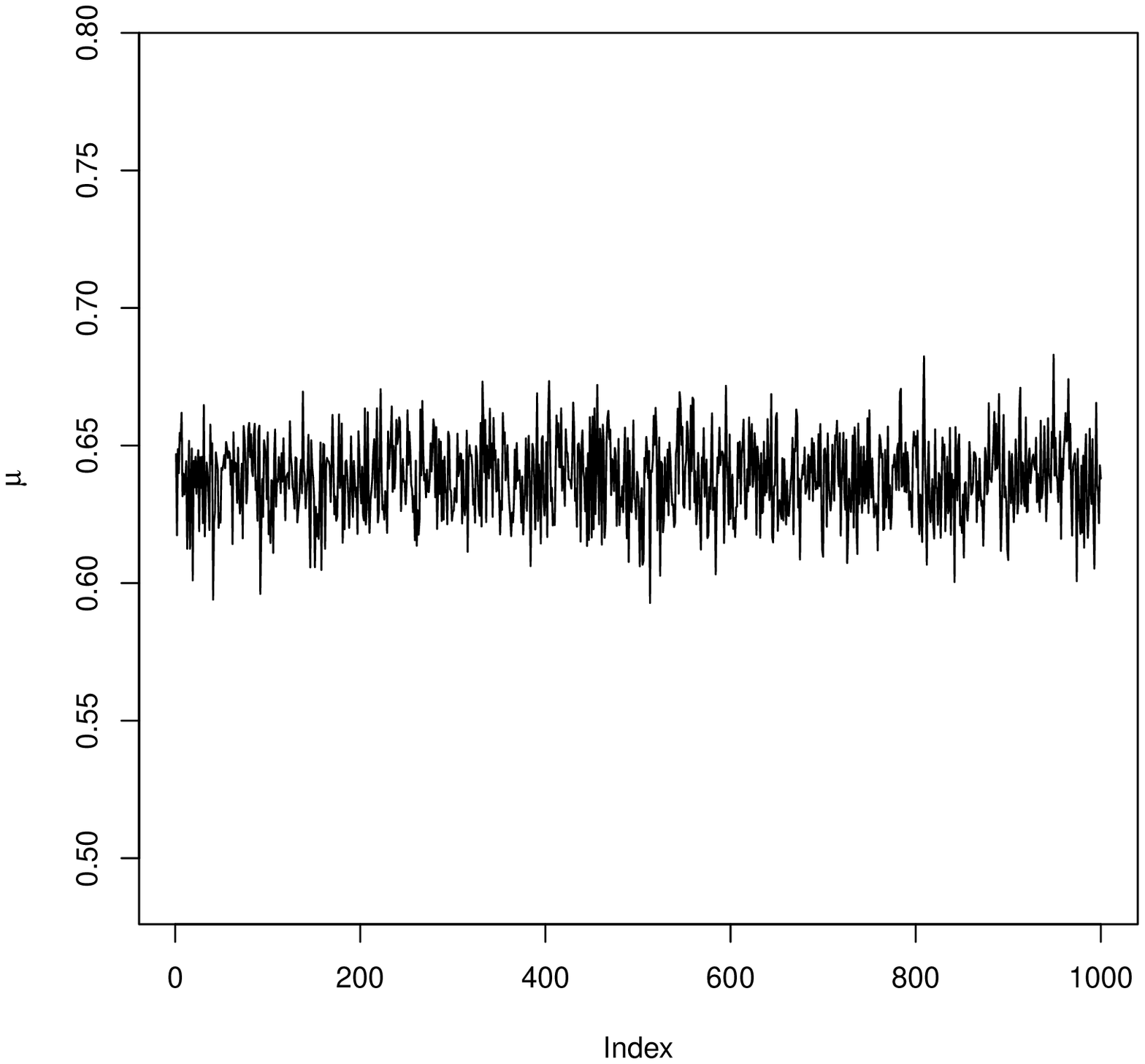}}
    \subfloat[Peak  parameter  (Data set 1).]{\label{figurchain:2}\includegraphics[width=50mm,height=33mm]{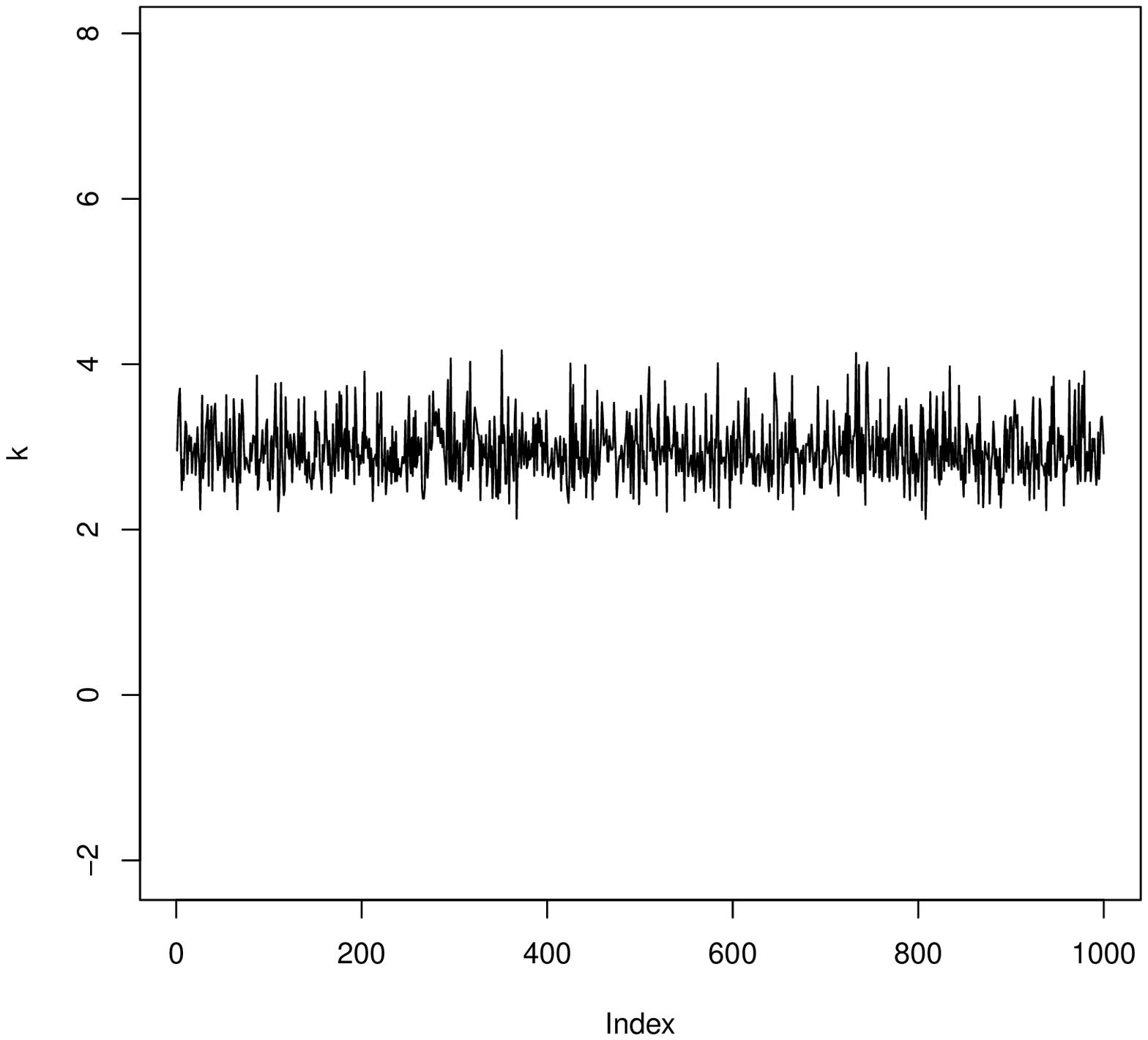}}
  \\
  \subfloat[Asymmetric  parameter (Data set 2). ]{\label{figurchain:3}\includegraphics[width=50mm,height=30mm]{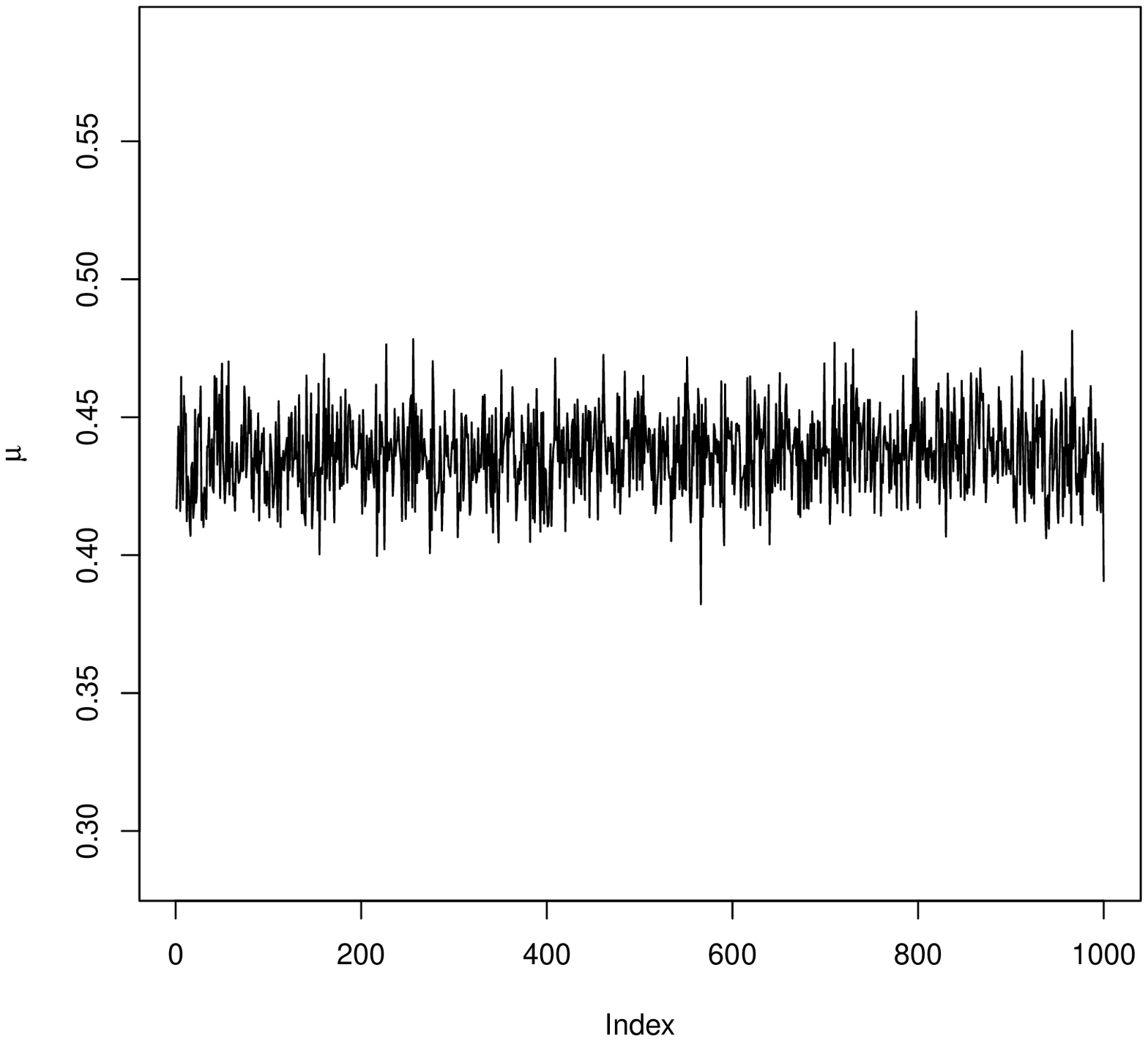}}
    \subfloat[Peak  parameter (Data set 2).]{\label{figurchain:4}\includegraphics[width=50mm,height=30mm]{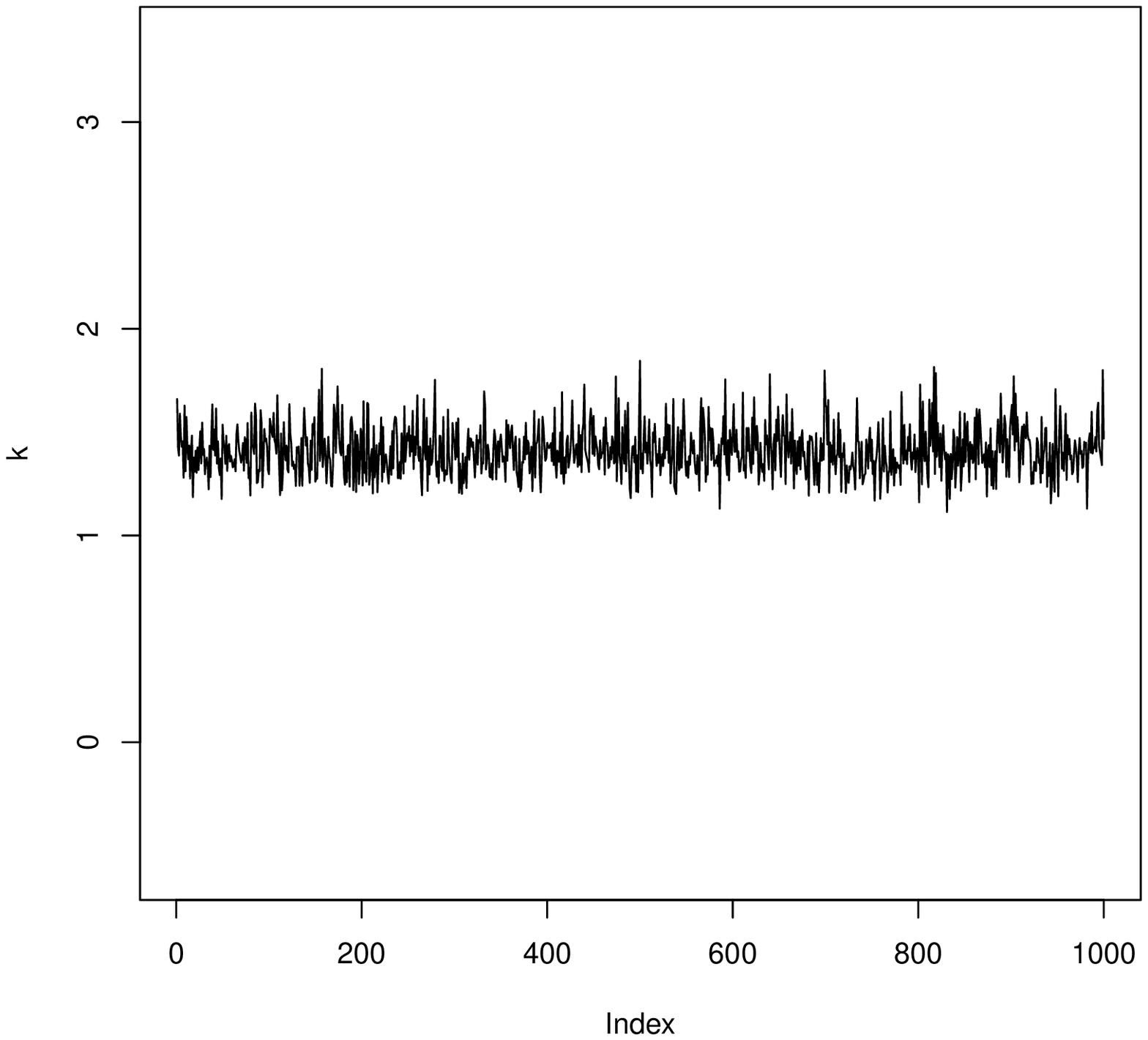}}
  \\
  \subfloat[Asymmetric  parameter (Data set 3). ]{\label{figurchain:5}\includegraphics[width=50mm,height=30mm]{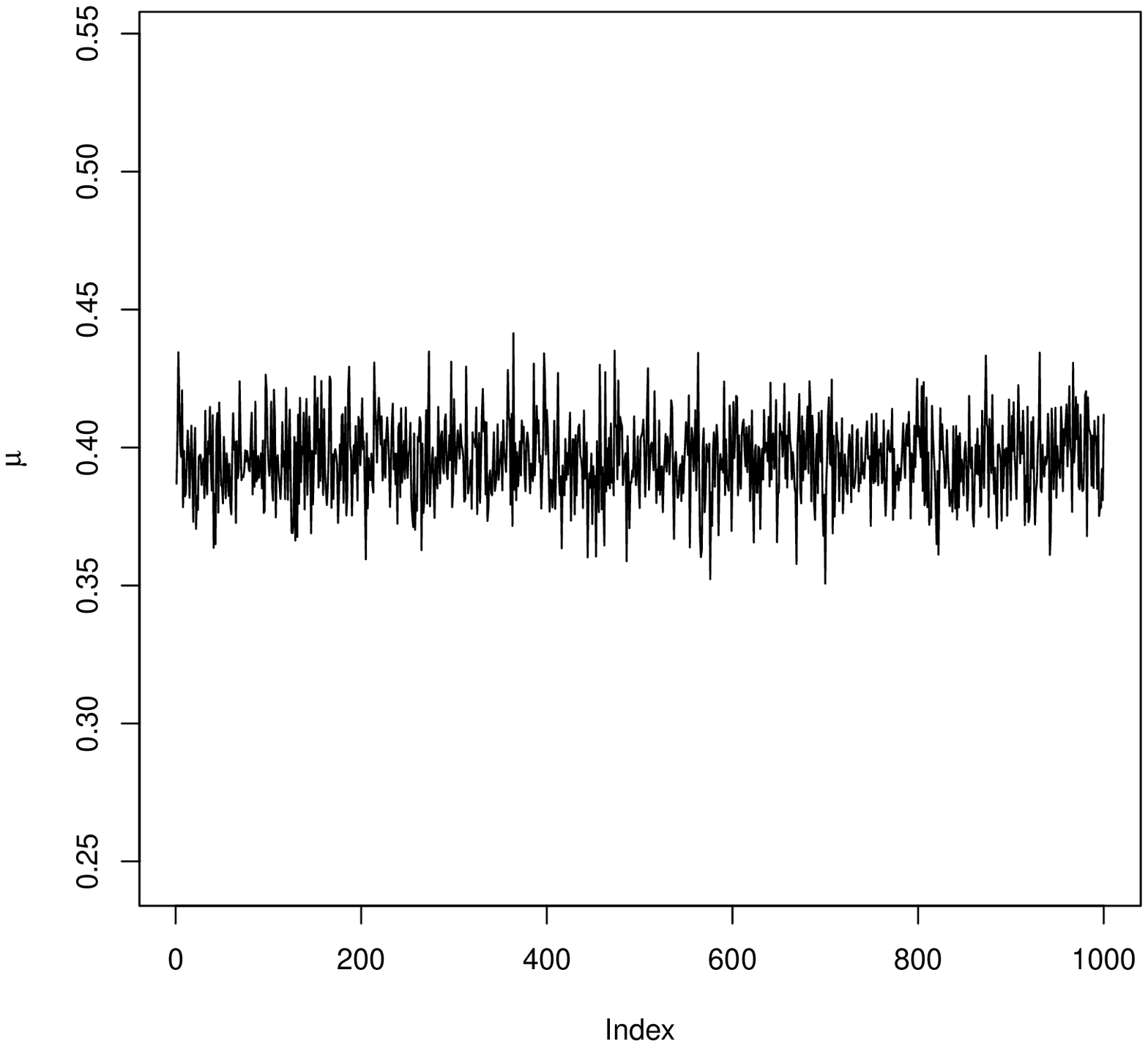}}
    \subfloat[Peak  parameter (Data set 3).]{\label{figurchain:6}\includegraphics[width=50mm,height=30mm]{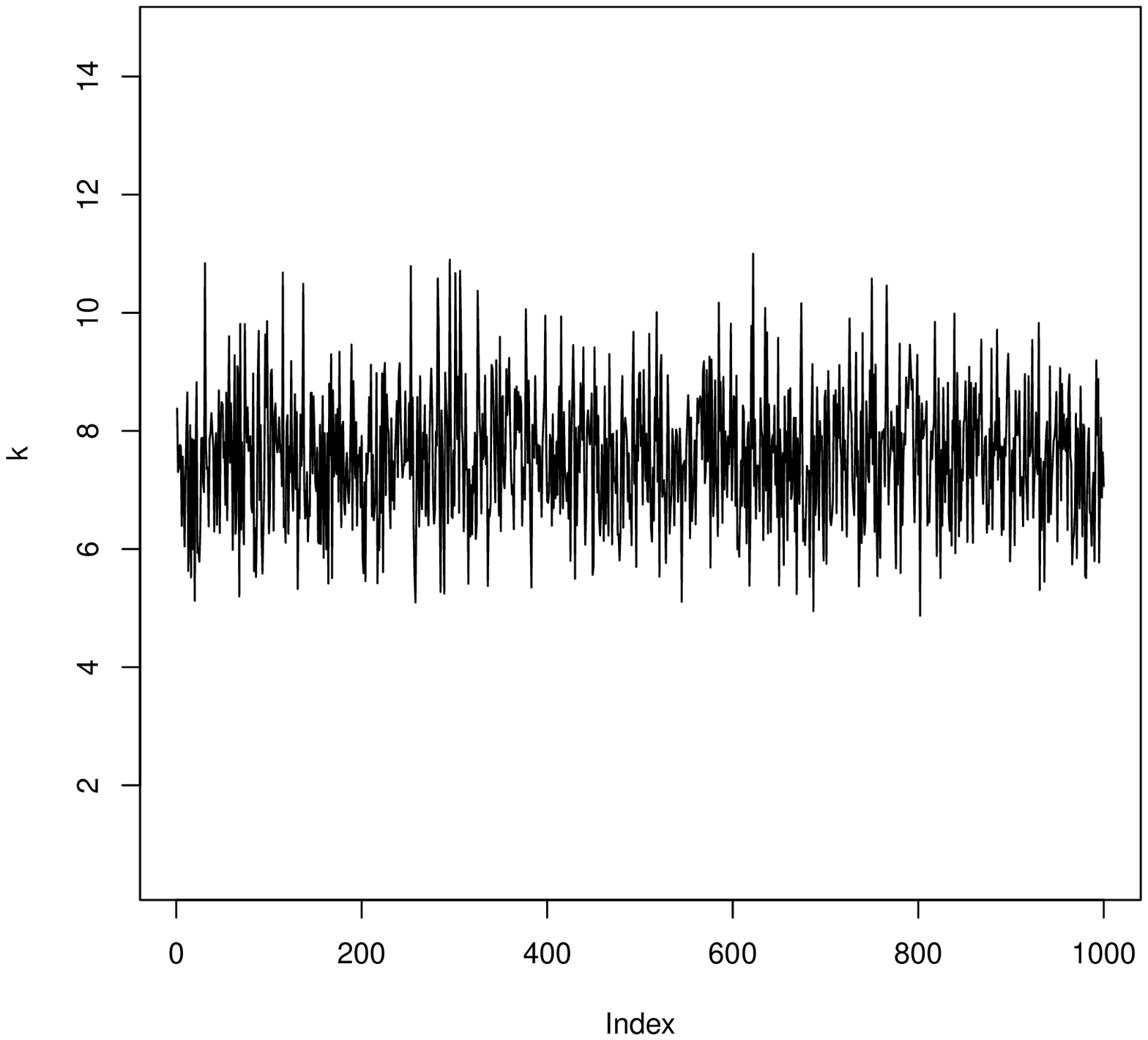}}
    \\
  \subfloat[Asymmetric  parameter (Data set 4) ]{\label{figurchain:7}\includegraphics[width=50mm,height=30mm]{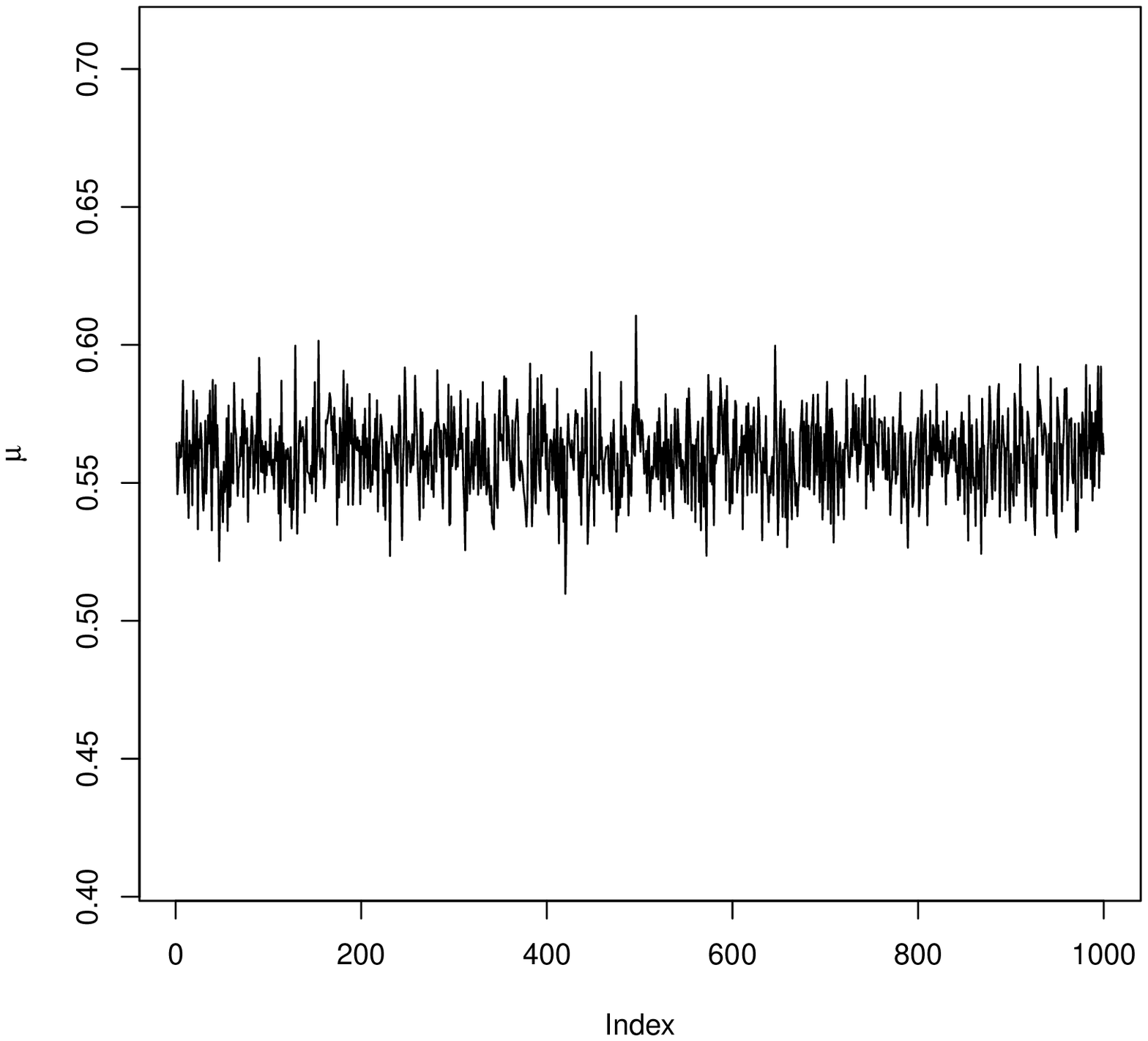}}
    \subfloat[Peak  parameter (Data set 4)]{\label{figurchain:8}\includegraphics[width=50mm,height=30mm]{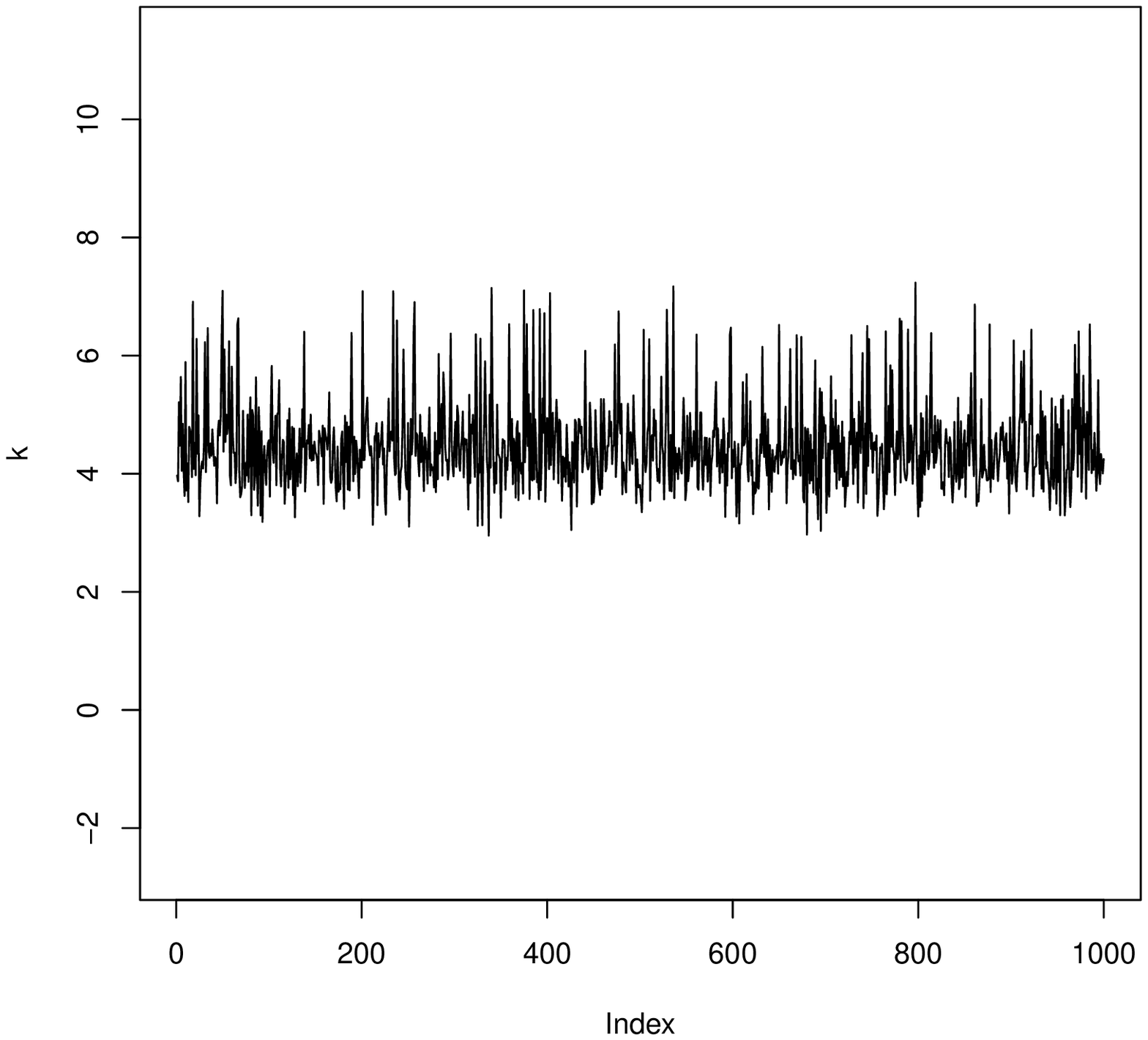}}
    \\
  \subfloat[Asymmetric  parameter (Data set 5). ]{\label{figurchain:9}\includegraphics[width=50mm,height=30mm]{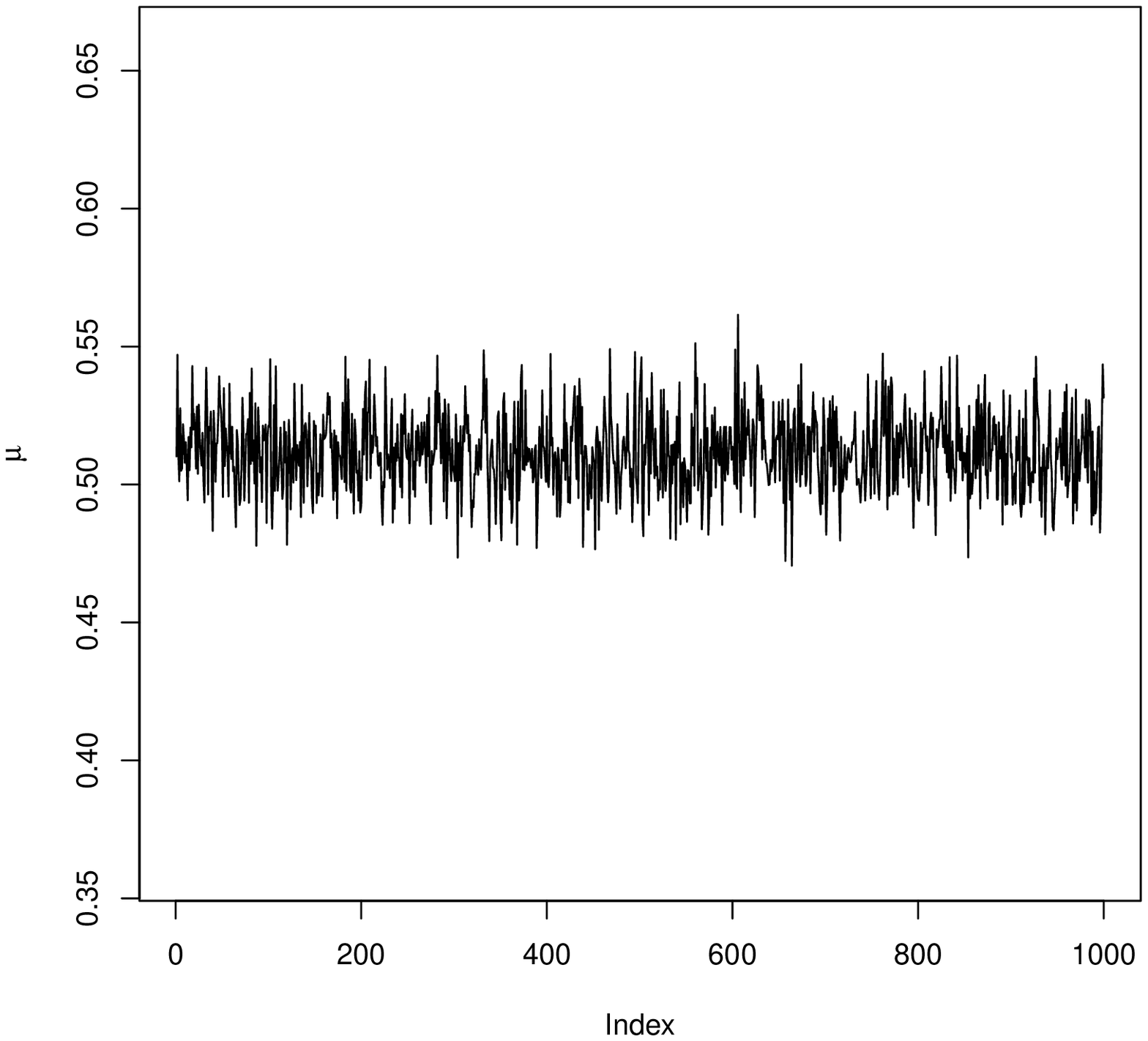}}
    \subfloat[Peak  parameter  (Data set 5).]{\label{figurchain:10}\includegraphics[width=50mm,height=30mm]{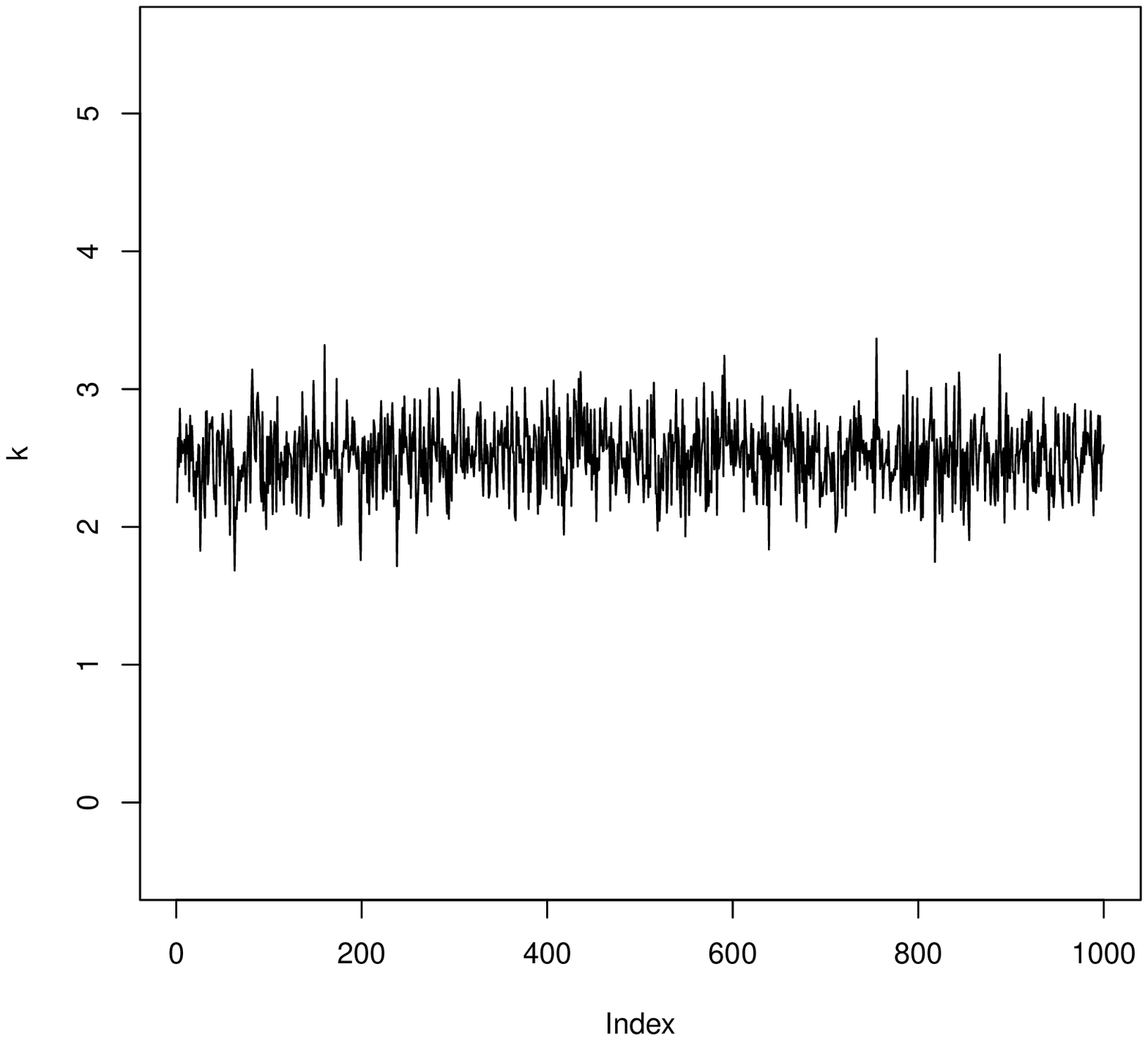}}

\end{figure}

\begin{figure}[htp]
  \centering
  \label{figurpost:0}\caption{Histograms of posterior samples for the  asymmetric and the peak parameters.}

  \subfloat[Asymmetric  parameter (Data set 1). ]{\label{figurpost:1}\includegraphics[width=55mm,height=35mm]{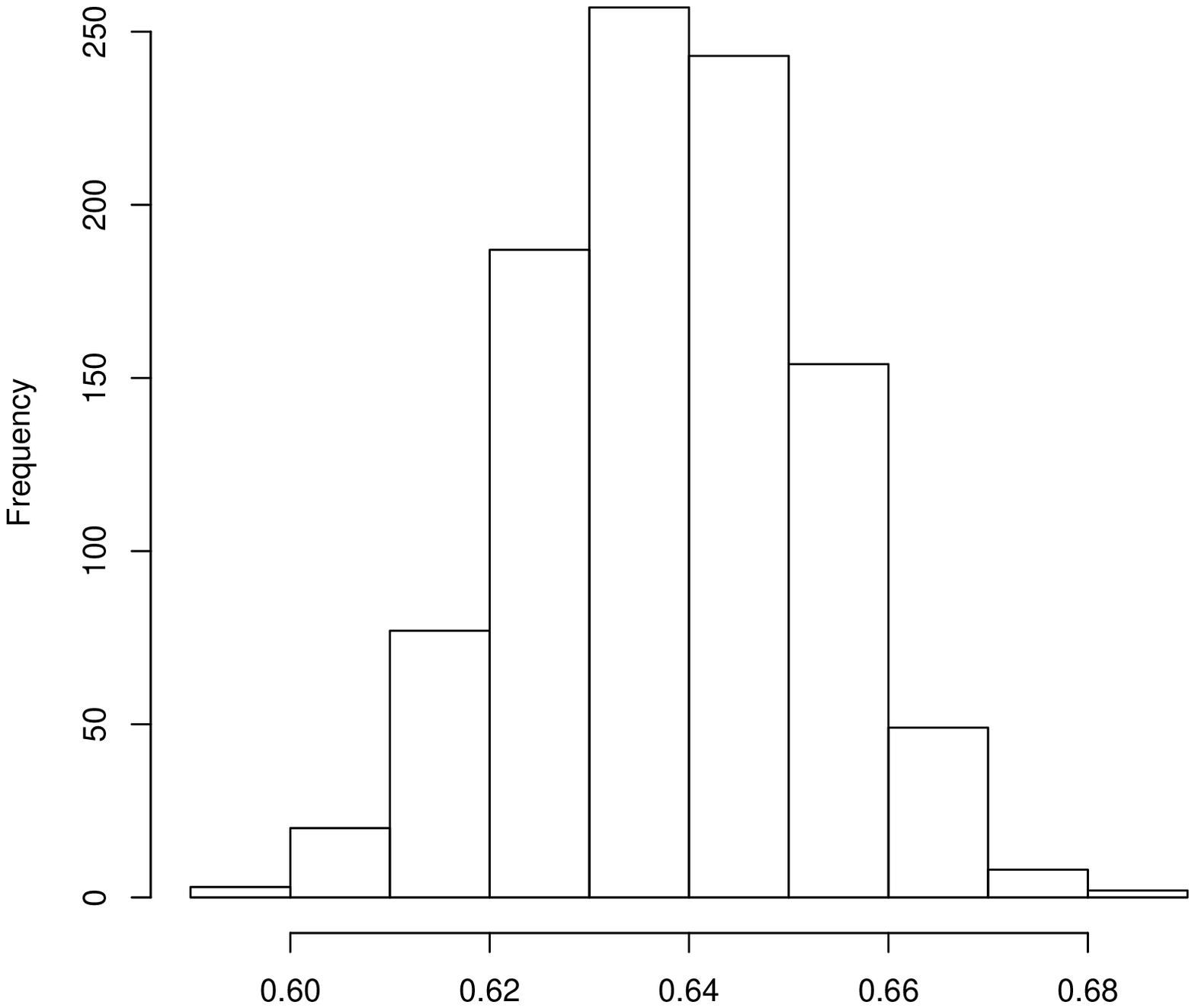}}
    \subfloat[Peak parameter (Data set 1).]{\label{figurpost:2}\includegraphics[width=55mm,height=35mm]{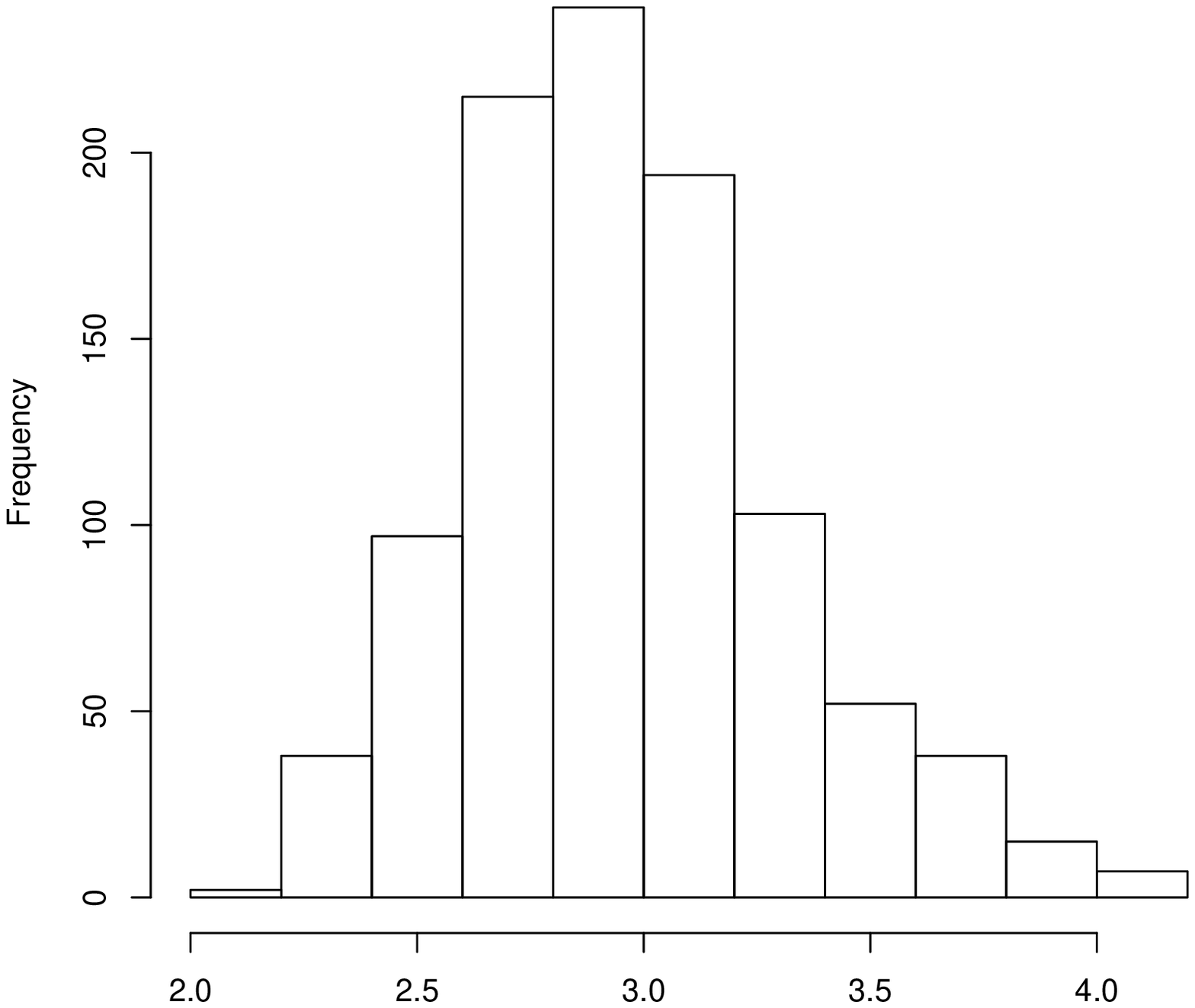}}
  \\
  \subfloat[Asymmetric  parameter (Data set 2). ]{\label{figurpost:3}\includegraphics[width=55mm,height=35mm]{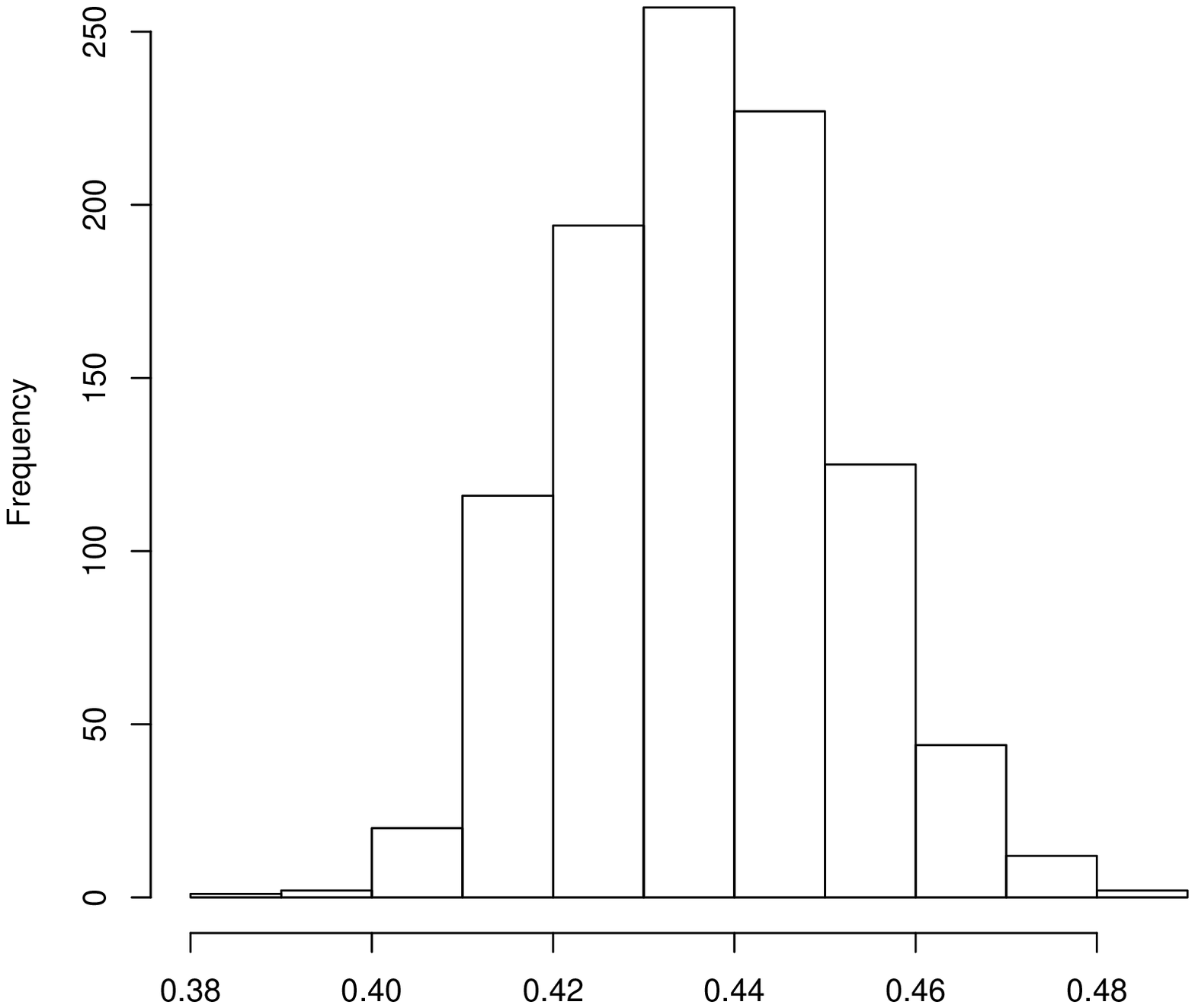}}
    \subfloat[Peak  parameter (Data set 2).]{\label{figurpost:4}\includegraphics[width=55mm,height=35mm]{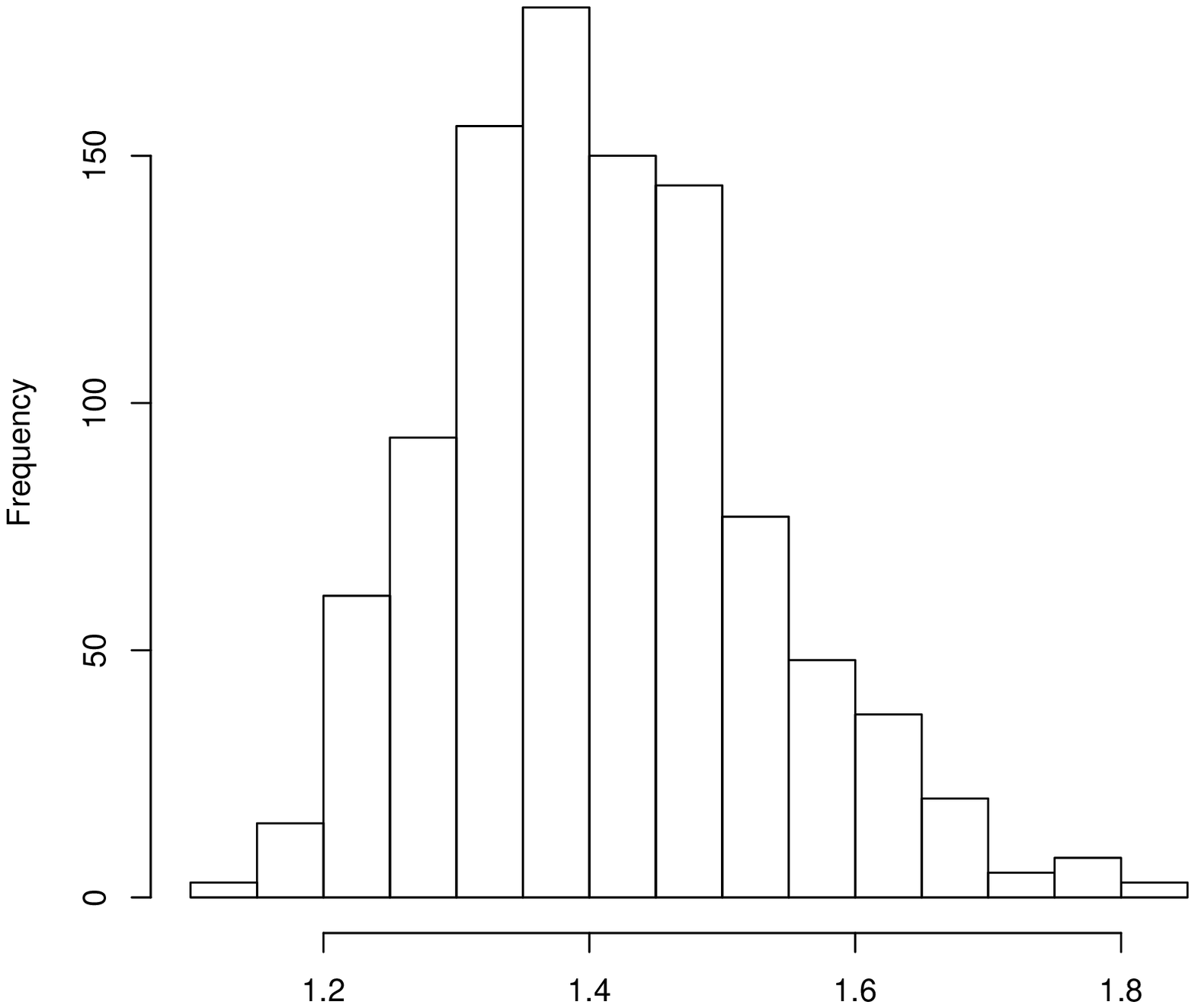}}
  \\
  \subfloat[Asymmetric  parameter (Data set 3). ]{\label{figurpost:5}\includegraphics[width=55mm,height=35mm]{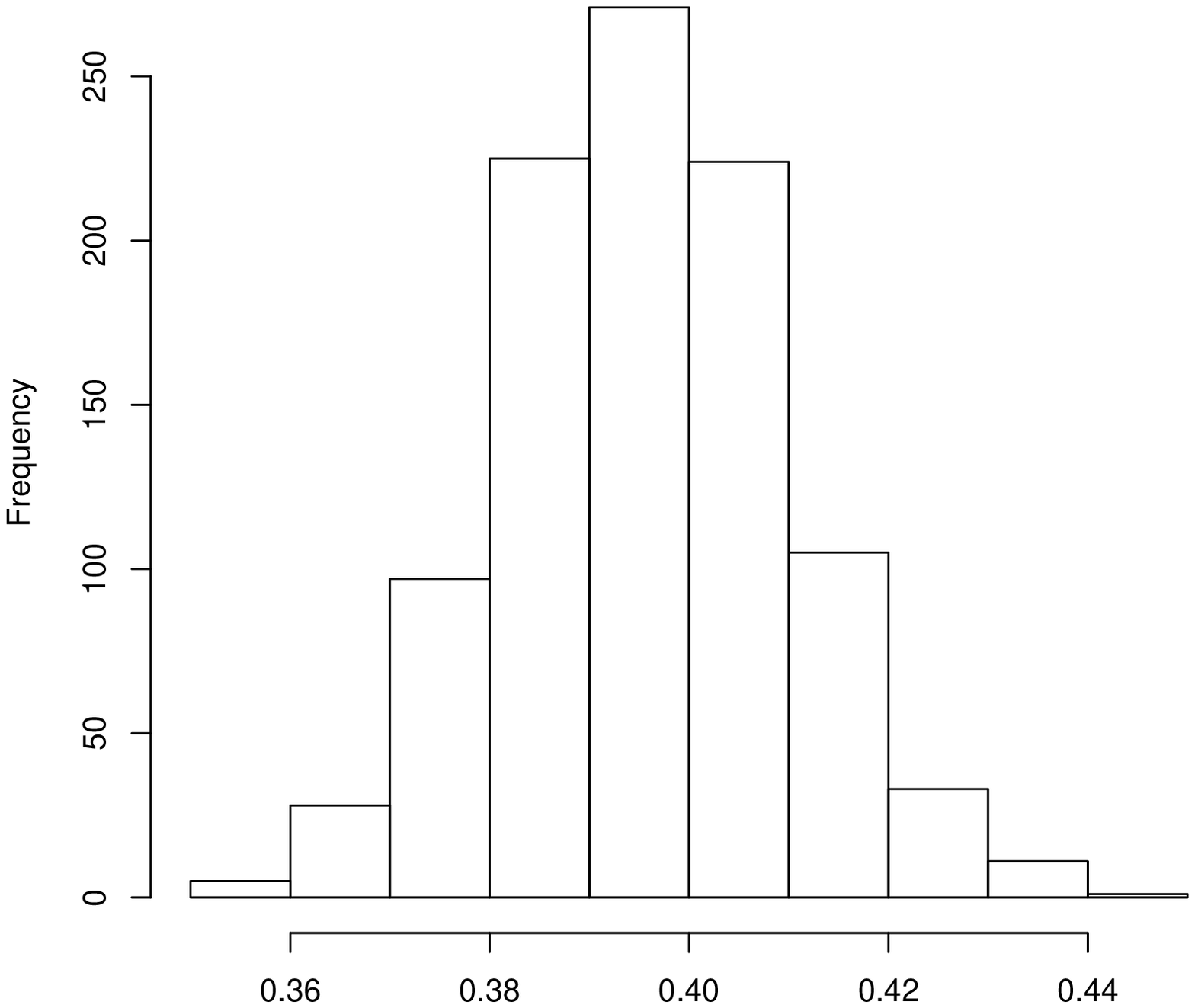}}
    \subfloat[Peak  parameter (Data set 3).]{\label{figurpost:6}\includegraphics[width=55mm,height=35mm]{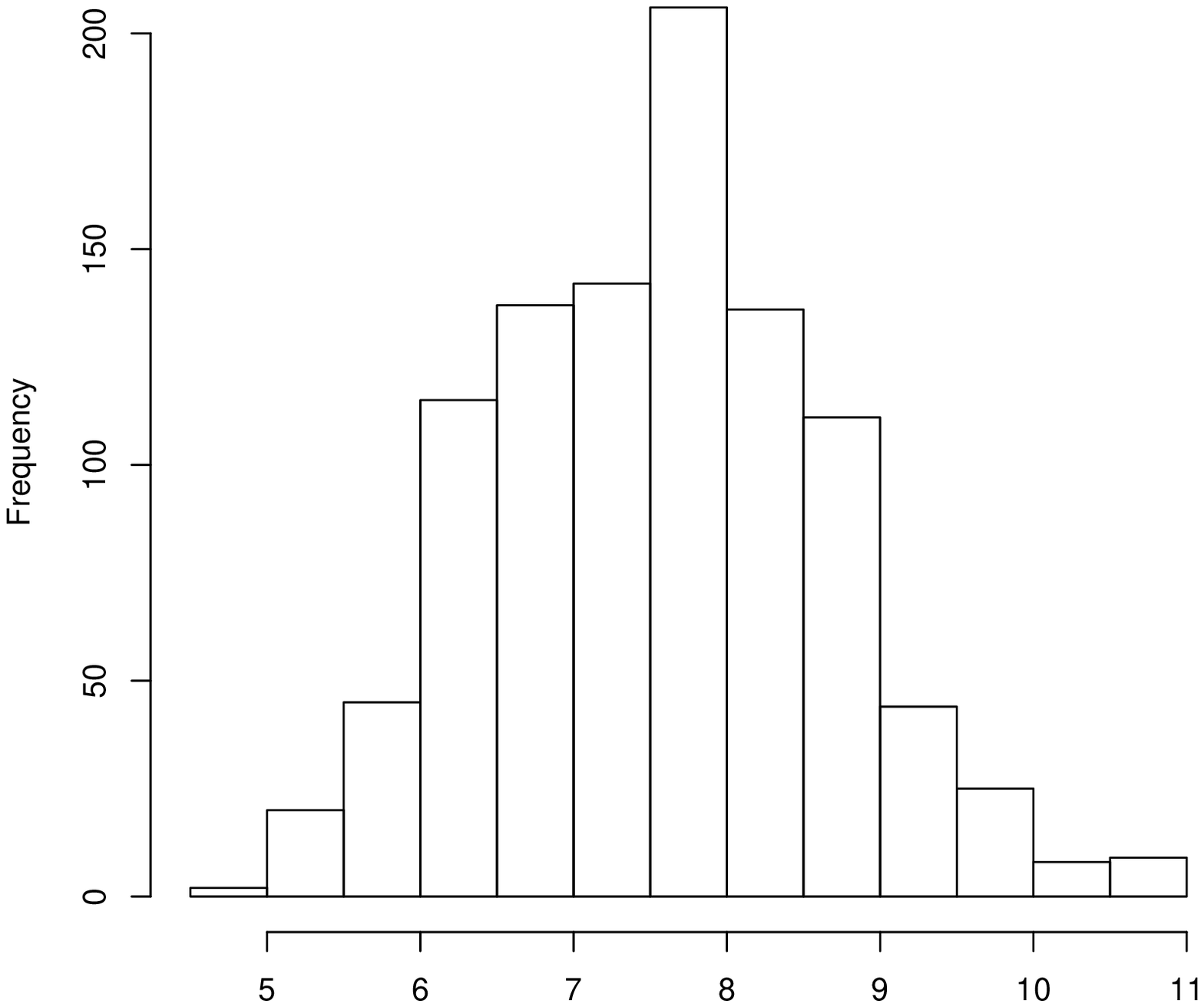}}
  \\
  \subfloat[Asymmetric  parameter (Data set 4). ]{\label{figurpost:7}\includegraphics[width=55mm,height=35mm]{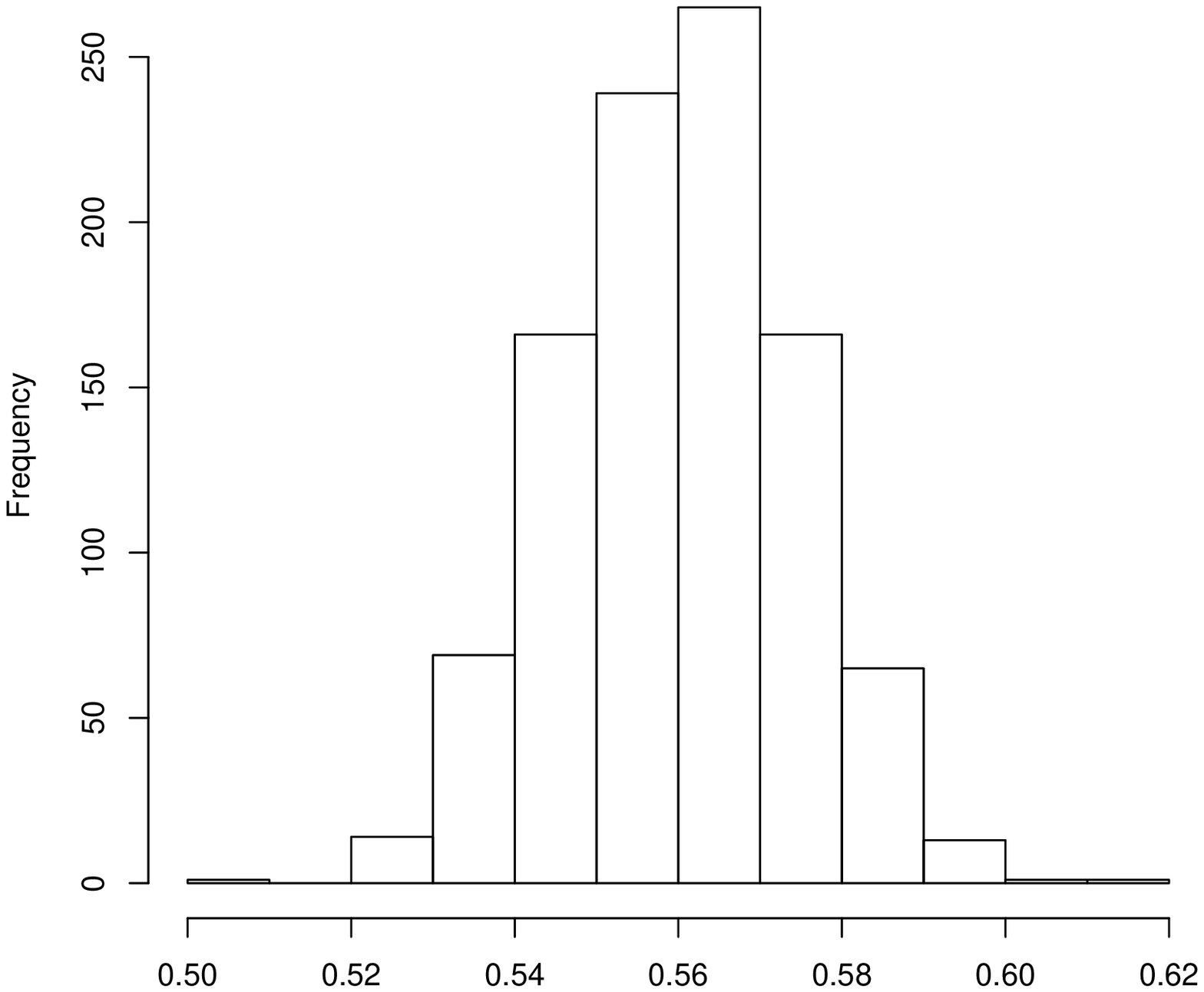}}
    \subfloat[Peak  parameter (Data set 4).]{\label{figurpost:8}\includegraphics[width=55mm,height=35mm]{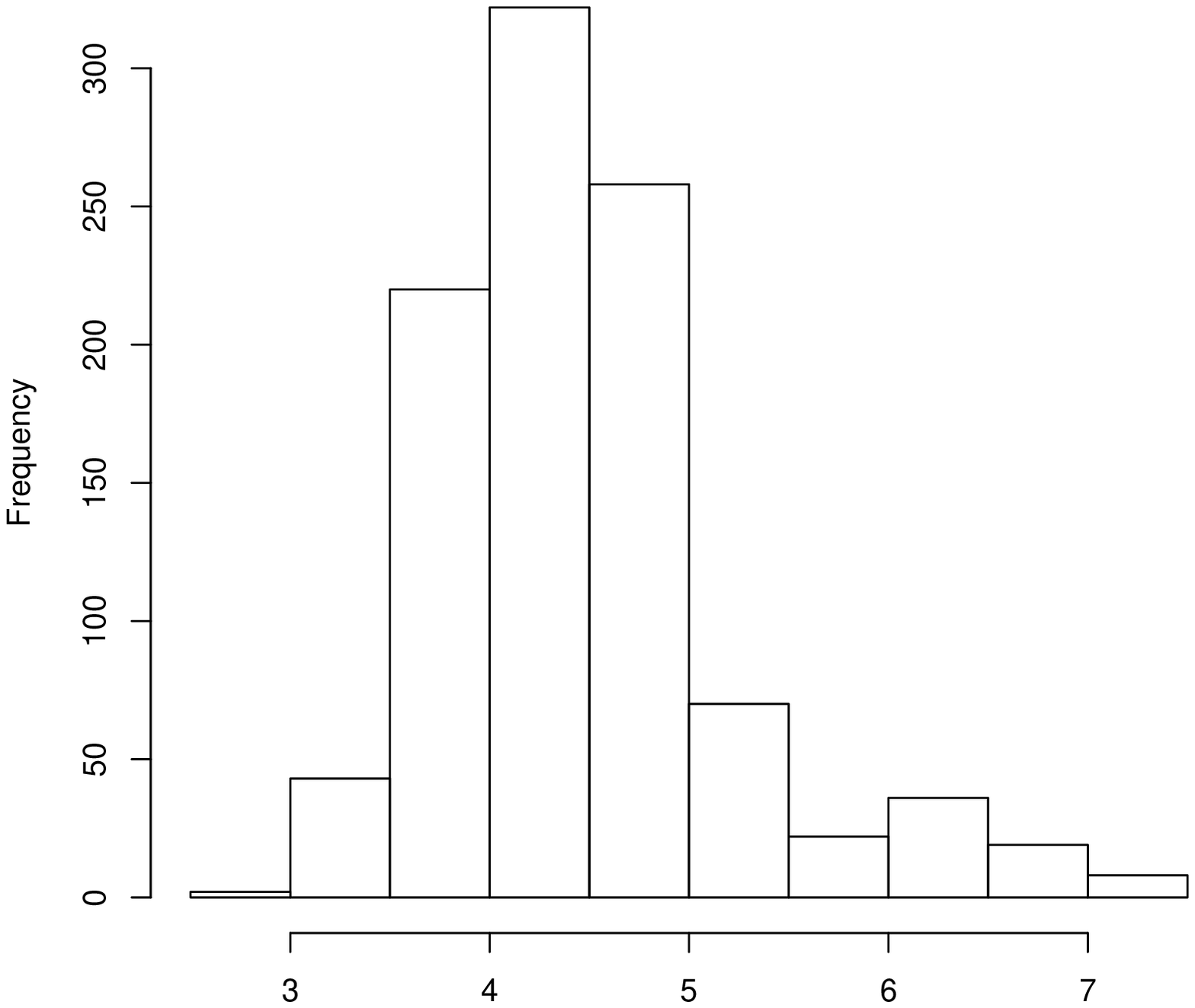}}
    \\
  \subfloat[Asymmetric  parameter (Data set 5). ]{\label{figurpost:9}\includegraphics[width=55mm,height=35mm]{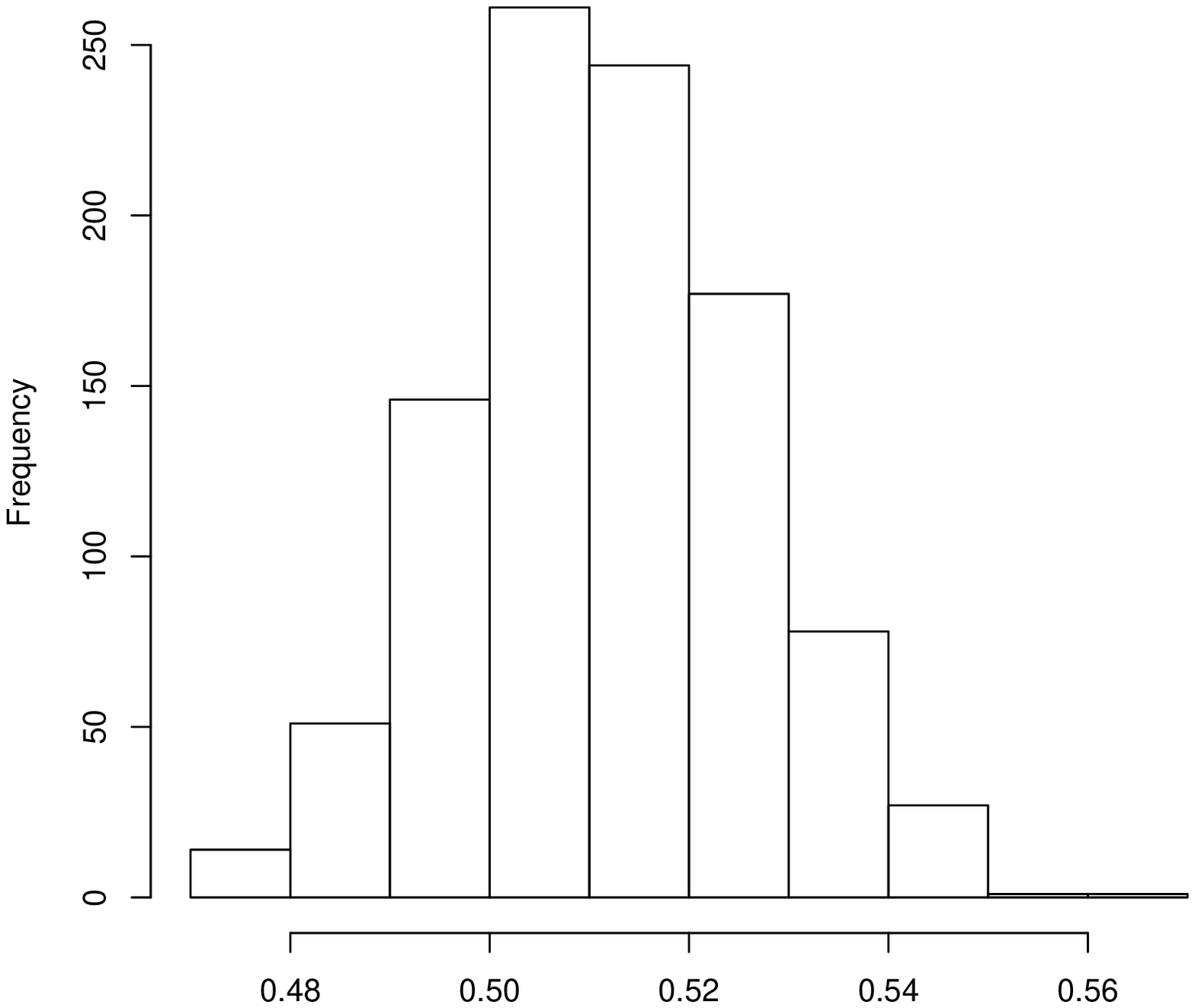}}
    \subfloat[Peak  parameter  (Data set 5).]{\label{figurpost:10}\includegraphics[width=55mm,height=35mm]{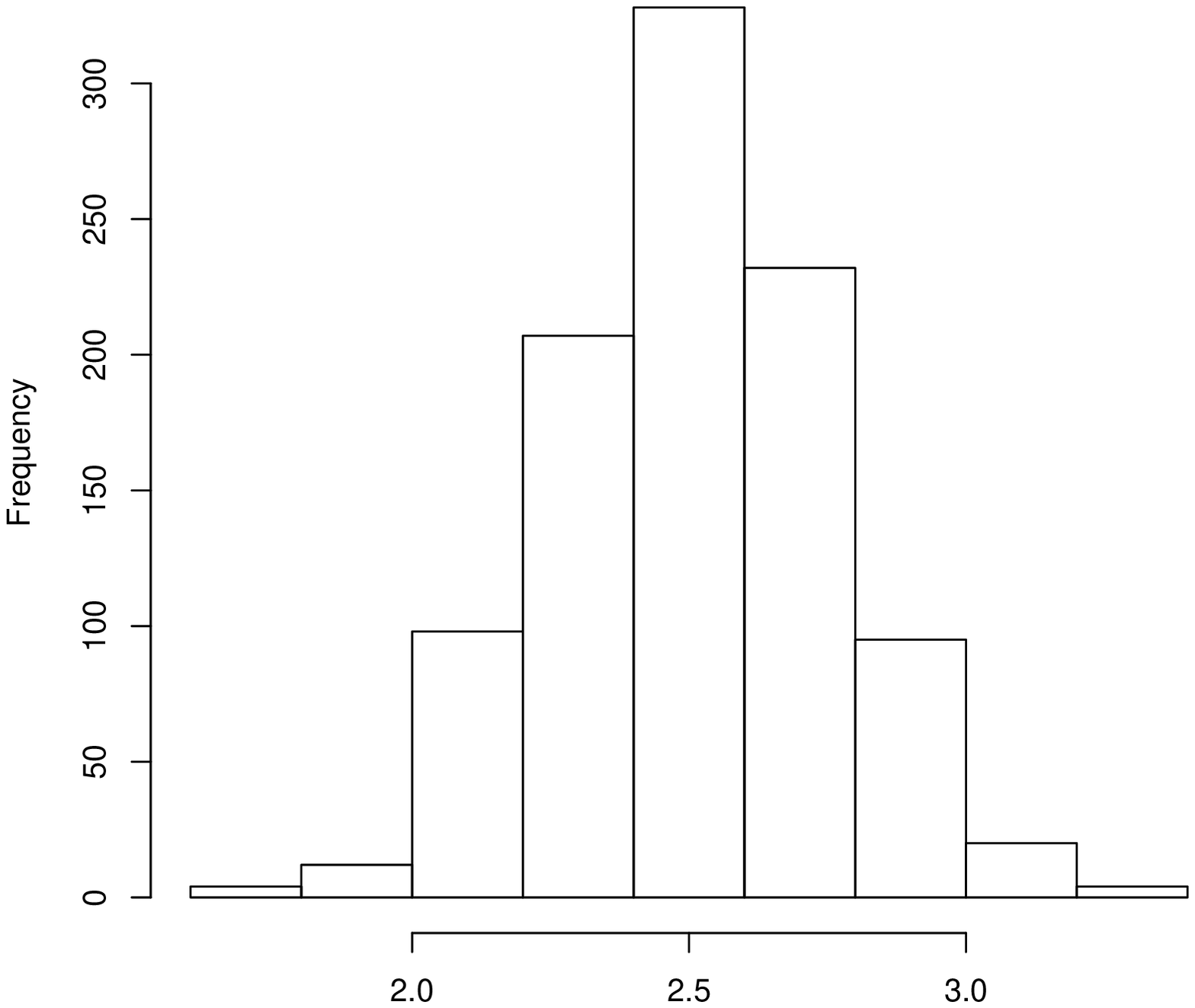}}

\end{figure}

\begin{figure}[htp]
  \centering
  \label{figur:0}\caption{Plots of observed and estimated densities for five data sets. }

  \subfloat[Gene $\sharp 202459$  (Data set 2). ]{\label{figur:1}\includegraphics[width=60mm]{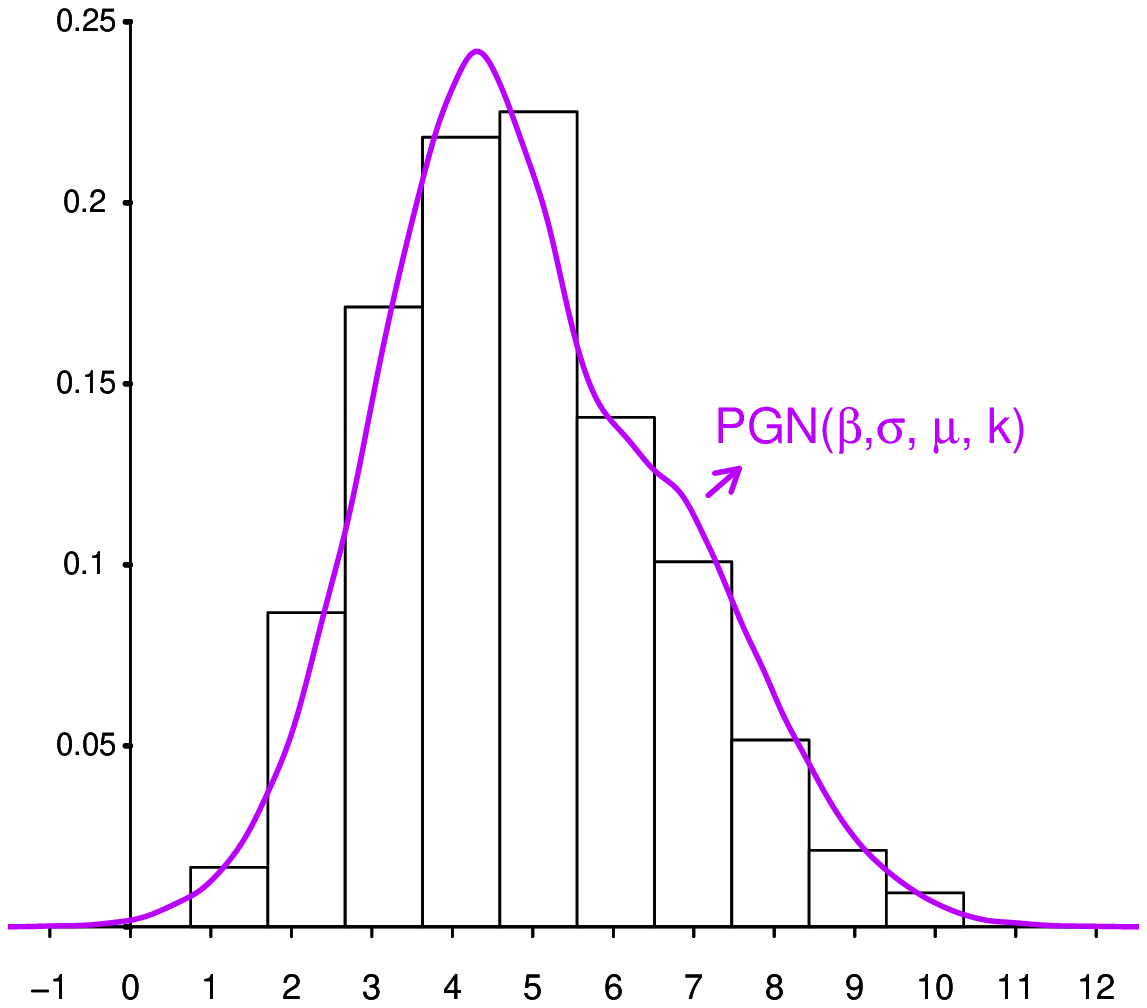}}
    \subfloat[Gene $\sharp 202459$ (Data set 2).]{\label{figur:2}\includegraphics[width=60mm]{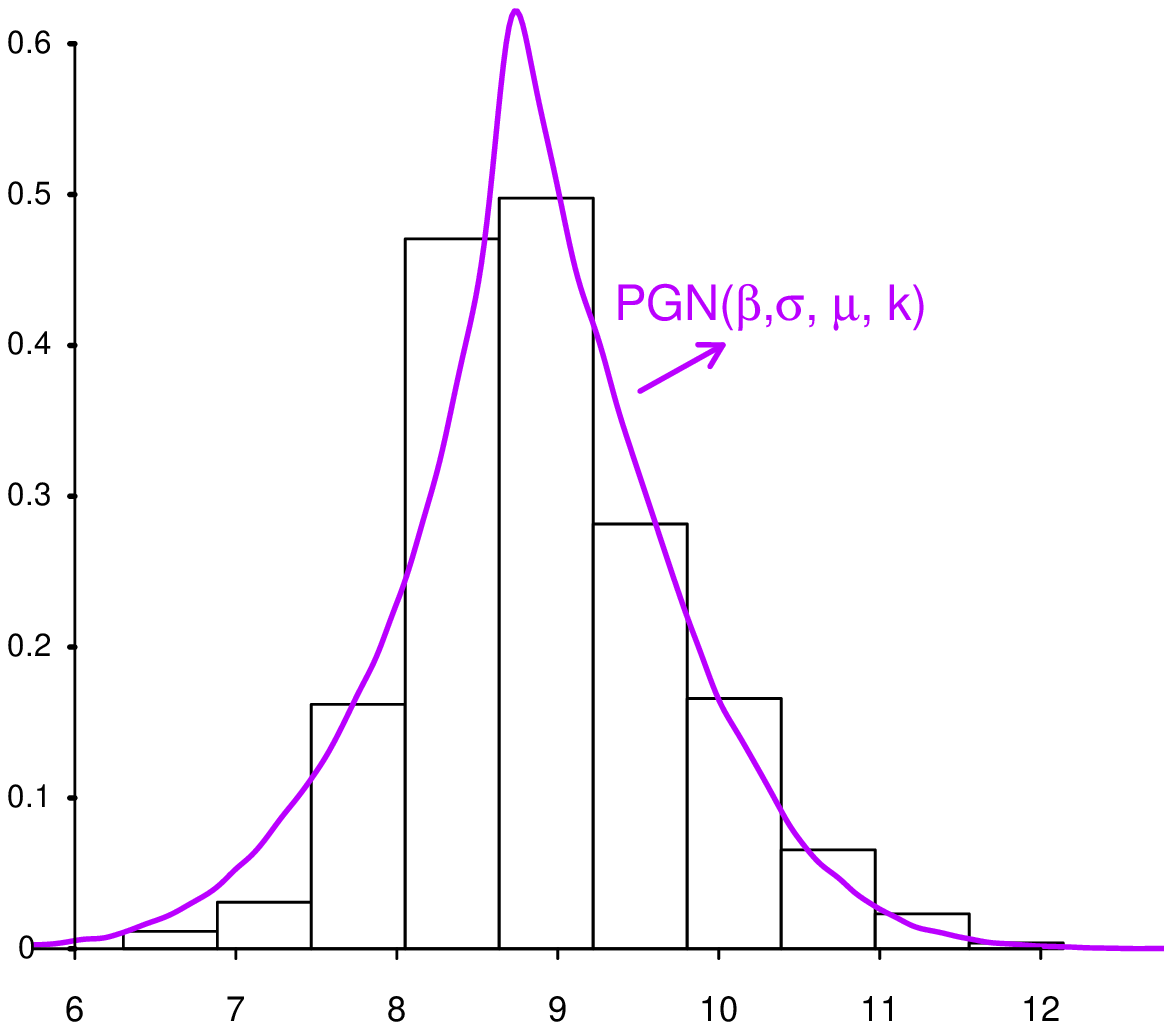}}
  \\
  \subfloat[Gene $\sharp 208288$ (Data set 3).]{\label{figur:3}\includegraphics[width=60mm]{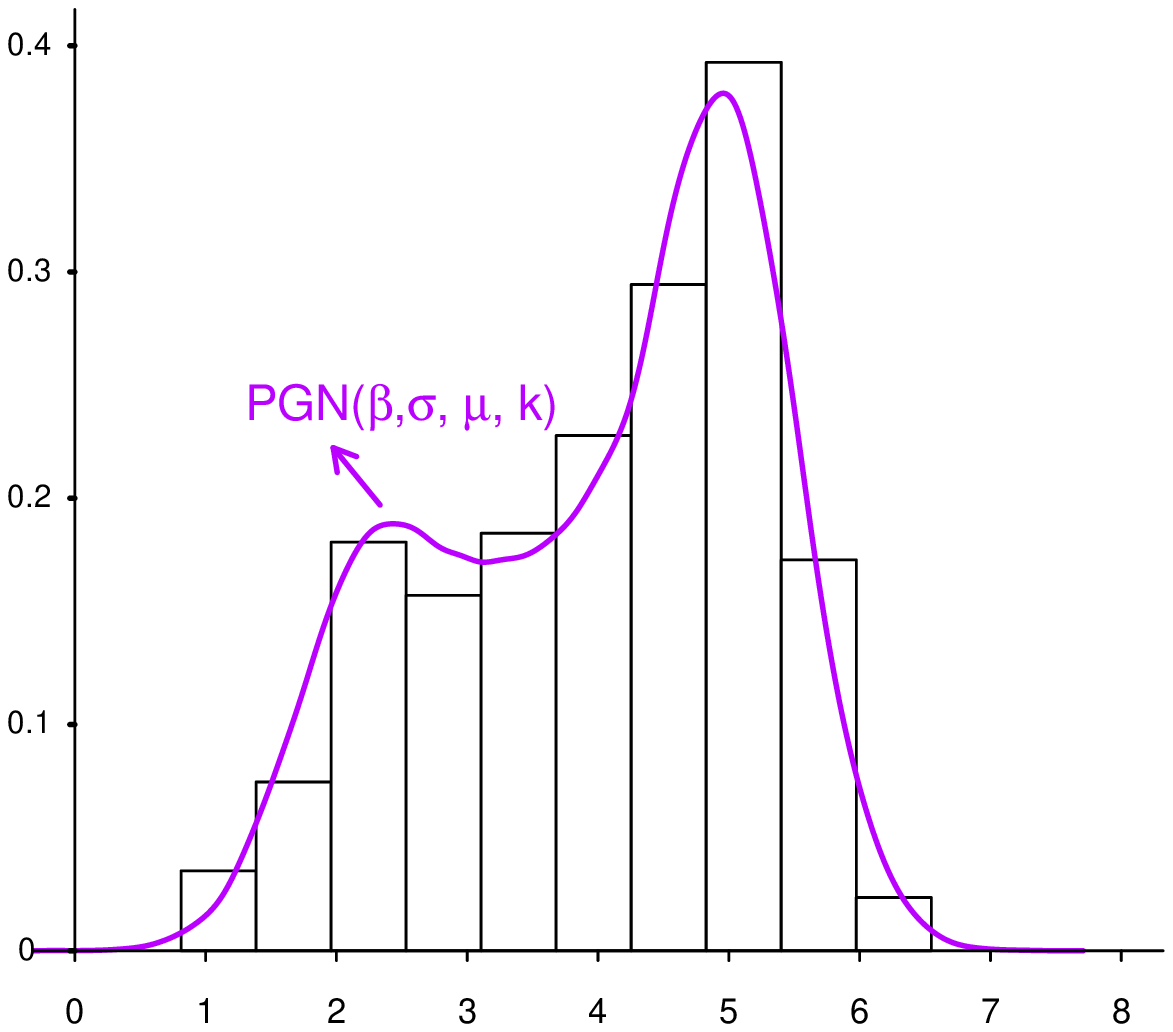}}
    \subfloat[Gene $\sharp 215456$ (Data set 4).]{\label{figur:4}\includegraphics[width=60mm]{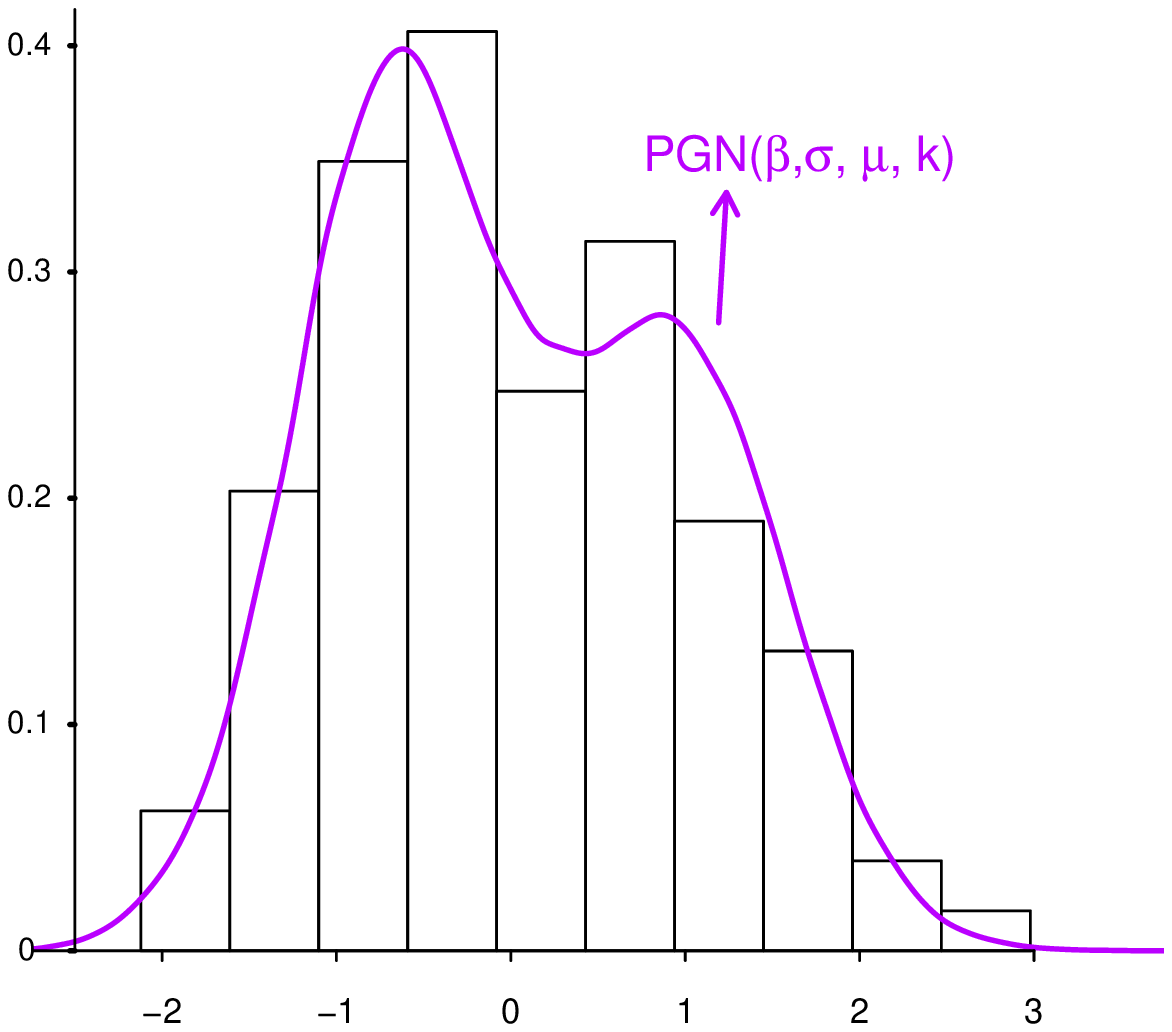}}
  \\
  \subfloat[Gene $\sharp 216437$ (Data set 5).]{\label{figur:5}\includegraphics[width=60mm]{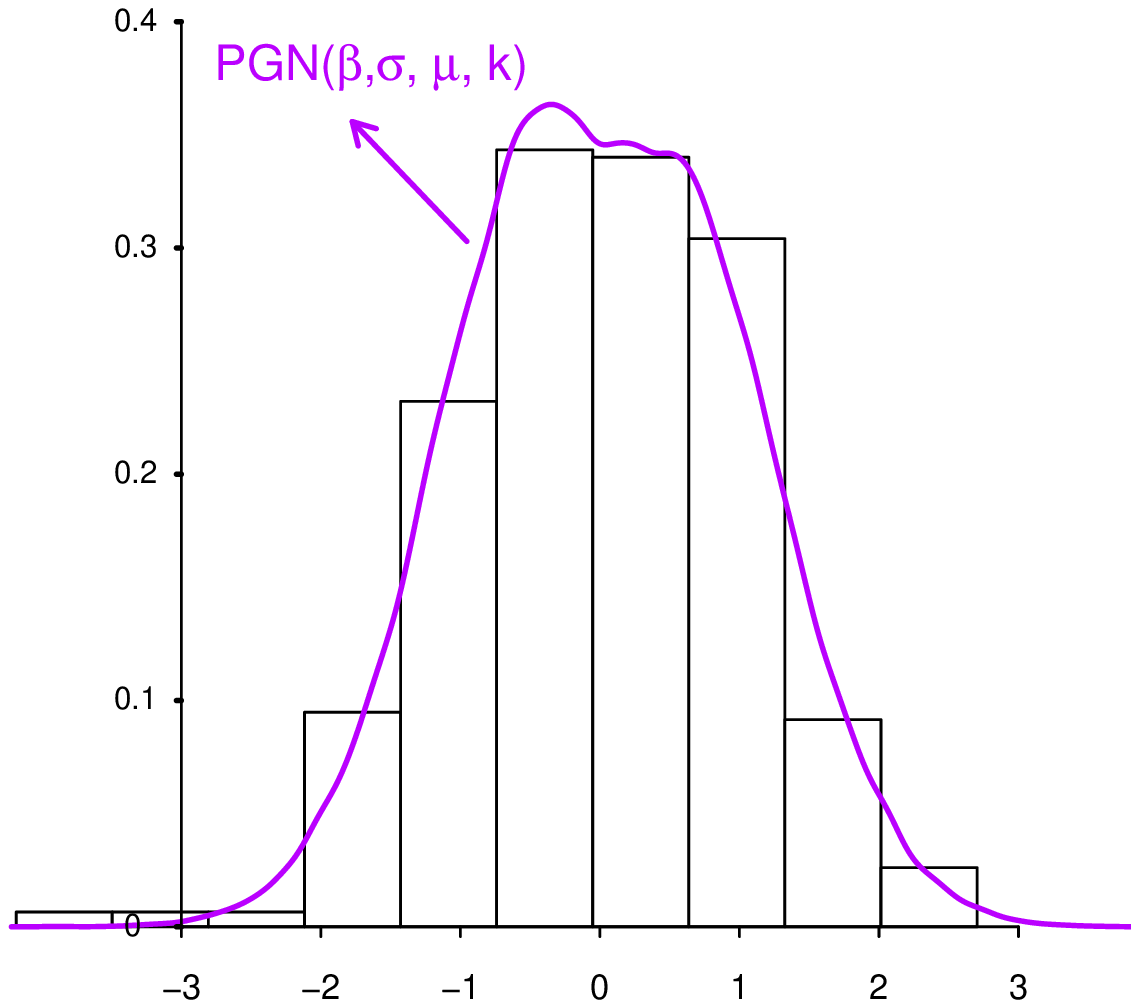}}

\end{figure}

\section{Conclusions}
This article proposed a new class of generalized normal distributions which includes as a special case the normal distribution. Two well known distributions,  i.e.
the Chi distribution and the half normal distribution, are limiting cases of the general construction. The proposed model provides a flexible approach to account
for  asymmetry and biomodality which are two pervasive features of the complex data.
The peak parameter and the asymmetric parameter  control the biamodality and asymmetry of the model, respectively. A Gibbs algorithm was developed for
the Bayesian inference. A simulation study is performed. The simulation results established that
the proposed model has an appropriate ability to estimate
model parameters. Several applications to real data have shown the usefulness of the $PGN$ model for applied statistical research. There is a
bias in the estimation of the parameters when $\mu$ is close $0$ or $1$, so further research could be conducted to obtain bias corrections
for these estimators, thus reducing their systematic errors in the estimation of the parameters.
It is also provide a test for symmetry and a test for normality in the context of Bayesian framework.
The class was in the context of univariate, so further research could be conducted to extend this class to multivariate contexts. It would be interesting if the properties and the flexibility of
this extension  were to be explored and compared in a further work.

%\bibliographystyle{Chicago}

%\bibliography{Bibliography-MM-MC}

\appendix

%\appendixone
\section*{Appendix 1}\label{Itable}

\begin{table}[!htb]
    \caption{Table of $I_m(2\mu,2(1-\mu))$}
    \begin{minipage}{.5\linewidth}
      \caption{$0\leq \mu \leq 0.50$}
      \centering
        \begin{tabular} {| l ||c |c |c |c |}
 \hline
   \  & m = 1 & m = 2 & m = 3 & m = 4 \\ \hline \hline
$\mu$ = 0 & 1 & 1 & 1 & 1 \\ \hline
$\mu$ = 0.01 & 0.9868 & 0.9856 & 0.9811 & 0.9803 \\ \hline
$\mu$ = 0.02 & 0.9733 & 0.9710 & 0.9620 & 0.9606 \\ \hline
$\mu$ = 0.03 & 0.9593 & 0.9563 & 0.9427 & 0.9408 \\ \hline
$\mu$ = 0.04 & 0.9449 & 0.9416 & 0.9231 & 0.9211 \\ \hline
$\mu$ = 0.05 & 0.9302 & 0.9268 & 0.9035 & 0.9015 \\ \hline
$\mu$ = 0.06 & 0.9151 & 0.9120 & 0.8836 & 0.8819 \\ \hline
$\mu$ = 0.07 & 0.8996 & 0.8971 & 0.8636 & 0.8624 \\ \hline
$\mu$ = 0.08 & 0.8837 & 0.8823 & 0.8436 & 0.8430 \\ \hline
$\mu$ = 0.09 & 0.8675 & 0.8676 & 0.8234 & 0.8238 \\ \hline
$\mu$ = 0.10 & 0.8510 & 0.8528 & 0.8031 & 0.8048 \\ \hline
$\mu$ = 0.11 & 0.8341 & 0.8382 & 0.7827 & 0.7860 \\ \hline
$\mu$ = 0.12 & 0.8169 & 0.8237 & 0.7623 & 0.7674 \\ \hline
$\mu$ = 0.13 & 0.7993 & 0.8092 & 0.7418 & 0.7491 \\ \hline
$\mu$ = 0.14 & 0.7815 & 0.7949 & 0.7213 & 0.7310 \\ \hline
$\mu$ = 0.15 & 0.7633 & 0.7808 & 0.7008 & 0.7132 \\ \hline
$\mu$ = 0.16 & 0.7449 & 0.7669 & 0.6802 & 0.6957 \\ \hline
$\mu$ = 0.17 & 0.7262 & 0.7531 & 0.6596 & 0.6785 \\ \hline
$\mu$ = 0.18 & 0.7071 & 0.7396 & 0.6390 & 0.6617 \\ \hline
$\mu$ = 0.19 & 0.6879 & 0.7263 & 0.6185 & 0.6452 \\ \hline
$\mu$ = 0.20 & 0.6683 & 0.7132 & 0.5979 & 0.6291 \\ \hline
$\mu$ = 0.21 & 0.6485 & 0.7004 & 0.5774 & 0.6134 \\ \hline
$\mu$ = 0.22 & 0.6284 & 0.6879 & 0.5569 & 0.5980 \\ \hline
$\mu$ = 0.23 & 0.6081 & 0.6756 & 0.5364 & 0.5831 \\ \hline
$\mu$ = 0.24 & 0.5876 & 0.6637 & 0.5159 & 0.5687 \\ \hline
$\mu$ = 0.25 & 0.5668 & 0.6521 & 0.4955 & 0.5547 \\ \hline
$\mu$ = 0.26 & 0.5459 & 0.6409 & 0.4752 & 0.5411 \\ \hline
$\mu$ = 0.27 & 0.5247 & 0.6300 & 0.4549 & 0.5280 \\ \hline
$\mu$ = 0.28 & 0.5033 & 0.6194 & 0.4346 & 0.5154 \\ \hline
$\mu$ = 0.29 & 0.4817 & 0.6093 & 0.4144 & 0.5033 \\ \hline
$\mu$ = 0.30 & 0.4600 & 0.5995 & 0.3942 & 0.4916 \\ \hline
$\mu$ = 0.31 & 0.4380 & 0.5901 & 0.3741 & 0.4805 \\ \hline
$\mu$ = 0.32 & 0.4160 & 0.5812 & 0.3540 & 0.4699 \\ \hline
$\mu$ = 0.33 & 0.3937 & 0.5727 & 0.3340 & 0.4598 \\ \hline
$\mu$ = 0.34 & 0.3713 & 0.5646 & 0.3140 & 0.4503 \\ \hline
$\mu$ = 0.35 & 0.3488 & 0.5569 & 0.2941 & 0.4413 \\ \hline
$\mu$ = 0.36 & 0.3261 & 0.5497 & 0.2743 & 0.4329 \\ \hline
$\mu$ = 0.37 & 0.3033 & 0.5430 & 0.2545 & 0.4250 \\ \hline
$\mu$ = 0.38 & 0.2804 & 0.5367 & 0.2347 & 0.4176 \\ \hline
$\mu$ = 0.39 & 0.2574 & 0.5309 & 0.2150 & 0.4109 \\ \hline
$\mu$ = 0.40 & 0.2343 & 0.5256 & 0.1953 & 0.4047 \\ \hline
$\mu$ = 0.41 & 0.2112 & 0.5208 & 0.1757 & 0.3991 \\ \hline
$\mu$ = 0.42 & 0.1879 & 0.5164 & 0.1560 & 0.3940 \\ \hline
$\mu$ = 0.43 & 0.1646 & 0.5126 & 0.1365 & 0.3896 \\ \hline
$\mu$ = 0.44 & 0.1412 & 0.5093 & 0.1169 & 0.3857 \\ \hline
$\mu$ = 0.45 & 0.1177 & 0.5064 & 0.0974 & 0.3825 \\ \hline
$\mu$ = 0.46 & 0.0942 & 0.5041 & 0.0779 & 0.3798 \\ \hline
$\mu$ = 0.47 & 0.0707 & 0.5023 & 0.0584 & 0.3777 \\ \hline
$\mu$ = 0.48 & 0.0471 & 0.5010 & 0.0389 & 0.3762 \\ \hline
$\mu$ = 0.49 & 0.0236 & 0.5003 & 0.0195 & 0.3753 \\ \hline
$\mu$ = 0.50 & 0.0000 & 0.5000 & 0.0000 & 0.3750 \\ \hline
\end{tabular}
    \end{minipage}%
    \begin{minipage}{.5\linewidth}
      \centering
        \caption{$0.51 \leq \mu \leq 1$}
        \begin{tabular} {| l ||c |c |c |c |}
 \hline
   \  & m = 1 & m = 2 & m = 3 & m = 4 \\ \hline \hline
$\mu$ = 0.51 & -0.0236 & 0.5003 & -0.0195 & 0.3753 \\ \hline
$\mu$ = 0.52 & -0.0471 & 0.5010 & -0.0389 & 0.3762 \\ \hline
$\mu$ = 0.53 & -0.0707 & 0.5023 & -0.0584 & 0.3777 \\ \hline
$\mu$ = 0.54 & -0.0942 & 0.5041 & -0.0779 & 0.3798 \\ \hline
$\mu$ = 0.55 & -0.1177 & 0.5064 & -0.0974 & 0.3825 \\ \hline
$\mu$ = 0.56 & -0.1412 & 0.5093 & -0.1169 & 0.3857 \\ \hline
$\mu$ = 0.57 & -0.1646 & 0.5126 & -0.1365 & 0.3896 \\ \hline
$\mu$ = 0.58 & -0.1879 & 0.5164 & -0.1560 & 0.3940 \\ \hline
$\mu$ = 0.59 & -0.2112 & 0.5208 & -0.1757 & 0.3991 \\ \hline
$\mu$ = 0.60 & -0.2343 & 0.5256 & -0.1953 & 0.4047 \\ \hline
$\mu$ = 0.61 & -0.2574 & 0.5309 & -0.2150 & 0.4109 \\ \hline
$\mu$ = 0.62 & -0.2804 & 0.5367 & -0.2347 & 0.4176 \\ \hline
$\mu$ = 0.63 & -0.3033 & 0.5430 & -0.2545 & 0.4250 \\ \hline
$\mu$ = 0.64 & -0.3261 & 0.5497 & -0.2743 & 0.4329 \\ \hline
$\mu$ = 0.65 & -0.3488 & 0.5569 & -0.2941 & 0.4413 \\ \hline
$\mu$ = 0.66 & -0.3713 & 0.5646 & -0.3140 & 0.4503 \\ \hline
$\mu$ = 0.67 & -0.3937 & 0.5727 & -0.3340 & 0.4598 \\ \hline
$\mu$ = 0.68 & -0.4160 & 0.5812 & -0.3540 & 0.4699 \\ \hline
$\mu$ = 0.69 & -0.4380 & 0.5901 & -0.3741 & 0.4805 \\ \hline
$\mu$ = 0.70 & -0.4600 & 0.5995 & -0.3942 & 0.4916 \\ \hline
$\mu$ = 0.71 & -0.4817 & 0.6093 & -0.4144 & 0.5033 \\ \hline
$\mu$ = 0.72 & -0.5033 & 0.6194 & -0.4346 & 0.5154 \\ \hline
$\mu$ = 0.73 & -0.5247 & 0.6300 & -0.4549 & 0.5280 \\ \hline
$\mu$ = 0.74 & -0.5459 & 0.6409 & -0.4752 & 0.5411 \\ \hline
$\mu$ = 0.75 & -0.5668 & 0.6521 & -0.4955 & 0.5547 \\ \hline
$\mu$ = 0.76 & -0.5876 & 0.6637 & -0.5159 & 0.5687 \\ \hline
$\mu$ = 0.77 & -0.6081 & 0.6756 & -0.5364 & 0.5831 \\ \hline
$\mu$ = 0.78 & -0.6284 & 0.6879 & -0.5569 & 0.5980 \\ \hline
$\mu$ = 0.79 & -0.6485 & 0.7004 & -0.5774 & 0.6134 \\ \hline
$\mu$ = 0.80 & -0.6683 & 0.7132 & -0.5979 & 0.6291 \\ \hline
$\mu$ = 0.81 & -0.6879 & 0.7263 & -0.6185 & 0.6452 \\ \hline
$\mu$ = 0.82 & -0.7071 & 0.7396 & -0.6390 & 0.6617 \\ \hline
$\mu$ = 0.83 & -0.7262 & 0.7531 & -0.6596 & 0.6785 \\ \hline
$\mu$ = 0.84 & -0.7449 & 0.7669 & -0.6802 & 0.6957 \\ \hline
$\mu$ = 0.85 & -0.7633 & 0.7808 & -0.7008 & 0.7132 \\ \hline
$\mu$ = 0.86 & -0.7815 & 0.7949 & -0.7213 & 0.7310 \\ \hline
$\mu$ = 0.87 & -0.7993 & 0.8092 & -0.7418 & 0.7491 \\ \hline
$\mu$ = 0.88 & -0.8169 & 0.8237 & -0.7623 & 0.7674 \\ \hline
$\mu$ = 0.89 & -0.8341 & 0.8382 & -0.7827 & 0.7860 \\ \hline
$\mu$ = 0.90 & -0.8510 & 0.8528 & -0.8031 & 0.8048 \\ \hline
$\mu$ = 0.91 & -0.8675 & 0.8676 & -0.8234 & 0.8238 \\ \hline
$\mu$ = 0.92 & -0.8837 & 0.8823 & -0.8436 & 0.8430 \\ \hline
$\mu$ = 0.93 & -0.8996 & 0.8971 & -0.8636 & 0.8624 \\ \hline
$\mu$ = 0.94 & -0.9151 & 0.9120 & -0.8836 & 0.8819 \\ \hline
$\mu$ = 0.95 & -0.9302 & 0.9268 & -0.9035 & 0.9015 \\ \hline
$\mu$ = 0.96 & -0.9449 & 0.9416 & -0.9231 & 0.9211 \\ \hline
$\mu$ = 0.97 & -0.9593 & 0.9563 & -0.9427 & 0.9408 \\ \hline
$\mu$ = 0.98 & -0.9733 & 0.9710 & -0.9620 & 0.9606 \\ \hline
$\mu$ = 0.99 & -0.9868 & 0.9856 & -0.9811 & 0.9803 \\ \hline
$\mu$ = 1 & -1 & 1 & -1 & 1 \\ \hline
\end{tabular}
    \end{minipage}
\end{table}

%\bibliographystyle{Chicago}

%\bibliography{Bibliography-MM-MC}

\end{document}